\documentclass[dvips]{amsart}	           
	
\usepackage[dvips]{graphicx}

\usepackage{epstopdf}
\usepackage{psfrag}
\usepackage{amsmath,amsthm,amscd,color}
\usepackage[small,nohug,PostScript=dvips]{diagrams}
\usepackage{hyperref}
\diagramstyle[labelstyle=\scriptstyle]

\newtheorem{thm}{Theorem}[section]
\newtheorem{cor}[thm]{Corollary}
\newtheorem{lem}[thm]{Lemma}
\newtheorem{prop}[thm]{Proposition}
\newtheorem*{result}{Main Result}

\theoremstyle{definition}
\newtheorem{defn}[thm]{Definition}

\theoremstyle{remark}
\newtheorem{rem}[thm]{Remark}
\numberwithin{equation}{section}

\newcommand{\R}{\mathbb R}
\newcommand{\Z}{\mathbb Z}

\renewcommand{\P}{\mathcal{P}}
\newcommand{\U}{\mathcal{U}}
\newcommand{\LL}{\mathcal{L}}
\newcommand{\ra}{\rightarrow}

\newcommand{\sz}{SL(2,\Z)}

\newcommand{\Id}{\rm Id}

\newcommand{\OO}{{\mathcal O}}
\newcommand{\E}{{\mathcal E}}
\newcommand{\T}{{\mathcal T}}
\newcommand{\af}{{\mathfrak{a}}}
\newcommand{\bb}{{\mathfrak{b}}}

\begin{document}

%\title{Structure of attractors for $(a,b)$-continued fractions}
\title [Attractors of $(a,b)$-continued fractions]{Structure of attractors for $(a,b)$-continued fraction transformations}
\author{Svetlana Katok}
\address{Department of Mathematics, The Pennsylvania State University,
University Park, PA 16802} \email{katok\_s@math.psu.edu}
\author{Ilie Ugarcovici}
\address{Department of Mathematical Sciences,
DePaul University, Chicago, IL 60614} \email{iugarcov@depaul.edu}
\thanks{We are grateful to Don Zagier for helpful discussions and the Max Plank Institute for Mathematics in Bonn for its hospitality and support. The second author is partially supported by the NSF grant DMS-0703421}

\date{\today}

%----------------------------------------------------------------
\begin{abstract} We study a two-parameter family of one-dimensional maps and related $(a,b)$-continued fractions suggested for consideration by Don Zagier. We prove
 that the associated natural extension maps
have attractors with  finite rectangular structure for the entire parameter set except for a Cantor-like set of one-dimensional Lebesgue zero measure that we completely describe.
We show that the structure of these attractors can be ``computed" from the data $(a,b)$, and that for a dense open set of parameters the Reduction theory conjecture holds, i.e. every point is mapped to the attractor after finitely many iterations. We also show how this theory can be applied to the study of invariant measures and ergodic properties of the associated Gauss-like maps.
\end{abstract}
\maketitle
\tableofcontents

\section{Introduction}\label{s:1}
The standard generators $T(x)=x+1$, $S(x)=-1/x$ of the modular group  $\sz$ were used classically to define piecewise continuous maps acting on the extended real line $\bar \R=\R\cup \{\infty\}$  that led to well-known continued fraction algorithms. 
In this paper we present a general method of constructing  such maps suggested by Don Zagier, and study their dynamical properties  and associated generalized continued fraction transformations. 

Let $\mathcal P$ be the two-dimensional parameter set
\[
\mathcal P=\{(a,b)\in \R^2\, |\, a\leq 0\leq b,\,b-a\geq 1,\,-ab\leq 1\}
%   \supset
%\Delta=\{(a,b)\,|\,-1\leq a\leq 0\leq b\leq 1, b-a\geq 1\},
\]
and consider the map $f_{a,b}:\bar\R\ra\bar\R$ defined as
\begin{equation}\label{fab}
f_{a,b}(x)=\begin{cases}
x+1  &\text{ if }  x< a\\
-\displaystyle\frac{1}{x} &\text{ if } a\le x<b\\
x-1  &\text{ if } x\ge b\,.
\end{cases}
\end{equation}
Using the first return map of $f_{a,b}$ to the interval $[a,b)$, denoted by $\hat f_{a,b}$, we introduce a two-dimensional family of %foot{SK13: change}
%various 
continued fraction algorithms and study their properties. We mention here three classical examples: the case $a=-1/2$, $b=1/2$ gives the ``nearest-integer" continued fractions considered first by Hurwitz in \cite{Hurwitz1}, the case $a=-1$, $b=0$  described in \cite{Z, K3} gives the ``minus" (backward) continued fractions, while the situation $a=-1$, $b=1$ was presented in \cite{S1, KU2} in connection with a method of symbolically coding the geodesic flow on the modular surface following Artin's pioneering work \cite{Artin}. Also, in the case $b-a=1$, the class of one-parameter maps $f_{b-1,b}$ with $b\in [0,1]$ is conceptually similar to the ``$\alpha$-transformations" introduced by Nakada in \cite{N1} and studied subsequently in \cite{LM, MCM, NN1,NN2, Sw}. %\foot{IU: added references}
%Also, in the case $b-a=1$, the class of one-parameter maps $f_{b-1,b}$, $b\in [0,1]$ is very similar (at least in spirit) to the ``$\alpha$-transformations" introduced by Nakada in \cite{N1} and studied by several other authors.

 The main object of our study is a two-dimensional realization of the natural extension map of $f_{a,b}$, $F_{a,b}:\bar\R^2\setminus\Delta\ra \bar\R^2\setminus\Delta$,  $\Delta=\{(x,y)\in \bar\R^2| x=y\}$, %\foot{IU: changed}
 %(also called the \emph{reduction map}) 
 defined by 
  \begin{equation}\label{Fab}
F_{a,b}(x,y)=\begin{cases}
(x+1,y+1)  &\text{ if } y<a\\
\left(-\displaystyle\frac{1}{x},-\displaystyle\frac{1}{y}\right) &\text{ if } a\le y<b\\
(x-1,y-1)  &\text{ if } y\ge b\,.
\end{cases}
\end{equation}
The map $F_{a,b}$ is also called the \emph{reduction map}.
%\foot{SK9: change!}
Numerical experiments led Don Zagier to conjecture that such a map $F_{a,b}$ has several interesting properties for all parameter pairs $(a,b)\in \mathcal P$ that we list under the
{\bf Reduction theory conjecture.}
 \begin{enumerate}
%For example, 
\item The map $F_{a,b}$ possesses a global attractor set $D_{a,b}=\cap_{n=0}^\infty F^n(\bar\R^2\setminus\Delta)$
%\subset \bar\R^2$
%=\cap_{n=0}^\infty F^n(\bar\R^2)$ 
on which $F_{a,b}$ is essentially bijective. 
\item The set $D_{a,b}$ consists of two (or one, in degenerate cases) connected components each having {\em finite rectangular structure}, i.e. bounded by non-decreasing step-functions with a finite 
number of steps.   
\item Every point $(x,y)$ of the plane ($x\ne y$) is mapped to $D_{a,b}$ after finitely many iterations of $F_{a,b}$.
\end{enumerate}
 \begin{figure}[htb]
 \psfrag{a}[l]{\small $a$}
\psfrag{b}[l]{\small $b$}
\includegraphics[scale=1.1]{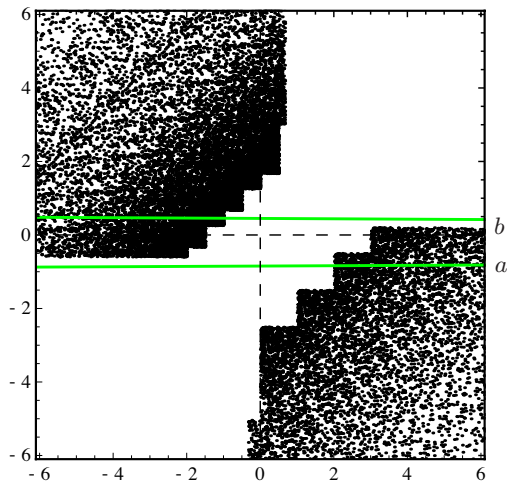}
\caption{Attracting domain for Zagier's example: $a=-\frac{4}{5},\,b=\frac{2}{5}$}
\label{don-a}
\end{figure}
%\foot{SK13: for some reason there is too much space after the picture and section 2 starts on a new page - need to do something}
Figure \ref {don-a} shows the  computer picture of such a the set $D_{a,b}$ with $a=-4/5$, $b=2/5$.  
It is worth mentioning that the complexity of the domain $D_{a,b}$
%of the map $F_{a,b}$ 
increases as $(a,b)$ approach the line segment $b-a=1$ in $\mathcal P$, a situation fully analyzed in what follows.
%\foot{SK9: added}
The main result of this paper is the following theorem.
%We say that a map $f_{a,b}$ {\em has a reduction theory} if properties (1)-(3) hold.
\begin{result} There exists an explicit one-dimensional Lebesgue measure zero, uncountable set $\mathcal E$ that lies on the diagonal boundary $b= a+ 1$ of $\mathcal P$ such that:
\begin{itemize}
\item[(a)] for all $(a,b)\in\mathcal P\setminus\mathcal E$ the map $F_{a,b}$ has an attractor $D_{a,b}$ satisfying properties (1) and (2) above;
\item[(b)] for an open and dense  set in $\mathcal P\setminus\mathcal E$ property (3), and hence the Reduction theory conjecture, holds. For the rest of $\mathcal P\setminus\mathcal E$ property (3) holds for almost every point of $\bar\R^2\setminus \Delta$.
\end{itemize}
\end{result}
%\foot{SK12: next paragraph was completely changed}

We point out that this approach gives explicit conditions for the set $D_{a,b}$ to have finite rectangular structure that are satisfied, in particular, for  all pairs $(a,b)$ in the interior of  the maximal parameter set $\mathcal P$. 
At the same time, it provides an effective algorithm for finding $D_{a,b}$, independent of the complexity of its boundary (i.e., number of horizontal segments). The simultaneous properties satisfied by $D_{a,b}$, attracting set and bijectivity domain for $F_{a,b}$, is an essential feature that has not been exploited in earlier works. This approach  makes the notions of reduced geodesic and dual expansion natural and transparent, with a potential for generalization to other Fuchsian groups. We remark that 
for ``$\alpha$-transformations" \cite{N1,LM}, explicit descriptions of the domain of the natural extension maps have been obtained only for a subset of the parameter interval $[0,1]$ (where the boundary has low complexity). 

The paper is organized as follows. In Section \ref{s:2} we  develop the theory of $(a,b)$-continued fractions associated to the map $f_{a,b}$. In Section \ref{s:3} we prove that the natural extension map $F_{a,b}$ possesses a {\em trapping region}; it will be used in Section \ref{s:6} to study the attractor set for $F_{a,b}$.
In Section \ref{s:4} we further
study the map $f_{a,b}$. Although it
is discontinuous at $x=a$, $b$, one can look at two orbits of each of the discontinuity points. For generic $(a,b)$, these orbits meet after finitely many steps, forming a {\em cycle} that can be {\em strong} or {\em weak}, depending on whether or not the product over the cycle is equal to the identity transformation.
The values appearing in these cycles play a crucial role in the theory. Theorems \ref{b-cycle} and \ref{a-cycle} give necessary and sufficient conditions for $b$ and $a$ to have the {\em cycle property}. In Section \ref{s:5} we introduce the {\em finiteness condition} using the notion of {\em truncated orbits} and prove that under this condition the map $F_{a,b}$ has a bijectivity domain $A_{a,b}$ with a finite rectangular structure that can be ``computed" from the data $(a,b)$ (Theorem \ref{thm:recstructure}). In Section \ref{s:6} we define the attractor for the map $F_{a,b}$ by iterating the trapping region, and
identify it with the earlier constructed set $A_{a,b}$ assuming the finiteness condition (Theorem \ref{attractor}). In Section \ref{s:7}
we prove that the Reduction theory conjecture holds under the assumption that both $a$ and $b$ have the strong cycle property, and that under the finiteness condition property, (3) holds for almost every point of $\bar\R^2\setminus\Delta$.
In Section \ref{s:8} we prove that the finiteness condition holds for all $(a,b)\in\P$ except for
an uncountable set of one-dimensional Lebesgue measure zero that lies on the boundary $b=a+1$ of $\P$, and we present a complete description of this exceptional set. We conclude by showing that the set of $(a,b)\in\mathcal P$ where $a$ and $b$ have the strong cycle property is open and dense in $\mathcal P$.
And, finally, in Section \ref{s:9} we show how these results can be applied to the study of invariant measures and ergodic properties of the associated Gauss-like maps.

\section{Theory of $(a,b)$-continued fractions}\label{s:2}
Consider $(a,b)\in \P$. The map $f_{a,b}$ defines what we call \emph{$(a, b)$-continued fractions} using 
a generalized integral part function $\lfloor x\rceil_{a,b}$ : 
%As mentioned in the previous section,  using the first return map of $f$ to the interval $(a,b)$ one can introduce a two-parameter family of continued fractions which we describe in what
%follows.  Consider the set
%$$\Delta=\{(a, b): -1\le  a\le 0\le b\le 1,b-a\ge
%1\}\,.$$ We define for $(a,b)\in \Delta$
%the following integer-valued function: 
for any real $x$, let
\begin{equation}
\lfloor x\rceil _{a,b}=\begin{cases}
\lfloor x-a \rfloor &\text{if } x< a\\
0 & \text{if } a\le x<b\\
\lceil x-b \rceil & \text{if } x\ge b\,,
\end{cases}
\end{equation}
where $\lfloor x\rfloor$ denotes the integer part of $x$ and $\lceil x\rceil=\lfloor x\rfloor+1$.

Let us remark that the first return map of $f_{a,b}$ to the interval $[a,b)$, $\hat f_{a,b}$,  is given by the function 
\[
\hat f_{a,b}(x)=-\frac{1}{x}-\left\lfloor-\frac{1}{x}\right\rceil_{a,b}=T^{-\lfloor-1/x\rceil_{a,b}}S(x) \text{ if } x\ne 0, f(0)=0.
\]

We prove that any irrational number $x$ can be expressed in a unique way as an infinite $(a,b)$-continued fraction
\[
x=n_0-\cfrac{1}{n_1 -\cfrac{1}{n_2-\cfrac{1}{\ddots}}}
\]
which we will denote by  $\lfloor n_0,n_1,\dots \rceil_{a,b}$ for short. %\foot{SK: why do we use ";" after the first digit? in section 8 we don't and in RA we don't - need to change}
The
``digits" $n_i$, $i\ge 1$,  are non-zero integers determined recursively by
\begin{equation}\label{a-b}
n_0=\lfloor x\rceil_{a,b},\,x_1=-\frac1{x-n_0},\text{ and }
n_{i}=\lfloor x_i \rceil_{a,b},\,x_{i+1}=- \frac1{x_i-n_i}.
\end{equation}
In what follows, the notation $(\alpha_0,\alpha_1,\dots,\alpha_k)$ is used to write formally a ``minus" continued fraction expression, where $\alpha_i$ are real numbers.%\foot{IU: added}
%$n_{i+1}=( x_{i+1}),\,x_{i+1}=- \frac1{x_i-n_i}$, starting with
%$n_0=\lfloor x\rceil_{a,b}$ and $x_1=-\frac1{x-n_0}$.%\foot{Replace $n_1$ by $x_1$!}
%$$r_k=\lfloor n_0,n_1,\dots, n_n\rceil_{a,b}:=n_0-\cfrac{1}{n_1-\cfrac{1}{n_2-\cfrac{1}{\ddots-\cfrac{1}{n_n}}}}$$

\begin{thm}\label{convergence}
Let $x$ be an irrational number, $\{n_i\}$ the associated sequence of integers defined by (\ref{a-b})
and
\[
r_k=( n_0,n_1,\dots, n_k)\,.%:=a_0-\cfrac{1}{a_1-\cfrac{1}{a_2-\cfrac{1}{\ddots-\cfrac{1}{a_n}}}}.
\]
Then the sequence $r_k$ converges to $x$. 
\end{thm}

\begin{proof}\footnote{The authors proved initially the convergence statement assuming $-1\le a \le 0 \le b\le 1$, and two Penn State REU students, Tra Ho and Jesse Barbour, worked on the proof for $a,b$ outside of this compact triangular region.  The unified proof presented here uses  some of their ideas.} We start by proving that none of the pairs of type  $(p,1)$, $(-p,-1)$, with $p\ge 1$
are allowed to appear as consecutive entries of the sequence $\{n_i\}$. Indeed, if $n_{i+1}=1$, then 
$$b\le x_{i+1}=-\frac{1}{x_i-n_i}<b+1\,,$$
therefore $-\displaystyle\frac{1}{b}\le x_i-n_i<-\frac{1}{b+1}\le (b-1)$, and $n_i<0$.  If $n_{i+1}=-1$, then
$$a-1\le x_{i+1}=-\frac{1}{x_i-n_i}<a\,,$$
so $-\displaystyle\frac{1}{a-1}\le x_i-n_i<-\frac{1}{a}$. But $a+1\le -\frac{1}{a-1}$, thus $n_i>0$.

%To show that the pair $(1,p)$ is not allowed, notice that if $n_i=1$, then $b-1\le x_i-n_i<b$. Clearly $n_{i+1}<0$ if $b\ge 1$. Assuming that $b<1$, one has $x_{i+1}<-1/b$ or $x_{i+1}\ge -1/(b-1)>(b+1)$. Both situations imply that $n_{i+1}\ne 2$.

With these two restrictions, the argument follows the lines of the proof for the classical case
of minus (backward) continued fractions  \cite{K3}, where $n_i\ge 2$, for all
$i\ge 1$. We define inductively two sequences of integers
$\{p_k\}$ and $\{q_k\}$ for $k\ge -2$:
\begin{equation}\label{pkqk}
\begin{split}
&p_{-2}=0\;,\;p_{-1}=1\;;\; p_{k}=n_{k}p_{k-1}-p_{k-2}\; \text{ for } k\ge 0\\
&q_{-2}=-1\;,\;q_{-1}=0\;;\;q_{k}=n_{k}q_{k-1}-q_{k-2}\; \text{ for } k\ge 0\;.
\end{split}
\end{equation}
We have the following properties: %are proved by induction (for $n\ge 0$):\foot{IU: changed property iv}
\begin{itemize}
\item[(i)] there exists $l\ge 1$ so that $|q_l|<|q_{l+1}|<\dots<|q_k|<\dots$;
\item[(ii)] $(n_0,n_1,\dots,n_k,\alpha)=\displaystyle\frac{\alpha p_k-p_{k-1}}{\alpha q_k-q_{k-1}}$, for any real number $\alpha$;%, $|\alpha|\ge 1$; 
\item [(iii)] $p_kq_{k+1}-p_{k+1}q_k=1$;
%\item[(iv)] $n_0-1\le \lfloor n_0;n_1,\dots,n_k,\dots\rceil_{a,b}\le n_0+1$.
\end{itemize}
Let us prove property (i). Obviously $1=q_0\le |q_1|=|n_1|$, $q_2=n_2q_1-q_0=n_2n_1-1$. Notice that $|q_2|>|q_1|$ unless
$n_1=1, n_2=2$ or $n_1=-1$, $n_2=-2$. We analyze the situation $n_1=1$, $n_2=2$. This implies that $q_3=n_3(n_2n_1-1)-n_1=n_3-n_1$, so $|q_3|>|q_2|$, unless $n_3=2$.  Notice that it is impossible to have $n_i=2$ for all $i\ge 2$, because $x$ is irrational and the minus continued fraction expression consisting only of two's, $(2,2,\dots)$, has numerical value $1$. Therefore,  there exists $l\ge 1$ so that
$n_{l+1}\ne 1,2$. This implies that $|q_{l+1}|>|q_l|$. We continue to proceed by induction. Assume that property (i) is satisfied up to $k$-th term, $k>l$. If $|n_{k+1}|\ge 2$, then %\foot{SK13: misprint - replaced $q_{k+1}$ by $n_{k+1}$!}
$$
|q_{k+1}|\ge |n_{k+1}|\cdot |q_k|-|q_{k-1}|\ge 2|q_k|-|q_{k-1}|>|q_{k}|\,.
$$
If $n_{k+1}=1$, then
$q_{k+1}=q_k-q_{k-1}$. Since $q_k=n_{k}q_{k-1}-q_{k-2}$ with $n_k<0$, one gets
$$q_{k-1}=\frac{q_k+q_{k-2}}{n_{k}}\,.$$
We analyze the two possible situations
\begin{itemize}
\item If $q_k>0$ then $|q_{k-2}|<q_k$, so $q_k+q_{k-2}>0$ and $q_{k-1}<0$. This implies that
$q_{k+1}=q_k-q_{k-1}\ge q_{k}>0\,.$
\item If $q_k<0$, then $|q_{k-2}|<-q_k$, so $q_k+q_{k-2}<0$ and $q_{k-1}>0$. This implies that
$q_{k+1}=q_k-q_{k-1}<q_k<0\,.$
\end{itemize}
Thus $|q_{k}|<|q_{k+1}|$. A similar argument shows that the inequality remains true if $n_{k+1}=-1$.

Properties (i)--(iii) show that $r_k=p_k/q_k$ for $k\ge 0$.
Moreover, the sequence $r_k$ is a Cauchy sequence because
$$|r_{k+1}-r_k|=\frac{1}{|q_kq_{k+1}|}\le \frac{1}{(k-l)^2}\, \text{ for } k>l.$$
Hence $r_k$ is convergent.

In order to prove that $r_k$ converges to $x$, we write $x=(n_0,n_1,\dots,n_k,x_{k+1})$, and
%where $\alpha=(n_{k+1},n_{k+2},\dots)$
look only at those terms $(n_0,n_1, \dots,n_k,x_{k+1})$ with $|x_{k+1}|\ge 1$. 
There are infinitely many such terms: indeed, if $-1\le a<b\le 1$, then $|x_{k+1}|\ge 1$ for all $k\ge 1$; if $a<-1$, and $|x_{k+1}|<1$, then $b\le x_{k+1}<1$, so $x_{k+2}=-1/(x_{k+1}-1)\ge 1$; if $b>1$, and $|x_{k+1}|<1$, then $-1<x_{k+1}<a$, so $x_{k+2}=-1/(x_{k+1}+1)\ge 1$. %\foot{IU: added}
Therefore, the corresponding subsequence $r_k=p_k/q_k$ satisfies
\begin{equation*}
\begin{split}\left|\frac{p_k}{q_k}-x\right|=&\left|\frac{p_k}{q_k}-\frac{p_k x_{k+1}-p_{k-1}}{q_k x_{k+1}-q_{k-1}}\right|
=\frac{1}{|q_k(q_k x_{k+1}-q_{k-1})|}\\ \le & \frac{1}{|q_k|(|q_k||x_{k+1}|-|q_{k-1}|)}\le \frac{1}{|q_k|}\rightarrow 0.
\end{split}
\end{equation*}
We showed that the convergent sequence $r_k=p_k/q_k$ has a subsequence convergent to $x$, therefore the whole sequence converges to $x$.
%Using property (iv), one has
%$$(n_0;2,2,\dots)\le \lfloor n_0,n_1,\dots\rceil\le (n_0;-2,-2,\dots)$$ which implies that
%$n_0-1\le \lfloor n_0,n_1,\dots\rceil_{a,b}\le n_0+1$. (We used here that
%$(2;2,2,\dots)=1$ and $(-2;-2,-2,\dots)=-1$.)
\end{proof}
\begin{rem}
One can construct $(a,b)$-continued fraction expansions for rational numbers, too. However, such expansions will terminate after finitely many steps if $b\ne 0$. If $b=0$, the expansions of rational numbers will end with a tail of $2$'s, since $0=(1,2,2,\dots)$.%\foot{SK13: modified remark}
\end{rem}
\begin{rem} \label{quadratic}It is easy to see that if the $(a,b)$-continued fraction expansion of a real number is eventually periodic, then the number is a quadratic irrationality.
\end{rem}
It is not our intention to present in this paper some of the typical number theoretical results that can be derived for the class of  $(a,b)$-continued fractions. However, we state and prove a simple version about $(a,b)$-continued fractions with ``bounded digits". For the regular continued fractions, this is a classical result due to Borel and Bernstein (see \cite[Theorem 196]{HW} for an elementary treatment).  We are only concerned with %\foot{SK13: replaced ``about" by ``with"}
$(a,b)$-expansions that are written with two consecutive digits, a result explicitly needed in Sections \ref{s:7} and \ref{s:8}.
%\foot{IU: added}

\begin{prop}\label{bdigits1}
The set $\Gamma^{(m)}_{a,b}=\{x=\lfloor 0, n_1,n_2,\dots\rceil_{a,b} \:|\: n_k\in\{m,m+1\}\}$ 
%and  $\Gamma^{-m}=\{x=\lfloor 0, n_1,n_2,\dots\rceil_{a,b}\:|\: n_i\in\{-m,-m-1\}\}$ 
has zero Lebesgue measure for every $m\ge 1$.
\end{prop}
\begin{proof}
First, notice that if $m=1$, then the set $\Gamma^{(1)}_{a,b}$  has obviously zero measure, since the pairs $(2,1)$ and $(-2,-1)$ are not allowed in the $(a,b)$-expansions. 

Assume $m\ge 2$. Notice that $\Gamma^{(m)}_{a,b}\subset \Gamma^{(m)}_{0,-1}$ since a formal continued fraction $x=(0, n_1,n_2,\dots)$ with $n_k\in\{m,m+1\}$ coincides with its ``minus" (backward) continued fraction expansion ($a=-1, b=0$), $x=\lfloor 0, n_1,n_2,\dots\rceil_{-1,0}$. The reason is that any sequence of digits $n_i\ge 2$ gives a valid ``minus" continued fraction expansion.

% \:|\: n_k\in\{m,m+1\}\} (backward continued fraction).
%Notice that $b<1$ in order for the set $\Gamma^{(m)}$ to be nonempty; otherwise, if $b\ge 1$ and  $n_k>0$, then $x_k-n_k>b-1$ so $x_{k+1}<0$ and $n_{k+1}<0$.  

In what follows, we study the set $\Gamma^{(m)}_{0,-1}$. For practical reasons we will drop the subscript $(0,-1)$.  
It is worth noticing that the result for $\Gamma^{(m)}_{0,-1}$ does not follow automatically from the result about regular continued fractions, since there are numbers for which the $(0,-1)$-expansion has only digits $2$ and $3$, while the regular continued fractions expansion has unbounded digits.
We follow the approach of \cite[Theorem 196]{HW} and estimate the size of the set $\Gamma^{(m)}_{n_1,n_2,\dots,n_k}\subset \Gamma^{(m)}$ with the digits $n_1,n_2,\dots,n_k\in \{m,m+1\}$ being fixed. In this particular case, the recursive relation \eqref{pkqk} implies that $1=q_1<q_2<\dots<q_k$.
%Any point $x\in\Gamma^m_1$ satisfies
%$$m+b-1\le -\frac{1}{x}<m+b+1, \text{ so } -\frac{1}{m+b-1}\le x<-\frac{1}{m+b+1}$$
If $x\in \Gamma^{(m)}_{n_1,n_2,\dots,n_k}$, then 
$$(0,n_1,n_2,\dots ,n_k-1)\le x<(0, n_1,n_2,\dots, n_k)\,.$$
Using property (iii), the endpoints of such an interval $I^{(m)}_{n_1,\dots,n_k}$ are given by 
$$
\frac{(n_k-1)p_{k-1}-p_{k-2}}{(n_k-1)q_{k-1}-q_{k-2}}\text{ , }
\frac{n_kp_{k-1}-p_{k-2}}{n_kq_{k-1}-q_{k-2}}
$$
and the length of this interval is
\[
%\begin{split}
l(I^{(m)}_{n_1,\dots,n_k})=\frac{1}{(n_kq_{k-1}-q_{k-2})((n_k-1)q_{k-1}-q_{k-2})}=\frac{1}{q_{k}(q_k-q_{k-1})}
%\end{split}
\]
by using that $p_{k-2}q_{k-1}-p_{k-1}q_{k-2}=1$ and $q_k=n_kq_{k-1}-q_{k-2}$. 

Denote by $\Gamma^{(m)}_k$ the set of numbers in $[-1,0)$ with $(-1,0)$-continued fraction digits $n_1, n_2, \dots$, $n_k\in\{m,m+1\}$. The set $\Gamma^{(m)}_k$ is part of the set
$$I^{(m)}_k=\bigcup_{n_1,\dots, n_k\in\{m,m+1\}}I^{(m)}_{n_1,\dots,n_k}\,.$$ We have the following relation:
$$
I^{(m)}_{k+1}=\bigcup_{n_1,\dots, n_k\in\{m,m+1\}}I^{(m)}_{n_1,\dots,n_k,m}\cup I^{(m)}_{n_1,\dots,n_k,m+1}
$$
If $x$ lies in $I^{(m)}_{n_1,\dots,n_k,m}\cup I^{(m)}_{n_1,\dots,n_k,m+1}$, then
$$(0,n_1,n_2,\dots ,n_k,m-1)\le x<(0, n_1,n_2,\dots, n_k,m+1)\,.$$
The length of this interval is 
$$l(I^{(m)}_{n_1,\dots,n_k,m}\cup I^{(m)}_{n_1,\dots,n_k,m+1})=\frac{2}{((m+1)q_{k}-q_{k-1})((m-1)q_{k}-q_{k-1})}
$$
Now we estimate the ratio
\begin{equation*}
\begin{split}
\frac{l(I^{(m)}_{n_1,\dots,n_k,m}\cup I^{(m)}_{n_1,\dots,n_k,m+1})} {l(I^{(m)}_{n_1,n_2,\dots,n_k})} &=  \frac{2q_{k}(q_k-q_{k-1})}{((m+1)q_{k}-q_{k-1})((m-1)q_{k}-q_{k-1})}\\
&\le \frac{2q_{k}}{(m+1)q_{k}-q_{k-1}}\\
&\le \frac{2q_{k}}{3q_{k}-q_{k-1}}=\frac{2}{3-q_{k-1}/q_k}\\
&\le\frac{2k}{2k+1}
%&= \frac{2(1+bq_{k-1}/q_k)}{(m+b+1)-q_{k-1}/q_k}
%&\le 2\frac{1+b}{(m+b+1)(m+b-1)}\le \frac{2}{m+b+1}
\end{split}
\end{equation*}
%\frac{2q_k}{3q_k-q_{k-1}}=\frac{2}{3-q_{k-1}/q_k}\le\frac{2k}{2k+1}
%&\le \frac{2(q_{k}+bq_{k-1})}{((m+b+1)q_{k}-q_{k-1})}
%&= \frac{2(1+bq_{k-1}/q_k)}{(m+b+1)-q_{k-1}/q_k}
%&\le 2\frac{1+b}{(m+b+1)(m+b-1)}\le \frac{2}{m+b+1}
%\end{split}
%\end{equation*}
since $\displaystyle\frac{q_{k-1}}{q_{k}}\le \frac{k-1}{k}$. Indeed, if $n_1=\dots=n_k=2$, then $q_{k-1}/q_k=(k-1)/k$; if some $n_j>2$, then $q_{k-1}/q_k\le 1/2$ from \eqref{pkqk}. This proves that for every $k\ge 1$
$$
I^{(m)}_{k+1}\le \frac{2k}{2k+1}I^{(m)}_k
$$
%A simple calculation also shows that 
%$$l(I^{(m)}_1)=\frac{2}{(m+b+1)(m+b-1)}\le\frac{2(1+b)}{m+b}\,.$$
so
$$l(I^{(m)}_{k})\le \frac{2\cdot 4 \cdots (2k-2)}{3\cdot 5\cdots(2k-1)} \cdot l(I^{(m)}_1)\longrightarrow 0 \text { as } k\rightarrow \infty.$$

Therefore, in all cases, $l(I^{(m)}_{k})\ra 0$ as $k\ra \infty$. Since $\Gamma^{(m)}\subset I^{(m)}_{k}$ for every $k\ge 1$, the proposition follows. 
\end{proof}
\begin{rem}\label{bdigits2}
By a similar argument, the set $\Gamma^{(-m)}_{a,b}=\{x=\lfloor 0, n_1,n_2,\dots\rceil_{a,b}\:|\: n_k\in\{-m,-m-1\}\}$ 
has zero Lebesgue measure for every $m\ge 1$.
\end{rem}

\section{Attractor set for $F_{a,b}$}\label{s:3}
%We study the reduction map (a  {\em natural extension} of some sort) $f:\bar\R^2\ra \bar\R^2$ defined by
%\begin{equation}
%f(x,y)=\begin{cases}
%(x-1,y-1)  &\text{ if } y\ge b\\
%(-1/x,-1/y) &\text{ if } a<y<b\\
%(x+1,y+1)  &\text{ if } y\le a
%\end{cases}\,.
%\end{equation}
The reduction map $F_{a,b}$ defined by \eqref{Fab} has a trapping domain, i.e. a closed
set $\Theta_{a,b}\subset \bar\R^2\setminus \Delta$ with the following properties:
\begin{itemize}
\item[(i)] for every pair $(x,y)\in \bar\R^2\setminus \Delta$, there exists a positive integer $N$
such that $F_{a,b}^N(x,y)\in \Theta_{a,b}$;
\item[(ii)] $F_{a,b}(\Theta_{a,b})\subset \Theta_{a,b}$.
\end{itemize}

\begin{thm}\label{Delta-trapping} The region $\Theta_{a,b}$ consisting of two connected components (or one if $a=0$ or $b=0$)
defined as
\begin{equation*}
\Theta^u_{a,b}=
\begin{cases}
[-\infty,-1]\times[b-1,\infty]\cup [-1,0]\times[-\frac{1}{a},\infty] & \text{ if } b\ge 1, a\ne 0\\
\emptyset & \text{ if } a=0\\
\begin{split}
[-\infty,-1]\times[b-1,\infty]&\cup [-1,0]\times[\min  (-\frac{b}{b-1}, -\frac1 a),\infty]\\
&\cup[0,1]\times[-\frac{1}{b-1},\infty]
\end{split} & \text{ if } 0<b<1
\end{cases}
\end{equation*}

\begin{equation*}
\Theta^l_{a,b}=
\begin{cases}
[0,1]\times[-\infty,-\frac{1}{b}]\cup [1,\infty]\times[-\infty,a+1] & \text{ if } a\le -1, b\ne 0\\
\emptyset & \text{ if } b=0\\
\begin{split}
[-1,0]\times[-\infty,-\frac{1}{a+1}]&\cup [0,1]\times[-\infty,\max (\frac{a}{a+1}, -\frac{1}{b})]\\
&\cup[1,\infty]\times[-\infty,a+1]
\end{split} &\text{ if } a>-1
\end{cases}
\end{equation*}
is the trapping region for the reduction map $F_{a,b}$.
\end{thm}

\begin{figure}[htb]
\includegraphics[scale=0.6]{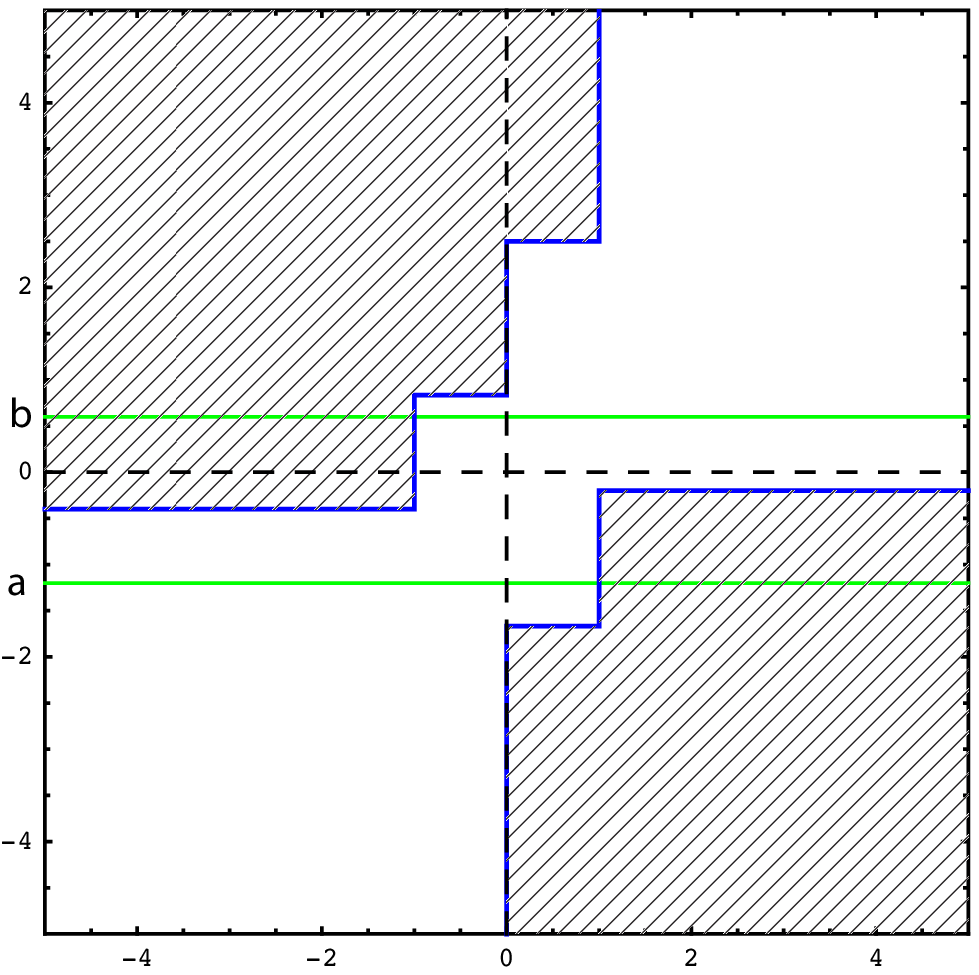} $\quad$ \includegraphics[scale=0.6]{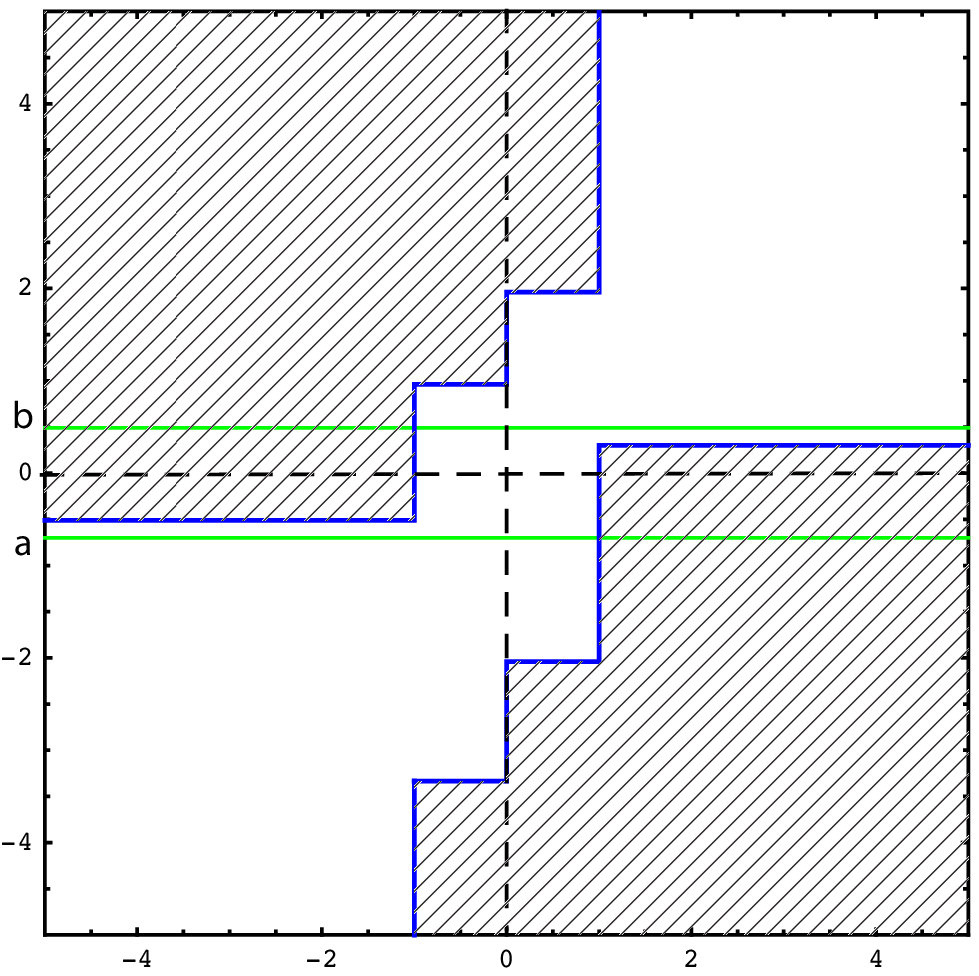}
\caption{Typical trapping regions: case $a<-1, 0<b<1$ (left); case $-1<a<0<b<1$ (right)}
\label{trap}
\end{figure}

\begin{proof}
%From here on, we will drop the subscripts for $F$ and $\Theta_{a,b}$.
The fact that the region $\Theta_{a,b}$ is $F_{a,b}$-invariant is verified by a direct calculation.
We focus our attention on the attracting property of $\Theta_{a,b}$. Let $(x,y)\in \R^2\setminus\Delta$, write 
$y=\lfloor n_0,n_1,\dots\rceil_{a,b}$, and construct the following sequence of
real pairs $\{(x_k,y_k)\}$ ($k\ge 0$) defined by $x_0=x$, $y_0=y$
and:
$$y_{k+1}=ST^{-n_k}\dots ST^{-n_1}ST^{-n_0}y\,,\quad x_{k+1}=ST^{-n_k}\dots ST^{-n_1}ST^{-n_0}x\,.$$
If $y$ is rational and its $(a,b)$-expansion terminates $y=\lfloor n_0,n_1,\dots,n_l\rceil_{a,b}$, then $y_{l+1}=\pm \infty$, 
so $(x,y)$ lands in $\Theta_{a,b}$ after finitely many iterations. If $y$ has an infinite $(a,b)$-expansion, then 
$y_{k+1}=\lfloor n_{k+1},n_{k+2},\dots\rceil_{a,b}$, and $y_{k+1}\ge -1/a$ or $y_{k+1}\le -1/b$ for $k\geq 0$. Also,
\begin{equation*}
\begin{split}
y&=T^{n_0}ST^{n_1}S\dots T^{n_k}S(y_{k+1})=\frac{p_{k}y_{k+1}-p_{k-1}}{q_{k}y_{k+1}-q_{k-1}}\\
x&=T^{n_0}ST^{n_1}S\dots T^{n_k}S(x_{k+1})=\frac{p_{k}x_{k+1}-p_{k-1}}{q_{k}x_{k+1}-q_{k-1}}\,,
\end{split}
\end{equation*}
hence
\begin{equation*}\label{eq:uk}
x_{k+1}=\frac{q_{k-1}x-p_{k-1}}{q_kx-p_k}=\frac{q_{k-1}}{q_k}+\frac{1}{q_k^2(p_k/q_k-x)}=\frac{q_{k-1}}{q_k}+\varepsilon_{k}
\end{equation*}
where $\varepsilon_{k}\ra 0$. This shows that for $k$ large enough $x_{k+1}\in [-1,1]$. We proved that there exists $N>0$, such that
$$F_{a,b}^{N}(x,y)=ST^{-n_k}\dots ST^{-n_1}ST^{-n_0}(x,y)\in [-1,1]\times ([-1/a,\infty]\cup [-\infty,-1/b])\,.$$
The point $F^N_{a,b}(x,y)=:(\tilde x,\tilde y)$ belongs to $\Theta_{a,b}$, unless $b<1$ and $(\tilde x,\tilde y)\in [0,1]\times [-1/a,-1/(b-1)]$ or
$a>-1$ and $(\tilde x,\tilde y)\in [-1,0]\times [-1/b,-1/(a+1)]$.

Let us study the next iterates of  $(\tilde x,\tilde y)\in [0,1]\times [-1/a,-1/(b-1)]$. If $\tilde y\ge  b+1$ then 
$$F^2_{a,b}(\tilde x,\tilde y)=(\tilde x-2,\tilde y-2)\in [-1,1]\times [b-1,\infty]\,,$$ so $F^2_{a,b}(\tilde x,\tilde y)\in \Theta_{a,b}$. If it so happens that $-1/a\le \tilde y<b+1$, then 
$$F_{a,b}(\tilde x,\tilde y)=(\tilde x-1,\tilde y-1)\in [-1,0]\times [0,b]$$ and
$$F^2_{a,b}(\tilde x,\tilde y)=ST^{-1}(\tilde x,\tilde y)\in [0,\infty]\times [-1/b,\infty] \subset \Theta_{a,b}\,.$$ 
Similarly,  if $(x,y)\in [-1,0]\times [-1/b,-1/(a+1)]$, then $F^2_{a,b}(x,y)\in \Theta_{a,b}$.

Notice that if $a=0$, then $y_{k+1}\le -1/b$ for all $k\ge 0$ (so $\Theta_{a,b}^u=\emptyset$) and if $b=0$, then $y_{k+1}\ge -1/a$ for al $k\ge 0$ (so $\Theta_{a,b}^l=\emptyset$).%\foot{SK20: added $_{a,b}$ to $\Theta$}
\end{proof}

Using the trapping region described in Theorem \ref{Delta-trapping} we define the associated \emph{attractor set} 
\begin{equation}\label{def-atrac}
D_{a,b}=\bigcap_{n=0}^\infty D_n,
\end{equation}
where $D_n=\bigcap_{i=0}^n F_{a,b}^i(\Theta_{a,b})$.

\begin{rem}
In the particular cases when $a=0$ and $b\ge 1$, or $b=0$ and $a\le -1$ or $(a,b)=(-1,1)$ the trapping regions 
%\foot{IU: new remark}
\begin{eqnarray*}
\Theta_{0,b}&=&[-1,0]\times [-\infty,-1]\cup [0,1]\times[-\infty,0]\cup[1,\infty]\times [-\infty,1]\\
\Theta_{a,0}&=&[-\infty,-1]\times[-1,\infty]\cup [-1,0]\times [0,\infty]\cup [0,1]\times [1,\infty]\\\
\Theta_{-1,1}&=&
[-\infty,-1]\times[-1,\infty]\cup [-1,0]\times [1,\infty]\\
& &\!\!\!\cup \,\,[0,1]\times [-\infty,-1]\cup [1,\infty]\times [-\infty,0]
\end{eqnarray*}
are also bijectivity domains for the corresponding maps $F_{a,b}$. Therefore, in these cases the attractor $D_{a,b}$ coincides with the trapping region $\Theta_{a,b}$, so the properties mentioned in the introduction are obviously satisfied. In what follows, all our considerations will exclude these degenerate cases. 
\end{rem}
%\foot{SK20: replaced ``particular" by ``degenerate"}

\section{Cycle property}\label{s:4}
%We will construct an explicit cycle sequences for a generic pair  $a,b\in \Delta$, each cycle having $T$-leg and $T^{-1}$-leg.

In what follows, we simplify the notations for $f_{a,b}$, $\lfloor,\cdot\rceil_{a,b}$, $\hat f_{a,b}$ and $F_{a,b}$ to $f$, $\lfloor,\cdot\,\rceil$, $\hat f$ and $F$, respectively, assuming implicitly their dependence on parameters $a,b$. We will use the notation $f^n$ (or $\hat f^n$) for the $n$-times composition operation of $f$ (or $\hat f$). Also, 
for a given point $x\in (a,b)$ the notation $\hat f^{(k)}$ means the transformation of type $T^iS$ ($i$ is an integer) such that
\[
\hat f^k(x)=\hat f^{(k)}\hat f^{(k-1)}\cdots \hat f^{(2)}\hat f^{(1)}(x),
\]
where $\hat f^{(1)}(x)=\hat f(x)$.

The map $f$ is discontinuous at $x=a$, $b$, however, we can associate to each $a$ and $b$ two forward orbits: to $a$ we associate the {\em upper orbit} $\OO_u(a)=\{f^{n}(Sa)\}$, and the {\em lower orbit}  $\OO_\ell(a)=\{f^n(Ta)\}$, and to $b$ --- the {\em lower orbit} $\OO_\ell(b)=\{f^n(Sb)\}$ and the {\em upper orbit} $\OO_u(b)=\{f^n(T^{-1}b)\}$. We use the convention 
%\foot{IU: resonance description} 
that if an orbit hits one of the discontinuity points $a$ or $b$, then the next iterate is computed according to the lower or upper location: for example, 
if the lower orbit of $b$ hits $a$, then the next iterate is $Ta$, %followed by $S$, $STa$;\foot{SK: why not to say the next iterate is $Ta$?}
%$Ta=a+1$ followed by $STa=-\frac{1}{a+1}$; 
if the upper orbit of $b$ hits $a$ then the next iterate is $Sa$.  
%\foot{SK: change; maybe it is better to delete $STa=-\frac{1}{a+1}$ and $S(a)=-1/a$ and just to leave the generators?}

Now we explore  the patterns in the above orbits. The following property plays an essential role in studying the map $f$.
\begin{defn}\label{def:cycles}
We say that the point $a$ has  the {\em cycle property} if for some non-negative integers $m_1,k_1$
\[
f^{m_1}(Sa)=f^{k_1}(Ta)=c_a.
\]
We will refer to 
the set
\[
\{Ta, fTa,\dots ,f^{k_1-1}Ta\}
\]
as the {\em lower side of the $a$-cycle}, to the set
\[
\{Sa, fSa,\dots ,f^{m_1-1}Sa\}
\]
as the {\em upper side of the $a$-cycle}, and to $c_a$ as the {\em end of the $a$-cycle}.
If the product over the $a$-cycle equals the identity transformation, i.e. %\foot{SK: why not ``identity transformation"? We say so in RA.}
\[
T^{-1}f^{-k_1}f^{m_1}S=\Id,
\]
we say that $a$ has {\em strong cycle property}, otherwise, we say that  $a$ has {\em weak cycle property}.

Similarly, we say that $b$ has  {\em cycle property} if for some  non-negative integers $m_2,k_2$
\[
f^{k_2}(Sb)=f^{m_2}(T^{-1}b)=c_b.
\]
We will refer to 
the set
\[
\{Sb, fSb,\dots ,f^{k_2-1}Sb\}
\]
as the {\em lower side of the $b$-cycle}, to the set
\[
\{T^{-1}b, fT^{-1}b,\dots ,f^{m_2-1}T^{-1}b\}
\]
as the {\em upper side of the $b$-cycle}, and to $c_b$ as the {\em end of the $b$-cycle}.
If the product over the $b$-cycle equals the identity transformation, i.e. %\foot{SK: why not ``identity transformation"? We say so in RA.}
\[
Tf^{-m_2}f^{k_2}S=\Id,
\]
%\foot{\textcolor{red}{ interchanged $k_2$ and $m_2$ to comply with previous notations}}
%and $c_b>Sa$ or $<Sb$,
we say that $b$ has {\em strong cycle property}, and otherwise we say that $b$ has {\em weak cycle property}.
\end{defn}

\begin{figure}[htb]
\psfrag{A}[l]{\small $a$}
\psfrag{B}[l]{\small $b$}
\psfrag{a}[c]{\small $B_{-1}$}
\psfrag{b}[c]{\small $B_{-2}$}
\psfrag{c}[l]{\small $B_3$}
\psfrag{d}[l]{\small $A_1$}
\psfrag{e}[c]{\small $A_2$}
\psfrag{f}[c]{\small $A_3$}
\psfrag{u}[l]{\small $\vdots$}
\psfrag{v}[l]{\small $\dots$}
\includegraphics[scale=1]{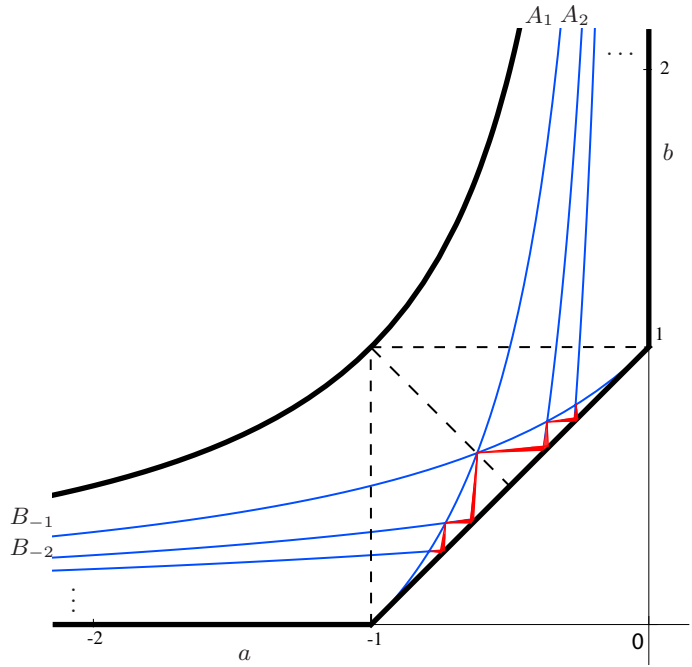}
\caption{The parameter set $\P$ and its partition}
\label{fig-P}
\end{figure}

It turns out that the cycle property is the prevalent pattern. It can be analyzed and described explicitly by partitioning the parameter set $\mathcal P$ based on the first digits of $Sb$, $STa$, and $Sa$, $ST^{-1}b$, respectively.  Figure \ref{fig-P} shows a part of the countable partitions, with $B_{-1},B_{-2},\dots$ denoting the regions where $Sb$ has the first digit $-1, -2,\dots$,  and $A_1,A_2,\dots$, denoting the regions where $Sa$ has the first digit $1,2,\dots$. 
For most of the parameter region, the cycles are short: everywhere except for the very narrow  triangular regions shown in Figure \ref{fig-P} the cycles for both $a$ and $b$ end after the first return to $[a,b)$.
However, there are Cantor-like recursive sets where the lengths of the cycles can be arbitrarily long. Part of this more complex structure, studied in details in Section \ref{s:8}, can be seen as very narrow triangular regions close to the boundary segment $b-a=1$.  

By symmetry of the parameter set $\P$ with respect to the line $b=-a$, $(a,b)\mapsto (-b,-a)$, we may assume that $b\leq -a$ and concentrate our attention to this subset of $\P$.

The structure of the set where the cycle property holds for $b$ is described next for the part of the parameter region with $0<b\le -a<1$. We make use extensively of the first return map $\hat f$.

\begin{thm}\label{b-cycle}%\foot{IU: $a>-1$}
Let $(a,b)\in\P$, $0<b\leq -a<1$ and $m\ge 1$ such that $a\le T^mSb<a+1$.
\begin{itemize}
\item[\bf{(I)}] Suppose that there exists $n\geq 0$ such that
\[
\quad\quad\quad \hat f^{k}T^mSb\in \Big(\frac{b}{b+1},a+1\Big) \text{ for } k<n, \text{ and } \hat f^{n}T^mSb\in \Big[a,\frac{b}{b+1}\Big].
\]
\begin{itemize}
\item[(i)] If $\hat f^{n}T^mSb\in (a,\frac{b}{b+1})$, then $b$
has the cycle property; the cycle property is strong if and only if $\hat f^{n}T^mSb\neq 0$. 
\item[(ii)] If $\hat f^{n}T^mSb=a$, then $b$ has the cycle property if and only if $a$ has the cycle property.
\item[(iii)] $\hat f^{n}T^mSb=b/(b+1)$, then $b$ does not  have the cycle property, but the orbits of $Sb$ and $T^{-1}b$ are periodic.
\end{itemize}
\item[\bf{(II)}] If $\hat f^{k}T^mSb\in (\frac{b}{b+1},a+1)$ for all $k\geq 0$, then $b$ does not have the cycle property.
\end{itemize}
\end{thm}

\begin{proof} {\bf{(I)}} 
In the case $m=1$, and assuming $a<TSb<a+1$ we have
\begin{equation}\label{m=1}
a<1-\frac1{b}< \frac{b}{b+1},
\end{equation}
and the cycle relation for $b$ can be explicitly described as
\begin{equation}\label{m1}
\begin{diagram}
        &             & \boxed{b-1}       &\rTo{S}     & \boxed{-\frac 1{b-1}}\\
        &\ruTo{T^{-1}}&                   &            &                         &\rdTo{T^{-1}}                 \\
\boxed{b}&             &                   &            &                         &         &\quad\boxed{c_b=\frac{b}{1-b}}\\
         &\rdTo{S}     &                   &            &                         &\ruTo{S}                          \\
         &             &\boxed{-\frac 1{b}}&\rTo{T}& \boxed{\frac{b-1}{b}}
\end{diagram}
\end{equation}

In the particular situation that $TSb=a$, %\foot{SK: misprint corrected!}
the lower orbit of $b$ hits $a$ and continues to $a+1$, while the upper orbit hits $\frac{b}{1-b}=-1/a$. This means that the iterates will follow the lower and upper orbits of $a$, respectively, thus statement (ii) holds.
Since the second inequality (\ref{m=1}) is strict, the case (iii) cannot occur.

\medskip

For the case $m=2$ (and assuming $T^2Sb\ne a$) we analyze the following situations: if $b<\frac12$, then
$
2-\frac1{b}<0,
$
and the cycle relation is
\begin{equation}\label{m=2}
\begin{diagram}
         &        & \boxed{b-1}       &\rTo{S}     & \boxed{-\frac 1{b-1}}   &\rTo{ST^{-2}}& \boxed{1+\frac{b}{1-2b}}\\
         &\ruTo{T^{-1}}&                   &            &                         &                  &                          &\rdTo{\;\;T^{-1}}\\
\boxed{b}&        &                   &            &                         &                  &                          &        & \boxed{c_b=\frac{b}{1-2b}}\\
    &\rdTo{S}&                   &            &                         &   & \ruTo(3,2){S}       &            &        &\\
    &        &\boxed{-\frac 1{b}}&\rTo{T^{2}}& \boxed{-\frac{1-2b}{b}}
\end{diagram}
\end{equation}
%Same observation as above holds for the case $T^2SB=a$.
If $b>\frac12$ we have
\[
0<2-\frac1{b}\leq \frac{b}{b+1},
\]
since we must also have $2-\frac1{b}<a+1$, i.e. $b\leq \frac1{1-a}$,
and
the cycle relation is
\begin{equation}
\begin{diagram}
        &             & \boxed{b-1}       &\rTo{S}     & \boxed{-\frac 1{b-1}}   &                  &                       &          &\\
        &\ruTo{T^{-1}}&                   &            &                         &                                          \rdTo(4,2){ST^{-2}} & \\
\boxed{b}&        &                   &            &                         &                  &                       &          & \boxed{c_b=1+\frac{b}{1-2b}}\\
         &\rdTo{S}&                   &            &                         &                  &                       & \ruTo_{T} \\
         &        &\boxed{-\frac 1{b}}&\rTo{T^{2}}& \boxed{-\frac{1-2b}{b}} & \rTo^S           & \boxed{\frac{b}{1-2b}}&
\end{diagram}
\end{equation}
The above cycles are strong.
If $b=\frac12$ the cycle relation is %\foot{IU: weak cycle case}
\begin{equation}
\begin{diagram}
         &        & \boxed{b-1}       &\rTo{S}     & \boxed{-\frac 1{b-1}}\\
         &\ruTo{T^{-1}}&              &            &                &\rdTo(2,2){T^{-2}} \\
\boxed{b}&        &                   &            &                &           &\boxed{c_b=-\frac{1-2b}{b}=0}\\
         &\rdTo{S}&                   &          & \ruTo(3,2){T^2}   &\\
         &        &\boxed{-\frac 1{b}}&            &                &
\end{diagram}
\end{equation}
%\foot{SK: added $0$ at the end of the weak cycle}
It is easy to check that this cycle is weak. In the particular situation when 
%\foot{replaced``that" by ``when"} 
$T^2Sb=a$, the lower orbit of $b$ hits $a$, and continues with $a+1$, while the upper orbit still hits $\frac{b}{1-2b}=-1/a$. This means that the iterates will follow the lower and upper orbits of $a$, respectively, and statement (ii) holds. 
The relation $2-\frac1{b}=\frac{b}{b+1}$ implies $b=\frac{-1+\sqrt{5}}2$ that does not have the cycle property and the orbits
of $Sb$ and $T^{-1}b$ are periodic; this is the only possibility for (iii) to hold.
\medskip

The situation for $m\geq 3$ is more intricate. First we will need the following lemmas.

\begin{lem} \label{xy} Suppose $STSx=y$. The following are true:
\begin{itemize}
\item[(a)] if $TSb\leq x<a$, then $b-1\leq y<\frac{a}{1-a}$;
\item[(b)] if $a\le x<\frac{b}{b+1}$, then $\frac{a}{1-a}\le y< b$;%\foot{IU: inequalities changed}
\item[(c)] if $\frac{b}{b+1}\le x< a+1$, then $b\le y<\frac{a}{1-a}+1$;
\item[(d)] if $x=0$, then $y=0$. %\foot{IU: added this line}
\end{itemize}
\end{lem}
%\foot{For weak cycles if $x=0$, then $y=0$.}
\begin{proof}  Applying $STS$ to the corresponding inequalities
we obtain

\smallskip

(a)
\(\qquad
b-1=STSTSb\leq y<STSa=\displaystyle\frac{a}{1-a}.
\)

\smallskip
(b)
\(\qquad
\displaystyle \frac{a}{1-a}=STSa\le y<STSTSTb=b
\)

\smallskip

(c)
\(\qquad 
b=STSTSTb\le y<STSTa=T^{-1}Sa\le \displaystyle\frac1{1-a}=\frac{a}{1-a}+1,
\)

\smallskip

\noindent where the last inequality is valid for $a\le \frac{1-\sqrt{5}}{2}$,
%\approx -0.61803$, 
which is true in the considered region $b\le\frac1{2-a}$. 
%Notice that this inequality is strict on both sides. 
Relation (d) is obvious.
\end{proof}

%Let us write
%\[
%\hat f^k=\hat f^{(k)}\hat f^{(k-1)}\cdots \hat f^{(1)}\hat f^{(0)}.
%\]
%Now we return to main proof.  Now we need the following lemma.
\begin{lem}\label{main-argument}
Suppose that for all $k<n$
\begin{equation}\label{continuation}
\frac{b}{b+1}<\hat f^{k}T^mSb<a+1.
\end{equation}
Then 
\begin{enumerate}
\item for $0\leq k\leq n$, in the lower orbit of $b$, $\hat f^{(k)}=T^mS$ or $T^{m+1}S$; in the upper orbit of $b$, $\hat f^{(k)}=T^{-i}S$ with $i=2$ or $3$;
\item
there exists $p>1$ such that
\begin{equation}\label{STS}
(STS)\hat f^{n}T^mS=(T^{-2}S)\hat f^{p}T^{-1}.
\end{equation}
\end{enumerate}
\end{lem}
\begin{proof} 
(1) Applying $T^mS$ to the inequality (\ref{continuation}), we obtain
\[
a-1\leq T^{m-1}Sb=T^mSTSTb<T^mS\hat f^k T^mSb\leq T^mSTa\leq T^mSb< a+1,
\]
therefore $\hat f^{(k+1)}=T^mS$ or $T^{m+1}S$. Since $\hat f^{(0)}=T^mS$, we conclude that $\hat f^{(k)}=T^mS$ or $T^{m+1}S$ for $0\leq k\leq n$.

(2) In order to determine the upper side of the $b$-cycle,
%\foot{SK: terminology changed}
we will use the following relation in the group $\sz$ obtained by concatenation of the ``standard" relations (from right to left)
\begin{equation}\label{stand}
(STS)T^iS=(T^{-2}S)^{i-1}T^{-1}\quad  (i\geq 1),
\end{equation}
and Lemma \ref{xy} repeatedly.

The proof is by  induction on $n$. For the base case $n=1$ we have
\[
\frac{b}{b+1}< T^mSb< a+1.
\]
Then for $1\leq i\leq m-1$ $T^iSb$ satisfies condition (a) of  Lemma \ref{xy}, hence
\[
b-1<(T^{-2}S)^{i-2}T^{-1}b<\frac{a}{1-a},
\]
 which means that on the  upper side of the $b$-cycle $\hat f^{(1)}=T^{-1}$ and $\hat f^{(i)}=T^{-2}S$ for $1<i\leq m-1$.
 Using (\ref{stand}) for $i=m$ we obtain
 \[
 (STS)T^mS=(T^{-2}S)^{m-1}T^{-1}=(T^{-2}S)\hat f^{m-2}T^{-1},
 \]
 i.e. (\ref{STS}) holds with $p=m-2$.
 Now suppose the statement holds for $n=n_0$, and for all $k<n_0+1$ we have
 \[
\frac{b}{b+1}<\hat f^{k}T^mSb< a+1.
\]
By the induction hypothesis, 
there exists $p_0>1$ such that
\begin{equation}\label{cont}
(STS)\hat f^{n_0}T^mS=(T^{-2}S)\hat f^{p_0}T^{-1}.
\end{equation}
But since
\[
\frac{b}{b+1}<\hat f^{n_0}T^mSb< a+1,
\]
condition (c) of  Lemma \ref{xy} is satisfied, and hence
\[
b<(T^{-2}S)\hat f^{p_0}T^{-1}b<\frac{a}{1-a}+1,
\]
which is equivalent to
\[
b-1<(T^{-3}S)\hat f^{p_0}T^{-1}b<\frac{a}{1-a},
\]
i.e. $\hat f^{p_0+1}=T^{-3}S$. Using the relation
\(
(STS)T^2S=T^{-1}(STS),
\)
we can rewrite (\ref{cont}) as
\begin{equation}\label{ind}
(STS)T^2S\hat f^{n_0}T^mS=(T^{-3}S)\hat f^{p_0}T^{-1}=\hat  f^{p_0+1}T^{-1}.
\end{equation}
Let $\hat f^{(p_0+1)}=T^qS$. We have proved in (1) that $q=m$ or $m+1$, hence $q\geq 3$.
Let 
\[
b_0=T^2S\hat f^{n_0}T^mSb \text{ and } c_0=(T^{-3}S)\hat f^{p_0}T^{-1}b.
\]
 Then by (\ref{ind})
$(STS)b_0=c_0$.
Using the relation $(STS)T=T^{-2}S(STS)$, we obtain
\[
(STS)T^i=(T^{-2}S)^i(STS),
\]
and therefore,
\begin{equation}\label{STSi}
(STS)T^ib_0=(T^{-2}S)^i(STS)b_0=(T^{-2}S)^ic_0
\end{equation}
Since for $0\leq i<q-2$ $T^ib_0$ satisfies condition (a) of  Lemma \ref{xy},
we conclude that 
\[b-1< (T^{-2}S)^ic_0<\frac{a}{1-a}.
\]
Therefore $\hat f^{(i)}=T^{-2}S$ for $p_0+1<i\leq p_0+q$, 
and (\ref{STSi}) for $i=q-2$ gives us the desired relation
\[
(STS)\hat f^{n_0+1}T^mS=(T^{-2}S)\hat f^{p_0+q}T^{-1}
\]
with $p=p_0+q$.
\end{proof}
Now we complete the proof of  the theorem. In what follows we introduce the notations
\[
I_\ell=\Bigl( a,\frac{b}{b+1}\Bigr)\,,\; I_u=\Bigl(\frac{a}{1-a},b\Bigr)
\]
% and
%\[
%I^-_u=(\frac{a}{1-a},0),\,I^+_u=(0,b),\,\,\,I_u=I^-_u\cup\{0\}\cup I^+_u.
%\] 
and write $\overline I_\ell$, $\overline I_u$ for the corresponding closed intervals.

(I) If
\(
\hat f^nT^mSb\in I_\ell
\), then
condition (b) of Lemma \ref{xy} is satisfied, and
\[
(T^{-2}S)\hat f^{p}T^{-1}b\in I_u.
\]
It follows that $\hat f^{(p+1)}=T^{-2}S$,
therefore (\ref{STS}) can be rewritten as
\[
(STS)\hat f^{n}T^mS=\hat f^{p+1}T^{-1},
\]
which means that we reached the end of the cycle. More precisely,

(i) if  $\hat f^{n}T^mSb\in (0,\frac{b}{b+1})$, then
\[
TS\hat f^{n}T^mSb=S\hat f^{p}T^{-1}b=c_b;
\]
$b-1<\hat f^{j}T^{-1}b<\frac{a}{1-a}$ for $j<p$,  and $\hat f^{p}T^{-1}b\in (0,b)$. In this case $c_b<Sb$.

If $\hat f^{n}T^mSb\in (a,0)$, then
\[
S\hat f^{n}T^mSb=T^{-1}S\hat f^{p}T^{-1}b=c_b;
\]
 $b-1<\hat f^{j}T^{-1}b<\frac{a}{1-a}$  for $j<p$, and $\hat f^{p}T^{-1}b\in (\frac{a}{1-a},0)$. 
 In this case $c_b>Sa$.
Since the cycle relation in both cases is equivalent to the identity (\ref{STS}), the cycle property is strong, and (i) is proved.
  
 If $\hat f^{n}T^mSb=0$, then
\[\hat f^{n}T^mSb=\hat f^{p}T^{-1}b=0
\]
is the end of the cycle; for $j<p$, $b-1<\hat f^{j}T^{-1}b<\frac{a}{1-a}$. In this case the cycle ends ``before" the identity (\ref{STS}) is complete, therefore the product over the cycle is not equal to identity, and the cycle is weak.

\medskip

(ii) If $\hat f^{n}T^mSb=a$, then following the argument in (i) and using relation \eqref{STS} we obtain that the upper orbit of $b$ hits $T^{-1}S\hat f^{p}T^{-1}b=S\hat f^{n}T^mSb=Sa=-1/a$, while the lower orbit hits the value $a+1$, hence $b$ satisfies the cycle property if and only if $a$ does. 
%\foot{IU: resonance case}

\medskip

(iii) If $\hat f^{n}T^mSb=\frac{b}{b+1}$, then following the argument in (i) we obtain
\[
(T^{-2}S)\hat f^{p}T^{-1}b=b.
\]
However, one needs to apply one more $T^{-1}$ to follow the definition of the map $f$, hence $\hat f^{(p+1)}=T^{-3}S$, not $T^{-2}S$, and the cycle will not close.
One also observes that in this case the $(a,b)$-expansions of $Sb$ and $T^{-1}b$ will be periodic, and therefore the cycle will never close.

(II) If
\[
\hat f^{k}T^mSb\notin \overline {I_\ell}
\]
 for all $k\geq 0$, by the argument in the part (I) of the proof, on the %$Sb$-leg,
lower orbit of $b$ %\foot{SK: terminology}
each $\hat f^{(k)}=T^q S$, where $q=m$ or $m+1$, and on the %$T^{-1}$-leg
upper orbit of $b$ %\foot{SK: terminology}
each $\hat f^{(p)}=T^{-r}S$, where $r=2$ or $3$, and
for all $p\geq 1$
\[
\hat f^pT^{-1}b\notin\overline{I_u}.
\]
This means that for all images under the original  map $f$ on the
%$Sb$-leg
lower orbit of $b$ %\foot{SK: terminology}
we have
\[
f^kSb\in\left(-1-\frac1{b},a\right)\cup \left(\frac{b}{b+1},a+1\right)
\]
while for the images on the
%$T^{-1}b$-leg,
upper orbit of $b$ %\foot{SK: terminology}
\[
f^kT^{-1}b\in\left(b-1,\frac{a}{1-a}\right)\cup\left(b,1-\frac1{a}\right).
\]
Since these ranges do not overlap, the cycle cannot close, and $b$ has no cycle property.
%\foot{A weak cycle property is an artifact that I want to get rid of. It corresponds to the first alternative, and one can extend it a little to get the strong cycle.}
\end{proof}

A similar result holds for the $a$-cycles. First, if $Sa$ has the first digit $1$, i.e. $b\le Sa<b+1$, then one can easily write the $a$-cycle, similarly to \eqref{m=1}. For the rest of the parameter region we have:%\foot{IU25: added}

%\foot{SK20: added space after $\P$}
\begin{thm}\label{a-cycle}
Let $(a,b)\in\P$, $0<b\leq -a<1$ with $Sa\ge b+1$ and $m\ge 1$ such that $a\le T^mSTa<a+1$.
\begin{itemize}
\item[\bf{(I)}] Suppose that there exists $n\geq 0$ such that
\[
\quad\quad\quad \hat f^{k}T^mSTa\in \Big(\frac{b}{b+1},a+1\Big) \text{ for } k<n, \text{ and } \hat f^{n}T^mSTa\in \Big[a,\frac{b}{b+1}\Big].
\]
\begin{itemize}
\item[(i)] If $\hat f^{n}T^mSTa\in (a,\frac{b}{b+1})$, then $a$
has the cycle property; the cycle property is strong if and only if $\hat f^{n}T^mSTa\neq 0$.
\item[(ii)] If $\hat f^{n}T^mSTa=a$, then $a$ does not  have the cycle property, but the $(a,b)$-expansions of $Sa$ and $Ta$ are eventually periodic.
\item[(iii)] $\hat f^{n}T^mSTa=b/(b+1)$, then $a$ has the cycle property if and only if $b$ has the cycle property.
\end{itemize}
\item[\bf{(II)}] If $\hat f^{k}T^mSTa\in (\frac{b}{b+1},a+1)$ for all $k\geq 0$, then $a$ does not have the cycle property.
\end{itemize}
\end{thm}

\begin{proof} The proof follows the proof of Theorem \ref{b-cycle} with minimal modifications. In particular, the relation 
(\ref{STS}) should be replaced by relation
\begin{equation}\label{STS*}
(STS)\hat f^nT^mST=(T^{-2}S)\hat f^p.
\end{equation}
For (iii), since $\hat  f^nT^mSTa=\frac{b}{b+1}$, on the lower side we have $TS f^nT^mSTa=Sb$, and on the upper side, using (\ref{STS*}), $(T^{-2}S)\hat f^pb=b$. As in the proof of Theorem \ref{b-cycle}, $\hat f^{p+1}=T^{-3}S$, so $(T^{-3}S)\hat f^pb=T^{-1}b$. Therefore $a$ has (strong or weak) cycle property if and only if $b$ does.
\end{proof}

Let us now
%\foot{SK: added ``now"}
describe the situation when $a\le -1$. %\foot{IU: new part}

\begin{thm}\label{ab-cycle}%\foot{IU: changed a bit}
Let $(a,b)\in\P$ with $0<b\leq -a$ and $a\le -1$. Then $a$ and $b$ satisfy the cycle property.
\end{thm}
\begin{proof}
It is easy to see that $a=-1$ has the degenerate weak cycle:
\begin{equation}\label{a=-1}
\begin{diagram}
        &             & \boxed{1}       \\
        &\ruTo{S}&                                        &\rdTo{T^{-1}}                 \\
\boxed{a=-1} & & \rTo{T}       &   & \boxed{0} \\
     \end{diagram}
\end{equation}
while $a<-1$ satisfies the following strong cycle relation:
\begin{equation}\label{acycle}
\begin{diagram}
              &             & \boxed{-\frac{1}{a}}   &\rTo{T^{-1}}  & \boxed{-\frac{1}{a} - 1} & \rTo{S} &\boxed{\frac{a}{a+1}}\\
              &\ruTo{S} &                               &                  &                                 &            &                                  &\rdTo{T^{-1}}               \\
\boxed{a}&             &                               &                  &                                 &             &             &       & \boxed{c_a=-\frac{1}{a+1}} \\
             &\rdTo{T}  &                                &                 &                                 &             &             &      \ruTo(6,2){S}                            \\
             &             &\boxed{a+1}       
\end{diagram}
\end{equation}
In order to study the orbits of $b$, let $m\ge 0$ such that $a\le T^mSb<a+1$. If $m=0$, then $Sb=a$ (since $Sb\le a$), and the cycle of $b$ is identical to the one described by \eqref{a=-1}. If $m\ge 1$, then one can use relation \eqref{stand} to construct the $b$-cycle. More precisely, if $a<T^mSb<a+1$, then we have:
\begin{equation}
\begin{diagram}
         &        & \boxed{b-1}       &\rTo{S}     & \boxed{-\frac 1{b-1}}   &\rTo{(ST^{-2})^{m-1}}& \boxed{1+\frac{b}{1-mb}}\\
         &\ruTo{T^{-1}}&                   &            &                         &                  &                          &\rdTo{\;\;T^{-1}}\\
\boxed{b}&        &                   &            &                         &                  &                          &        & \boxed{c_b=\frac{b}{1-mb}}\\
    &\rdTo{S}&                   &            &                         &   & \ruTo(3,2){S}       &            &        &\\
    &        &\boxed{-\frac 1{b}}&\rTo{T^{m}}& \boxed{-\frac{1-mb}{b}}
\end{diagram}
\end{equation}
If $T^mSb=a$, then it happens again that the lower orbit of $b$ hits a, and then $Ta$, while the upper orbit hits $Sa$. Following now the cycle of $a$ described by \eqref{acycle}, we conclude that $b$ satisfies the strong cycle property.

If $T^mSb=0$, i.e. $b=1/m$, then a minor modification of the above b-cycle gives us the following weak cycle relation:
\begin{equation}
\begin{diagram}
         &        & \boxed{b-1}       &\rTo{S\;\;}     & \boxed{-\frac 1{b-1}}   &\rTo{T^{-1}(ST^{-2})^{m-2}}& \boxed{\frac{b}{1-mb+b}=1}\\
       & \ruTo_{T^{-1}} &                   &            &                         &                  &                          &\rdTo_{\;\;T^{-1}}\\
\boxed{b}&        &                   &            &                         &                  &                          &        & \boxed{c_b=0}\\
    &\rdTo{S\;\;}&                   &            &                         &   & \ruTo(3,2){T}       &            &        &\\
    &        &\boxed{-\frac 1{b}=-m}&\rTo{T^{m-1}}& \boxed{-1}
\end{diagram}
\end{equation}

\end{proof}

The following corollaries are immediate from the proof of Theorems \ref{b-cycle}, \ref{a-cycle}, \ref{ab-cycle}.

\begin{cor} \label{b-norepeats}If $b$ has the cycle property, then the upper side of the $b$-cycle
\[
\{T^{-1}b, fT^{-1}b,\dots ,f^{m_2-1}T^{-1}b\}
\]
and the lower side of the $b$-cycle
\[
\{Sb, fSb,\dots ,f^{k_2-1}Sb\}
\]
do not have repeating values.
\end{cor}
%\begin{cor} The cycle property for $b$ is strong if and only if $c_b\neq 0$
%\end{cor}

\begin{cor} \label{a-norepeats}If $a$ has the cycle property, then the upper side of the $a$-cycle
\[
\{Sa, fSa,\dots ,f^{m_1-1}Sa\}
\]
and the lower side of the $a$-cycle
\[
\{Ta, fTa,\dots ,f^{k_1-1}Ta\}
\]
do not have repeating values.
\end{cor}
%\begin{cor} The cycle property for $a$ is strong if and only if $c_a\neq 0$
%\end{cor}
%\foot{SK: added Corollaries}

%satisfy the conditions of Theorems \ref{a-cycle} and \ref{b-cycle}, and hence
\section{Finiteness condition implies finite rectangular structure}\label{s:5}
In order to state the condition under which the natural extension map $F_{a,b}$ has an attractor with finite rectangular structure mentioned in the Introduction, we follow the split orbits  of $a$ and $b$ 
%as long as they avoid the intervals $\overline I_u$ and $\overline I_\ell$ \foot{IU: I_u,I_l make sense only for $\Delta$}, 
%and  define the {\em truncated orbits} $\LL_a$ and $\U_a$ by
%\foot{IU: included "strong" for case 2}
\[
\LL_a=\begin{cases} \OO_\ell(Ta)&\text{ if $a$ has no cycle property}\\
\text{lower part of $a$-cycle}&\text{ if $a$ has strong cycle property}\\
\text{lower part of $a$-cycle $\cup\{0\}$}&\text{ if $a$ has weak cycle property},
\end{cases}
\]
\[
\U_a=\begin{cases} \OO_u(Sa)&\text{ if $a$ has no cycle property}\\
\text{upper part of $a$-cycle}&\text{ if $a$ has strong cycle property}\\
\text{lower part of $a$-cycle $\cup\{0\}$}&\text{ if $a$ has weak cycle property},
\end{cases}
\]
and,  similarly, $\LL_b$ and $\U_b$ by
\[
\LL_b=\begin{cases} \OO_\ell(Sb)&\text{ if $b$ has no cycle property}\\
\text{lower part of $b$-cycle}&\text{ if $b$ has strong cycle property}\\
\text{lower part of $b$-cycle $\cup\{0\}$}&\text{ if $b$ has weak cycle property},
\end{cases}
\]
\[
\U_b=\begin{cases} \OO_u(T^{-1}b)&\text{ if $b$ has no cycle property}\\
\text{upper part of $b$-cycle}&\text{ if $b$ has strong cycle property}\\
\text{lower part of $b$-cycle $\cup\{0\}$}&\text{ if $b$ has weak cycle property},
\end{cases}
\]

We find it useful to introduce the map $\rho_{a,b}:\bar\R\to\{T,S,T^{-1}\}$
\begin{equation}
\rho_{a,b}(x)=\begin{cases}
T  &\text{ if } x < a\\
S &\text{ if } a\le x<b\\
T^{-1}  &\text{ if } x\ge b
\end{cases}
\end{equation}
in order to write $f_{a,b}(x)=\rho_{a,b}(x)x$ and $F_{a,b}(x,y)=(\rho(y)x,\rho(y)y)$.
\begin{rem}It follows from the above definitions that $\rho(y)=S$ or $T$ if $y\in\LL_a\cup \LL_b$, and $\rho(y)=S$ or $T^{-1}$ if $y\in\U_a \cup\U_b$.%\foot{IU: remark simplified}
\end{rem}
%\foot{SK20: so for $y=0$ still  $\rho=S$ - correct. Then delete both margin noes.}
\begin{defn} We say that the map $f_{a,b}$ satisfies the {\em finiteness condition} if the sets of values in all four truncated orbits $\LL_a,\LL_b,\,\,\U_a,\,\,\U_b$ are finite.
\end{defn}
%The following proposition will be very useful in the further study.
\begin{prop} \label{shift} Suppose that the set $\LL_b$ is finite. Then
\begin{enumerate}
\item either
%\begin{enumerate}
$b$ has the cycle property or 
the upper and lower orbits of $b$ are eventually periodic. 
%\end{enumerate}
\item The finiteness of $\LL_b$ implies the finiteness of $\U_b$.
\end{enumerate}
 Similar statements hold for the sets $\LL_a$, $\U_a$ and $\U_b$ as well.%\foot{IU: CHANGES!}
\end{prop}
\begin{proof} The two properties follow from Theorem \ref{b-cycle} and its proof. If $b$ does not have the cycle property, but its lower orbit is eventually periodic, then one uses Lemma \ref{main-argument} to conclude that the upper orbit of $b$ has to be eventually periodic.
\end{proof}
\begin{rem}
If $b$ has the strong cycle property, then the set $\LL_b$  coincides with the lower side of the $b$-cycle and $\U_b$ coincides with the upper side of the $b$-cycle. If $b$ does not have the cycle property, but the lower and upper orbits of $b$ are eventually periodic then $\LL_b$ and $U_b$ are identified with these orbits accordingly, until the first repeat.
%if $a$ does not have the cycle property and the $(a,b)$-expansions of $Ta$ $Sa$ are periodic.
\end{rem}

%\begin{cor} 
%\end{cor}
\begin{thm} \label{thm:recstructure} Let $(a,b)\in \mathcal P$, $a\ne 0$, $b\ne 0$, and assume that the map $f_{a,b}$
satisfies the finiteness condition. Then there exists a set $A_{a,b}\underset{\neq}\subset\bar\R^2$ with the following properties:
\begin{enumerate}
\item[(A1)] The set $A_{a,b}$  consists of two connected components each having {\em finite rectangular structure},
i.e. bounded by non-decreasing step-functions with a finite number of steps.
\item[(A2)] $F_{a,b}: A_{a,b}\to A_{a,b}$  is a bijection except for some images of the boundary of $A_{a,b}$.
\end{enumerate}
\end{thm}
%\foot{For weak cycle property $c=0$ must be counted in both $\LL_{a,b}$ and $\U_{a,b}$.}
\begin{proof}
(A1) 
We will construct a set $A_{a,b}$ whose upper connected component is bounded by a step-function with values in the set $\U_{a,b}=\U_a\cup\U_b$ that we refer to as  {\em upper levels}), and whose lower connected component is bounded by a step-function with values in the set $\LL_{a,b}=\LL_a\cup\LL_b$ that we refer to as {\em lower levels}. Notice that each level in $\U_a$ and $\U_b$ appears exactly once, but if the same level appears in both sets, we have to count it twice in $\U_{a,b}$. The same remark applies to the lower levels. 

Now let $y_\ell\in\LL_{a,b}$ be the closest $y$-level to $Sb$ with $y_\ell\geq Sb$, and $y_u\in\U_{a,b}$ be the closest $y$-level to $Sa$ with $y_u\leq Sa$. Since each level in $\U_a$ and in $\LL_b$ appears only once,
 if $y_u=Sa$, $y_u$ can only belong to $\U_b$, and if $y_\ell= Sb$, $y_\ell$ can only belong to $\LL_a$. 
 
 We consider the rays $[-\infty,x_b]\times \{b\}$ and $[x_a,\infty]\times \{a\}$, where $x_a$ and $x_b$ are unknown, and
``transport" them (using the special form of the natural extension map $F_{a,b}$) along the
%the corresponding cycles or, in periodic case, along the sets
sets $\LL_b,\,\,\U_b,\,\,\LL_a$ and $\U_a$ respectively until we reach the levels $y_u$ and $y_\ell$ (see Figure \ref{fig-levels}).
Now we set-up a system of two fractional linear equations  by equating the right end of the segment at the level $Sb$ with the left end of the segment at the level $y_\ell$, and, similarly, the left end of the segment at the level $Sa$ and the right end of the level $y_u$. 

\begin{figure}[htb]
\psfrag{8}[c]{\tiny $STa$}
\psfrag{7}[l]{\tiny $ST^{-1}b$}
\psfrag{5}[c]{\tiny $Sa$}
\psfrag{6}[c]{\tiny $Sb$}
\psfrag{3}[c]{\tiny $y_u$}
\psfrag{4}[c]{\tiny $y_\ell$}
\psfrag{1}[c]{\tiny $x_b$}
\psfrag{2}[c]{\tiny $x_a$}
\psfrag{a}[c]{\small $a$}
\psfrag{b}[c]{\small $b$}
\includegraphics[width=2.7in, height=3.1in]{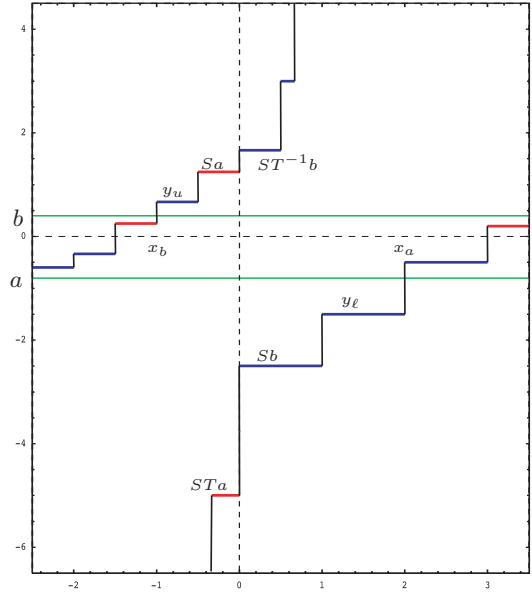}
\caption{Construction of the domain $A_{a,b}$}
\label{fig-levels}
%\label{g8}
\end{figure}
\begin{lem}\label{2-sys} The system of two equations at the consecutive levels $y_u$ and $Sa$, and $y_\ell$ and $Sb$, 
has a unique solution with $x_a\geq 1$ and $x_b\leq -1$. %\foot{SK: this does not work for Gauss; need to have a disclaimer of some kind}
\end{lem}
\begin{proof} 
In what follows, we present the proof assuming that $0<b\le -a<1$. The situation $a\le -1$ is less complex due to the explicit cycle expressions described in Theorem \ref{ab-cycle} and will be discussed at the end.
Let $m_a,m_b$ be positive integers such that $a\le T^{m_a}STa<a+1$ and $a\le T^{m_b}Sb<a+1$. For the general argument we assume that  $m_a,m_b\geq 3$, the cases $m_a$ or $m_b\in \{1,2\}$ being considered separately. %\foot{SK: need to write this part}
The level $y_u$ may belong to $\U_a$ or $\U_b$, and the level $y_\ell$ may belong to $\LL_a$ or $\LL_b$, therefore we need to consider  $4$ possibilities. 
%In what follows we assume that both $a$ and $b$ have cycle property, and  then explain which modifications are needed for the periodic situation.
%\foot{SK: lemma holds for periodic situation as well if we allow the rays $(-\infty,x_b)$ and $(x_a,\infty)$ to have a non-trivial intersection with the boundary of the constructed region $A_{a,b}$}

\noindent{\bf Case 1: $y_u\in\U_a,\,\,y_\ell\in\LL_a$.} Then we have
\[
%\begin{aligned}
Sx_a= T^{-1}S\hat f_-^{n_1}(\infty)\,,\quad
Sx_b=  TS   \hat f_+^{n_2}Tx_a,
%\end{aligned}
\]
where $\hat f_-^{n_1}$ is a product of factors $T^{-i}S$ (that appear on the upper orbit of $a$) with $i=2$ or $3$, and $\hat f_+^{n_2}$ is a product of factors $T^{i}S$ (that appear on the lower orbit of $a$) with $i=m$ or $m+1$. Using (\ref{STS*}) we rewrite the first equation as
\[
%\begin{aligned}
x_a=ST^{-1}S\hat f_-^{n_1}(\infty)=ST^{-1}SST^2STS\hat f_+^{k_1}T^mST(\infty)
=T^{-1}\hat f_+^{k_1}T^mST(\infty)\,.
%\end{aligned}
\]
Since $\hat f_+^{k_1}$ is a product of factors $T^{i}S$ with $i=m$ or $m+1$, $m\geq 3$, we conclude that $Tx_a$ has a finite formal continued fraction expansion starting with $m'\geq 3$, i.e. $Tx_a>2$, and $x_a>1$. 
%\foot{SK20: small change}
Furthermore, from the second equation
\[
x_b=  STS   \hat f_+^{n_2}Tx_a,
\]
hence $ \hat f_+^{n_2}Tx_a$ has a finite formal continued fraction expansion starting with $m'\geq 3$, i.e. $ \hat f_+^{n_2}Tx_a>2$, and $x_b<-2$.

\noindent{\bf Case 2: $y_u\in\U_a,\,\,y_\ell\in\LL_b$.} Then
\[
%\begin{aligned}
Sx_a= T^{-1}S\hat f_-^{n_1}(\infty)\,,\quad
Sx_b=  TS   \hat f_+^{n_2}(-\infty)\,.
%\end{aligned}
\]
Like in Case 1 we see that $x_a>1$, and
\[
x_b=  STS   \hat f_+^{n_2}(-\infty)<-2,
\]
since $ \hat f_+^{n_2}(-\infty)$ has a formal continued fraction expansion starting with $m'\geq 3$, and therefore is $>2$.

\noindent{\bf Case 3: $y_u\in\U_b,\,\,y_\ell\in\LL_a$.} Then
\[
%\begin{aligned}
Sx_a= T^{-1}S\hat f_-^{n_1}T^{-1}x_b\,,\quad
Sx_b=  TS   \hat f_+^{n_2}Tx_a\,.
%\end{aligned}
\]
Using (\ref{STS}) %\foot{SK23: ref changed!}we rewrite the first equation as
\[
x_a= ST^{-1}S\hat f_-^{n_1}T^{-1}x_b=ST^{-1}SST^2STS\hat f^{k_2}T^mSx_b,
\]
and using the second equation and simplifying, we obtain
\[
Tx_a=\hat f^{k_2}T^mSSTS\hat f_+^{n_2}(Tx_a)=\hat f^{k_2}T^{m+1}S\hat f_+^{n_2}(Tx_a).
\]
Since all its factors are of the form $T^iS$ with $i\geq 3$,  the matrix $\hat f^{k_2}T^{m+1}S\hat f_+^{n_2}$ is hyperbolic and its attracting fixed point $Tx_a$ has periodic formal continued fraction expansion starting with $m'\geq 3$ (see Theorem 3.1 of \cite{KU1}), hence $x_a>1$. Finally, as in Case 1,
\[
x_b=STS\hat f_+^{n_2}Tx_a<-2
\]
since $\hat f_+^{n_2}Tx_a$ has formal continued fraction expansion with $m'\geq 3$, hence $>2$.

\noindent{\bf Case 4: $y_u\in\U_b,\,\,y_\ell\in\LL_b$.} Then
\[
%\begin{aligned}
Sx_a= T^{-1}S\hat f_-^{n_1}T^{-1}x_b\,,\quad
Sx_b=  TS   \hat f_+^{n_2}S(-\infty)\,.
%\end{aligned}
\]
From the second equation we obtain 
\[
x_b=STS\hat f_+^{n_2}S(-\infty)<-2
\]
since $\hat f_+^{n_2}S(-\infty)$ has formal continued fraction expansion with $m'\geq 3$, hence $>2$.
Finally,
\[
x_a= ST^{-1}S\hat f_-^{n_1}T^{-1}x_b=T^{-1}\hat f^{k_2}T^{m+1}S\hat f_+^{n_2}S(-\infty),
\]
hence
\[
Tx_a=\hat f^{k_2}T^{m+1}S\hat f_+^{n_2}S(-\infty)>2
\]
since it has formal continued fraction expansion with $m'\geq 3$, therefore $x_a>1$.

Now we analyze the particular situations when $m_a$ or $m_b\in\{1,2\}$, using the explicit cycle descriptions that exist for these situations as described by Theorems \ref{b-cycle} and \ref{a-cycle}.
%\begin{itemize}

(i)  If $m_a=m_b=1$, then relation \eqref{m1} for the $b$-cycle and a similar one for the $a$-cycle shows that $y_\ell=-\frac1{b}+1$ and $y_u=-\frac1{a}-1$, therefore $x_a=1$ and $x_b=-1$.%\foot{SK20: attn: references corrected}
%\foot{SK20: values of $y_\ell$ and $y_u$ corrected}
%\foot{IU: needs more work}
\medskip

(ii) If $m_a=1$, $m_b=2$, following the explicit cycles given by \eqref{m=2} we obtain $y_\ell=-1/b +1$, and $y_u=-1/(b-1)-1$, therefore $x_a=2$, $x_b=-1$.

\medskip

(iii) If $m_a=1$, $m_b\ge 3$, using the cycle structure in Theorem \ref{b-cycle} we obtain $y_\ell=1/b+1$ and $y_u=T^{-1}(ST^{-2})^{m_b-2}ST^{-1}b$, therefore, $x_a=m_b$, and $x_b=-1$.

\medskip

(iv) If $m_a=2$, $m_b=2$, using the cycle structure in Theorems \ref{b-cycle} and \ref{a-cycle} we obtain $y_\ell=-\frac1{a+1}+1$ and $y_u=-\frac1{b-1}-1$, and a calculation in this particular case, like in Lemma \ref{2-sys}, Case 3 implies that  $x_a>1$ and $x_b<-1$.

\medskip

(v) if $m_a=2$, $m_b>2$, an analysis of the four cases above for this particular situation (with an explicit cycle relation for $a$) yields $x_a\ge 1$ and $x_b\le -1$. Indeed, in Case 1, we have $y_u=-1/a-1$, hence $x_a=1$ and $x_b=-2$. In Case 2, we get $x_a=1$ and $x_b<-2$. Cases 3 and 4 are treated similarly.
\end{proof}
%\foot{SK24: added}
Now, since $x_a$ and $x_b$ are uniquely determined, by ``transporting" the rays $[-\infty, x_b]\times\{b\}$ and $[x_a,\infty]\times\{a\}$ along the sets $\LL_b,\,\,\U_b,\,\,\LL_a$ and $\U_a$ we obtain the $x$-coordinates of the right and left end of the segments on each level.
\begin{defn} We say that two consecutive levels $y_1\leq y_2$ of $\LL_{a,b}$, respectively, $\U_{a,b}$, are called {\em connected by a vertical segment} (we will refer to this as {\em connected}) if the $x$-coordinate of the right end point of the horizontal segment on the level $y_1$ is equal to the the $x$-coordinate of the left end point of the horizontal segment on the level $y_2$.
\end{defn}
We will prove that  all levels of $\LL_{a,b}$ and all levels of $\U_{a,b}$ are connected.%\foot{SK25: change}
%, i.e. that the right end of the segment at a certain level is equal to the left end of the segment on the next level.
 We first look at the levels in $\LL_{a,b}$. By Lemma \ref{2-sys} the levels $y_u$ and $Sa$, and the levels $Sb$ and $y_\ell$ are connected. %Some levels of $\LL_{a,b}$ are automatically connected.
\begin{lem} \label{next} The levels $Sb\in\LL_b$ and $STa\in\LL_a$ are two consecutive levels of  $\LL_{a,b}$ connected by a vertical segment at $x=0$.
The levels $Sa\in\U_a$ and $ST^{-1}b\in\U_b$ are two consecutive levels of $\U_{a,b}$  connected by a vertical segment at $x=0$.
\end{lem}
\begin{proof} Suppose  there is $y\in \LL_{a,b}$ such that $STa\leq y\leq Sb$. Then $y\in \LL_a$ or $\LL_b$. In either case,
 since by Lemmas \ref{a-norepeats} and  \ref{b-norepeats} the truncated orbits $\LL_a,\LL_b$ do not have repeated values, neither $STa=y$ nor $y=Sb$ is possible. Thus the only case we need to consider is
%
%Suppose there is $y\in \LL_{a,b}$ such that
\[
STa<y<Sb.
\]
%But if $y\in\LL_a$, t
Then, either $y=Sy'$ for some $y'\in\LL_{a,b}$ ($0<y'\le a+1$) or $y=Ty''$ for some $y''\in\LL_{a,b}$. These would imply that either $y'>Ta$, which is impossible, or $Ty''<Sb$, i.e. $y''<T^{-1}Sb$, which is also impossible (if $y''<T^{-1}Sb$ then $y=Ty''$ must be the end of the $a$-cycle, by Theorem \ref{a-cycle}). The $x$-coordinate of the right end point of the segment at the level $STa$ and of the left end point of the segment at the level $Sb$ is equal to $0$.%\foot{SK25:added} 
The second part of the proof is similar.%\foot{IU: added reference}
\end{proof}

The following proposition will be used later in the proof.
\begin{prop} \label{tech} Suppose that the set $\LL_{a,b}$ is finite and $y\in \LL_{a,b}$ with $y>STa$. 
\begin{enumerate}
\item If $y\in\LL_a$, then  %either $a$ has no cycle property, the expansion of $STa$ is periodic and $y=STa$,or 
there exists $n_0>0$  such that
$\rho(f^{n}y)=\rho(f^n STa)$
 for all $0<n<n_0$ and  $\rho(f^{n_0}y)\neq \rho(f^{n_0}STa)$, or $f^{n_0}y=0$;
 %$f^{n_0}y$ is the end of the $a$-cycle, if we are in the cycle situation, or
 \item If $y\in \LL_b$, then $y>Sb$, and there exists $n_0>0$  such that
$\rho(f^{n}y)=\rho(f^n Sb)$
 for all $n<n_0$ and $\rho(f^{n_0}y)\neq \rho(f^{n_0}Sb)$, or $f^{n_0}y=0$.
 % or $f^{n_0}y$ is the end of the $a$-cycle, if we are in the cycle situation.
 \end{enumerate}
 \end{prop}
 \begin{proof}
Suppose that $y\in\LL_a$ and $a$ satisfies the cycle property. It follows that such an $n_0$ exists or $f^{n_0}y$
 is the end of the $a$-cycle. We will show that the latter is possible only if $f^{n_0}y=0$, i.e. it is the end of a weak cycle.
%the possibility we would like to exclude. We first give the argument for the $a$-cycle.
 Suppose $f^{n_0}y$ is the end of the $a$-cycle.
Then if
\[
\rho(f^{n_0-1}y)=\rho(f^{n_0-1}STa)=S,
\]
we must have
 $f^{n_0-1}y<0$ since otherwise the cycle would not stop at $S$, but $f^{n_0-1}(STa)>0$ since for $STa$ we have not reached the end of the cycle. This contradicts the monotonicity of $f^{n_0-1}$ and the original assumption $y>STa$, thus is impossible. The other possibility is
 \[
 \rho(f^{n_0-1}y)=\rho(f^{n_0-1}STa)=T.
 \]
 But this either implies that $f^{n_0-1}y<T^{-1}Sb$, and by monotonicity of $f^{n_0-1}$, $f^{n_0-1}(STa)<f^{n_0-1}y<T^{-1}Sb$, which implies that we have reached the end of the cycle of $STa$ as well, a contradiction, or, 
$f^{n_0}y=0$, i.e. it is the end of a weak cycle.

Now suppose $y\in\LL_b$. Then by Lemma \ref{next} $y\geq Sb$, but
since each level in $\LL_b$ appears only once, we must have but $y>Sb$. Now the argument that $f^{n_0}y$ cannot be the end of the $b$-cycle is exactly the same as for the $a$-cycle.
 
 In the periodic case, let us assume that no such $n_0$ exists. Then, in case (1) the $(a,b)$-expansions of $STa$ and $y$, which is the lower part of the former, are the same, i.e. $(a,b)$-expansions of $STa$ is invariant by a left shift.
In case (2), we have seen already that we must have $y>Sb$.%\foot{SK21: deleted ``but"}
%by Lemma \ref{next} $y\geq Sb$, but
%since each level in $\LL_b$ appears only once, we must have but $y>Sb$.
Then the $(a,b)$-expansions of $Sb$ and $y$, which is the lower part of the former, are the same, i.e. $(a,b)$-expansions of $Sb$ is invariant by a left shift.
 The proof that this is impossible 
is based on the following simple observation:  if $\sigma= (a_1,a_2,\dots , a_k, \overline{a_{k+1},a_{k+2},\dots a_{k+n}})$ is an eventually periodic symbolic sequence with the minimal period $n$  and invariant under a left shift by $m$, then $\sigma$ is purely periodic and $m$ is a multiple of $n$. 

By the uniqueness property of $(a,b)$-expansions, this would imply that
$y=STa$ or $y=Sb$, a contradiction.
\end{proof}
%\foot{make sure this is accurate for m=1 or 2}
Let $y_b^-,y_b^+\in\U_{a,b}$ be two consecutive levels with $y_b^-\leq b< y_b^+$, and  $y_a^-,y_a^+\in\LL_{a,b}$ be two consecutive levels with $y_a^-< a\leq y_a^+$.%\foot{SK23: IMPORTANT: inequalities changed}
\begin{lem}\label{a-connected}
%\begin{enumerate}
There is always one level connected with level $a+1$, and
the levels $y_a^-$ and $y_a^+$ are connected by the vertical segment at $x_a$.
 %\foot{SK24: change}
%\end{enumerate}
\end{lem}
\begin{proof}%\foot{SK25: change}
%Notice that the left end of the segment at the level $Sb$ and the right end of the segment at the level
% $STa$ are equal to $0$, hence 
 By Lemmas \ref{2-sys} and \ref{next}, we know that three consecutive levels $STa\leq Sb\leq y_\ell$ are connected. Moreover, their images remain connected under the same transformations in $\sz$.
Since each level in $\U_a$ and in $\LL_a$ appears only once, at least one of the two inequalities must be strict, i.e. if $STa= Sb$, then $STa=Sb< y_\ell$, and if $Sb= y_\ell$, then $STa< Sb= y_\ell$.

First we prove that
$y_\ell<TSb$. Suppose $y_\ell\geq TSb$. Its pre-image must be $y'_\ell=T^{-1}y_\ell$ since for any $y,\,0<y<Ta$,
$Sy<STa\leq Sb<TSb$, and
%The first is by the choice of $y_\ell$, and the second since otherwise 
%for its pre-image $y'_\ell=y_\ell-1$ 
we would have
$Sb\leq y'_\ell<y_\ell$ that contradicts the assumption that $y_\ell$ is {the next level above} $Sb$. Therefore, if the first digit in the $(a,b)$-expansion of $Sb$ is $-m$,
then the first digit of $y_\ell$ is $-(m-1)$ or $-m$. In the first case, the three levels%\foot{IU: checked}
\[
T^{m-1}Sb< a\le T^{m-1}y_{\ell}
\]
are connected and satisfy  $T^{m-1}Sb=y^-_a, T^{m-1}y_{\ell}=y^+_a$.
Therefore, the levels $T^{m}Sb$ and $a+1$ are connected. 

For the second case, we know that $Sb\leq y_\ell$ and
\[
a\le T^{m}Sb\leq T^{m}y_{\ell}< a+1.
\]
If $Sb=y_\ell$, then $y_\ell\in \LL_a$, and $STa<y_\ell$.
%and since $STa<y_\ell$, the $f$-orbit of $y_\ell$ is the lower part of the $f$-orbit of $STa$ in $\LL_a$.
If $Sb<y_\ell$, then $y_\ell\in  \LL_b$, or $y_\ell\in  \LL_a$ and  $STa<y_\ell$.
%If   $y_\ell$ belongs to $\LL_a$, since $STa<y_\ell$, the $f$-orbit of $y_\ell$ is the lower part of the $f$-orbit of 
%$STa$ in $\LL_a$, if $y_\ell$ belongs to $Sb$,  the $f$-orbit of $y_\ell$ is the lower part of the $f$-orbit of $Sb$ in $\LL_b$.

Let us assume that $y_\ell$ belongs to $\LL_a$.
%We should keep in mind two possible possibilities: (a) $a$ has a cycle property or (b) the $(a,b)$-expansion of $a$ is eventually periodic, as was explained above. 
Since $STa<y_\ell$, by Proposition \ref{tech}, there are two possibilities:
\begin{enumerate}
\item $f^{n_0}y_\ell$ is the end of a weak cycle. 
\item There exists 
$n_0$  such that
$\rho(f^{n}y_\ell)=\rho(f^n STa)$
 for all $n<n_0$, and  $\rho(f^{n_0}y_\ell)\neq \rho(f^{n_0}STa)$.
 \end{enumerate}
 In the first case, we have $f^{n_0}STa=y_a^-$ and $f^{n_0}Sb=y_a^+$, or $f^{n_0}Sb=y_a^-$ and $ f^{n_0}y_{\ell}=y_a^+$. Therefore, either
$f^{n_0+1}STa$ or $f^{n_0+1}Sb$ is connected with level $a+1$.
In the second case,
we notice that
\[
\rho(f^{n_0-1}y_\ell)=\rho(f^{n_0-1}STa)=T
\]
otherwise,  $\rho(f^{n_0-1}y_\ell)=\rho(f^{n_0-1}STa)=S$
%, then, since
%after $S$ we necessarily have $T$ in the set $\LL_a$, we 
would imply
 \[
 \rho(f^{n_0}y_\ell)=\rho(f^{n_0}STa)=T
 \]
in contradiction with the choice of $n_0$.
Further,
there are two possibilities:
\[
\text{(i) } \rho(f^{n_0}STa)=S,\,\,\rho(f^{n_0}y_\ell)=T\,,\quad \text{(ii) } \rho(f^{n_0}STa)=T,\,\,\rho(f^{n_0}y_\ell)=S.
\]
In case (i) we obtain %\foot{IU: checked}
\[
f^{n_0}y_\ell<a\le f^{n_0}STa
\]
which contradicts the monotonicity of $f$ and the original assumption $y_\ell>STa$. Thus the only possibility is
%the second case when
\[
f^{n_0}y_\ell \ge a>f^{n_0}STa.
\]
By using the monotonicity of  $f^{n_0}$ we have
\[
f^{n_0}y_\ell>f^{n_0}Sb>f^{n_0}STa
\]
and conclude that  $f^{n_0}STa=y_a^-$ and $f^{n_0}Sb=y_a^+$, or $f^{n_0}Sb=y_a^-$ and $ f^{n_0}y_{\ell}=y_a^+$. Therefore, either
$f^{n_0+1}STa$ or $f^{n_0+1}Sb$ is connected with level $a+1$.
The case when $y_\ell$ belongs to $\LL_b$ is very similar, and in this case $f^{n_0}Sb=y_a^-$, $ f^{n_0}y_{\ell}=y_a^+$, and  $f^{n_0+1}Sb$ is connected with $a+1$. By construction, in both cases  the common $x$-coordinate of the end points is equal to $x_a$. %\foot{SK24: change}
\end{proof}
After an application of $S$ the level connected with $a+1$ will be connected with $STa$, and now, instead of $3$ connected levels $STa\leq Sb\leq y_\ell$ (with at least one strict inequality) we have at least $4$ connected levels $y'\leq STa\leq Sb\leq y_\ell$ (with no more than two equalities in a row).

The process continues with a growing number of connected levels, the highest being $a+1$. Since on each step we cannot have more than two equalities in a row, the number of distinct levels in this sequence will also increase. Therefore, we obtain a sequence of connected levels
 \begin{equation}\label{lowersnake}
 a+1\geq y_1\geq \dots\geq y_{s}>\frac{b}{b+1}\geq y_{s+1}.
 \end{equation}%\foot{IU: fixed indices}
It is evident from the construction that 
there are no unaccounted  levels $y\in \LL_{a,b}$, $a+1\geq y\geq y_{s+1}$.%\foot{SK25: inequalities corrected}

Now we prove a similar result for $\U_{a,b}$.
\begin{lem}\label{b-connected}
%\begin{enumerate}
There is always one level connected with level $b-1$, and
the levels $y_b^-$ and $y_b^+$ are connected by a vertical segment at $x_b$.

%\end{enumerate}
\end{lem}
 \begin{proof}
% The proof that the levels of $\U_{a,b}$ are joint follows the same scheme.
By Lemmas \ref{2-sys} and \ref{next} %\foot{SK25: changed refs} 
we know that the three consecutive levels $y_u\leq Sa\leq ST^{-1}b$ are connected.
It is easy to see that the first digit in $(a,b)$-expansion of $ST^{-1}b$ is $2$,
and the first digit in $(a,b)$-expansion of $Sa$ is either $1$ or $2$. Therefore, the first digit in $(a,b)$-expansion of $y_u$ is either $1$ or $2$.
In the first case either %\foot{IU: checked}
\[
T^{-1}Sa< b\leq T^{-1}ST^{-1}b
\]
or
\[
T^{-1}y_u<b\leq T^{-1}Sa
\]
are the connected levels. Therefore either $T^{-1}Sa=y_b^-$ and $T^{-1}ST^{-1}b=y_b^+$, or $T^{-1}y_u=y_b^-$ and $T^{-1}Sa=y_b^+$ are connected. So either $T^{-2}ST^{-1}b$ or $T^{-2}Sa$ is connected with level $b-1$. 

In the second case, we know that $y_u\leq Sa$ and 
\[
b-1\leq T^{-2}y_u\leq T^{-2}Sa<b.
\]
If $y_u=Sa$, $y_u$ must belong to $\U_b$, in which case $y_u<ST^{-1}b$. If $y_u<Sa$, then $y_u\in \U_a$, or  $y_u\in \U_b$ and $y_u<ST^{-1}b$.
%

%
%notice that $y_u$ belongs to $\U_a$ or $\U_b$, and therefore the $f$-orbit of $y_u$
%is the lower part of the $f$-orbit of $Sa$ in $\U_a$ or $ST^{-1}b$ in $\U_b$. 
Let us assume that $y_u$ belongs to $\U_b$. Since $y_u<ST^{-1}b$, by
%We again have two possible possibilities to consider: (a) $b$ has a cycle property or (b) the $(a,b)$-expansion of $b$ is eventually periodic. 
%As in the proof of Lemma \ref{a-joint}, 
Proposition \ref{tech} there are two possibilities:
\begin{enumerate}
\item $f^{n_0}y_u$ is the end of a weak cycle,
\item there
exists $n_0$ such that $\rho(f^{n}y_u)=\rho(f^nST^{-1}b)$
for all $n<n_0$, and  $\rho(f^{n_0}y_u)\neq \rho(f^{n_0}ST^{-1}b)$.
\end{enumerate}
In the first case,  either $f^{n_0}ST^{-1}b=y_b^+$ and $f^{n_0}Sa=y_b^-$, or $f^{n_0}Sa=y_b^+$ and $ f^{n_0}y_u=y_b^-$, so either
$f^{n_0+1}ST^{-1}b$ or $f^{n_0+1}Sa$ is connected with level $b-1$.
In the second case, we first notice that
\[
\rho(f^{n_0-1}y_u)=\rho(f^{n_0-1}ST^{-1}b)=T^{-1}
\]
 since if we had $\rho(f^{n_0-1}y_u)=\rho(f^{n_0-1}ST^{-1}b)=S$, then
 %, since
%after $S$ we necessarily have $T^{-1}$ in this cycle, 
we would have
 \[
 \rho(f^{n_0}y_u)=\rho(f^{n_0}ST^{-1}b)=T^{-1}
 \]
in contradiction with the choice of $n_0$.
%Then
%\foot{SK: do we need this line?}
%\[
%\rho(f^{n_0-2}y_u)=\rho(f^{n_0-2}ST^{-1}b)=T^{-1}
% \]
%since a single
%$T^{-1}$ may only appear at the end of the cycle.
Further,
there are two possibilities:
\[
\text{(i) }  \rho(f^{n_0}ST^{-1}b)=S,\,\,\rho(f^{n_0}y_u)=T^{-1},\;
\text{(ii) } \rho(f^{n_0}ST^{-1}b)=T^{-1},\,\,\rho(f^{n_0}y_u)=S.
\]
In the first case we obtain
\[
f^{n_0}y_u>b>f^{n_0}ST^{-1}b
\]
which contradicts the monotonicity of $f^{n_0}$ and the original assumption $y_u<ST^{-1}b$. Thus the only possibility is 
\[
f^{n_0}y_u<b<f^{n_0}ST^{-1}b.
\]
By monotonicity of  $f^{n_0}$ we have
\[
f^{n_0}y_u<f^{n_0}Sa<f^{n_0}ST^{-1}b. 
\]
Therefore either $f^{n_0}y_u=y_b^-$ and $f^{n_0}Sa=y_b^+$,
 or $f^{n_0}Sa=y_b^-$ and $f^{n_0}ST^{-1}b=y_b^+$ are connected. So either $T^{-1}f^{n_0}ST^{-1}b$ or $T^{-1}f^{n_0}Sa$ is connected with level $b-1$. The case when $y_u$ belongs to the $a$-cycle is very similar, and in this case $f^{n_0}y_u=y_b^-$ and $f^{n_0}Sa=y_b^+$ and $T^{-1}f^{n_0}Sa$ is connected with level $b-1$. By construction, in both cases the common $x$-coordinate of the end points of the segments at the levels $y_b^-$ and $y_b^+$ is $x_b$.%\foot{SK24: change}
\end{proof}
%In all cases we obtain that
%\begin{enumerate}
%\item the levels $y_b^+$ and $y_b^-$ are joint;
%\item at least one level is joint with $b-1$.
%\end{enumerate}
After an application of $S$ the levels (2) will be connected with $ST^{-1}b$, and now, instead of $3$ connected levels $y_u\leq Sa\leq ST^{-1}b$ we have at least $4$ connected levels $y_u\leq Sa\leq ST^{-1}b\leq y''$.

The process continues with a growing number of connected levels, the lowest being $b-1$. Also the number of distinct levels will increase, and we obtain a sequence of connected levels
%\foot{IU: fixed}
 \begin{equation}\label{uppersnake}
 b-1\leq \bar y_1\leq \dots\leq \bar y_t<\frac{a}{1-a}\leq \bar y_{t+1}.
 \end{equation}
It is evident from the construction that 
there are no unaccounted levels $y\in \U_{a,b}$, $b-1\leq y\leq \bar y_{t+1}$.%%
%\foot{SK25: change}

Now we complete the proof that all levels of $\LL_{a,b}$ are connected. For that it is sufficient to find a sequence of  connected levels  with the distance between the highest and the lowest level $\geq 1$ and the lowest level $\geq T^{-1}Sb$. This is because the set of levels in $y\in\LL_{a,b}$ satisfying $T^{-1}Sb\leq y\leq a+1$ is periodic with period $1$, and each $y\in \LL_{a,b}$ uniquely determines a horizontal segment on level $y$, as was explained just before Lemma \ref{next}.
% starting with  the sequence of the connected levels (\ref{lowersnake}) and ``extending" it using the transformations $S$ and $T$. 

If $y_{s+1} \leq a$, then 
%$y_{s+2}=y_a^-$ and $y_{s+1}=y_a^+$, and all levels in $\LL_{a,b}$ are connected. If $y_{s+2} =a$,
%$y_{s+2}=y_a^+$ that is connected with $y_a^-$, thus 
all levels in $\LL_{a,b}$ are connected.
Suppose now that $y_{s+1}> a$. If
$y_{s+1}=y_a^+$, then, since $y_a^+$ is already connected with $y_a^-$, all levels of $\LL_{a,b}$ are connected.
Now assume that $y_{s+1}>y_a^+$. Then either
\[  y_{s+1}=\frac{b}{b+1}  \quad  \text{ or } \quad  y_{s+1}<\frac{b}{b+1}.
\]
%\end{enumerate}
In the first case either $TSy_{s+1}=y_\ell=Sb$ (this can only happen if $y_{s+1}\in\LL_a$), or $TSy_{s}>Sb$ is the next level above $Sb$, and hence  $TSy_{s}=y_\ell$. In either case
$Sy_{s+1}\leq Sy_{s}\leq \cdots \leq STa\leq Sb=TSy_{s+1}$ are the connected levels with the distance between the lowest and the highest equal to $1$,
%Since the levels $y_{s+2}$ and $y_{s+1}$ are connected, either $TSy_{s+2}=y_\ell=Sb$ is connected to the next level
%(this can only happen if $y_{s+2}\in\LL_a$), 
%or 
%$TSy_{s+1}>Sb$ is the next level above $Sb$, and hence  $TSy_{s+1}=y_\ell$, and we again have 
thus we conclude that all levels of $\LL_{a,b}$ are connected.

In the second case,  %let us first assume that $0<y_{s+2}<\frac{b}{b+1}$. If $y_a^+\neq 0$,
the two levels $y_a^+<y_{s+1}$ will produce the ends of the cycles (one of them can be weak if one of $y_a^+$ or $y_{s+1}$ is equal to $0$). 
%$TSy_{s+2}=c_1$, where $c_1=c_a$ or $c_b$, the end of one of the cycles. 
By the cycle property (Proposition \ref{main-argument}(ii)), there exists a level $z\in\U_{a,b}$, $\frac{a}{1-a}< z<b$ such that $z=(STS)y_{s+1}$. We claim that $z=y_b^-$. Suppose not, and $z<y_b^-$. Then $y_b^-$ gives rise to the second cycle, and again by the cycle property, there exists $y\in\LL_{a,b}$, $y<\frac{b}{b+1}$, such that $y_b^-=STSy$. Since $STS(z)=-\frac{z}{z-1}$
is monotone increasing for $z<1$, we conclude that 
%and therefore $y<\frac{b}{b+1}$. Thus we have
%\[
%TSy=Sy_b^->Sz=TSy_{s+2},
%\]
%which implies 
$y>y_{s+1}$ in contradiction with (\ref{lowersnake}). Thus $y_b^-=(STS)y_{s+1}$. 
%We notice that since $y_{s+2}<\frac{b}{b+1}$, $y_b^-<b$ (and  $b<y_b^+$), the right end of the segment at the level $y_b^-$ is equal to $x_b$
Then $TSy_{s+1}=Sy_b^-$ which implies that the right end of the segment at the level $Sy_b^-$, which is equal to the right  end of the segment at the level $Sb$, is equal to the right end of the segment at the level $TSy_{s+1}$ (notice that this level may belong to $\LL_{a,b}$, $\U_{a,b}$ or be at infinity if $y_{s+1}=0$). Since $y_{s}$ and $y_{s+1}$ were connected, the left end of the segment at the level $TSy_{s}$ is equal to the right end of the segment at the level $TSy_{s+1}$ even though they may belong to the boundaries of different connected components. Since $TSy_{s}\in \LL_{a,b}$,  we conclude that the segment at the level $TSy_{s}$ is adjacent to the segment at the level $Sb$, i.e. $TSy_{s}=y_\ell$. Thus  $Sy_{s}\leq Sy_{s-1}\leq \cdots \leq STa\leq Sb\leq TSy_{s}$ are the connected levels with the distance between the lowest and the highest equal to $1$, and therefore
all levels in $\LL_{a,b}$ are also connected. The proof for $\U_{a,b}$  follows exactly the same lines.

\medskip

(A2) In order to prove the bijectivity of the map $F$ on $A_{a,b}$ we write it as a union of the upper and lower connected components, $A_{a,b}=A_{a,b}^u\cup A_{a,b}^\ell$, and subdivide  each component into $3$ pieces: $A_{a,b}^u=\cup_{i=1}^3U_i$, and $A_{a,b}^\ell=\cup_{i=1}^3L_i$, where

\[\begin{aligned}
U_1=&\{(x,y)\in A_{a,b}^u\,:\, y\ge b\}\\
U_2=&\{(x,y)\in A_{a,b}^u\,:\, b-1\le y\le 0\}\\
U_3=&\{(x,y)\in A_{a,b}^u\,:\, 0\le y\le b\}\\
L_1=&\{(x,y)\in A_{a,b}^\ell\,:\, y\le a\}\\
L_2=&\{(x,y)\in A_{a,b}^\ell\,:\, 0\le y\le a+1\}\\
L_3=&\{(x,y)\in A_{a,b}^\ell\,:\, a\le y\le 0\},
\end{aligned}
\]
Now let
\[
U'_1=T^{-1}(U_1),\,\,
U'_2=S(U_2),\,\,
U'_3=S(U_3),\,\,
L'_1=T(L_1),\,\,
L'_2=S(L_2),\,\,
L'_3=S(L_3)
\]
be their images under the transformation $F$ (see Figure \ref{fig-bij}).
 
\begin{figure}[htb]
\psfrag{a}[l]{\tiny $a$}
\psfrag{b}[l]{\tiny $b$}
\psfrag{x}[c]{\tiny $x_b$}
\psfrag{y}[c]{\tiny $x_a$}
\psfrag{M}[c]{\tiny $F_{a,b}$}
\psfrag{A}[c]{\tiny $\mathbf{U_1}$}
\psfrag{C}[c]{\tiny $\mathbf{U_2}$}
\psfrag{E}[c]{\tiny $\mathbf{U_3}$}
\psfrag{B}[l]{\tiny $\mathbf{U'_1}$}
\psfrag{D}[c]{\tiny $\mathbf{U'_2}$}
\psfrag{F}[c]{\tiny $\mathbf{U'_3}$}
\psfrag{G}[c]{\tiny $\mathbf{L_2}$}
\psfrag{I}[c]{\tiny $\mathbf{L_3}$}
\psfrag{K}[c]{\tiny $\mathbf{L_3}$}
\psfrag{H}[c]{\tiny $\mathbf{L'_2}$}
\psfrag{J}[c]{\tiny $\mathbf{L'_3}$}
\psfrag{N}[c]{\tiny $\mathbf{L'_1}$}
\psfrag{L}[c]{\tiny $\mathbf{L_1}$}
\includegraphics[scale=0.5]{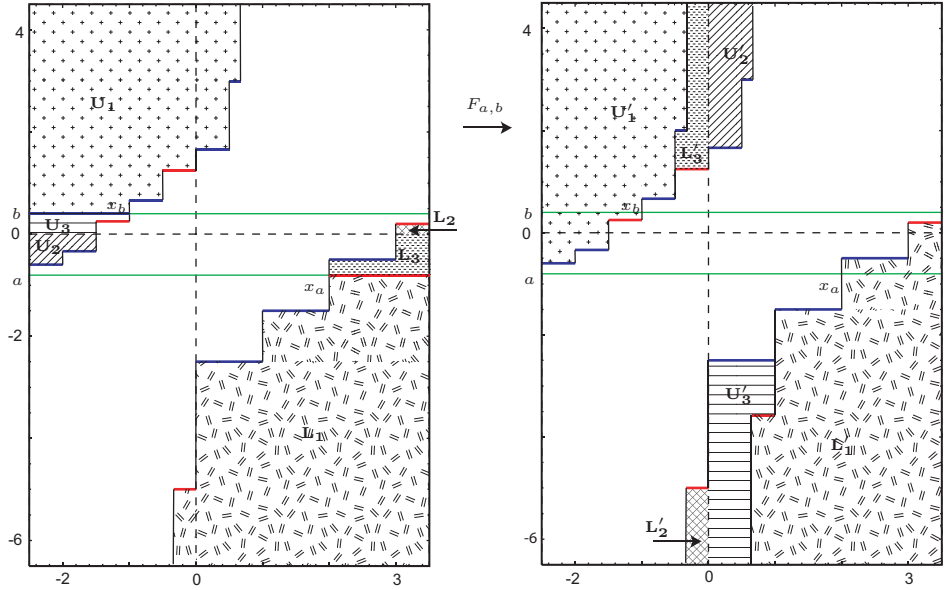}
\caption{Bijectivity of the map $F_{a,b}$}
%\center{The map $F$ is injective}
\label{fig-bij}
\end{figure}

Since the set $A_{a,b}$ is bounded by step-functions with finitely many steps, each of the pieces $U_i,L_i$ have the same property, and so do their images under $F$.
By the construction of the set $A_{a,b}$ we know that the levels corresponding to the ends of the cycles $c_a$ and $c_b$, if the cycles are strong, do not appear as horizontal boundary levels; the corresponding horizontal segments, let us call them the {\em locking segments}  lie in the interior of 
%\foot{SK: should we replace ``inside" by ``in the interior of"?} 
the set $A_{a,b}$. Furthermore, the images of all levels except for the levels next to the ends of the cycles, $f^{k_1-1}Ta$, $f^{m_1-1}Sa$, $f^{m_2-1}Sb$, and $f^{k_2-1}T^{-1}b$, also belong to $\U_{a,b}\cup\LL_{a,b}$.
The exceptional levels are exactly those between $0$ and $b$ and above $TSa$ in $\U_{a,b}$, and between $a$ and $0$ and below $T^{-1}Sb$ in $\LL_{a,b}$. The images of the horizontal segments belonging to these levels are the locking segments. Notice that the exceptional levels between $0$ and $b$ and between $a$ and $0$ constitute the horizontal boundary of the regions $U_3$ and $L_3$. 

Transporting the rays $[-\infty,x_b]$ and $[x_a,\infty]$ (with $x_a$ and $x_b$ uniquely determined by Lemma \ref{2-sys}),
along the corresponding cycles, and using the strong cycle property, we see that the ``locking segment" in the horizontal boundary of $U_1'$ coincides with the locking segment of the horizontal boundary of $L_3'$, and the locking segment in the horizontal boundary of $L_1'$ coincides with the locking segment of the horizontal boundary of $U_3'$. It can happen that both ``locking segments" belong to $A_{a,b}^u$ or $A_{a,b}^\ell$. If only one of the numbers $a$ or $b$ has the strong cycle property, then there will be only one locking segment. 

If the cycle property is weak or the $(a,b)$-continued fraction expansion of one or both $a$ and $b$ is periodic, then 
all levels of $\LL_a,\,\,\LL_b,\,\,\U_a$ and $\U_b$ will belong to the boundary of $A_{a,b}$, 
%\foot{SK20: is the case $b/(b+1)$ taken case of? I think, yes}
and there will be no locking segments. In these cases $L_3=[x_1,\infty]\times[a,0]$, and $L'_3=[-1/x_1,0]\times[-1/a,\infty]$, where $x_1=x_a$. Let $x_2$ be the $x$-coordinate of the  right vertical boundary segment of $U_2$. Then the $x$-coordinate of the  right vertical boundary segment of $U_1$ is $-1/x_2$. Let us denote the highest  level in $\U_{a,b}$ by $y_2$.
Since $y_2\leq -1/a+1$, $y_2-1\leq -1/a$ is the next level after $-1/a$ in $\U_{a,b}$. This is since if we had $y\in\U_{a,b}$ such that $y_2-1<y<-1/a$, its preimage $y'=Ty$ would satisfy  $y_2<y'<-1/a+1$, a contradiction. By construction of the region $A_{a,b}$ the segments at the levels $y_2-1$ and $-1/a$ are connected, therefore $Sx_1=T^{-1}Sx_2$.
%does not have horizontal segments between $\frac{a}{1-a}=STSa$ and $0$, 
%using the special form of the map $F$, we compute the $x$-coordinate of this vertical component to be $STSx_a$, and therefore the $x$-coordinate of the vertical boundary segment of $U'_1$ is $T^{-1}S(STSx_a)=-1/x_a$. 
This calculation shows that $L'_3$ and $U'_1$ do not overlap and fit together by this vertical ray.

Thus in all cases
the images $U'_i,\,L'_i$ do not overlap, and
$A_{a,b}=(\cup_{i=1}^3U'_i)\cup(\cup_{i=1}^3L'_i)$. This proves the bijectivity of the map $F$ on $A_{a,b}$ except for some images of its boundary. This completes the proof in the case $0<b\le -a<1$.
%\foot{SK20: added last sentence}

\medskip

%\noindent{\bf Case $a\le -1$:}
%\foot{IU: new case} 
%\foot{SK20: some changes in the last paragraph}
Now we return to the case $a\le -1$ dropped from consideration before Lemma \ref{2-sys}.
The explicit cycle relations for this case have been described in Theorem \ref{ab-cycle}. Notice that 
all lower levels are connected, and $T^mSb$ is connected with $a+1$. Therefore $y_\ell=TSb$, and this implies that $x_a=m$.
The upper levels in the positive part are
$$ST^{-1}b<ST^{-2}ST^{-1}b< ... <(ST^{-2})^{m-1}ST^{-1}b<a/(a+1)$$ 
and $y_u=T^{-1}(ST^{-2})^{m-2}ST^{-1}b$.
Lemma \ref{2-sys} in this case holds with $x_a=m$ and $x_b=-1$ since
the equation for adjacency of the levels $y_u$ and $Sa$ is
$$T^{-1}(ST^{-2})^{m-2}ST^{-1}x_b=ST^{m-1}Sx_b=-1/m\, ,$$ which implies $x_b=-1$. Lemma \ref{a-connected} also holds with $y_a^-=ST^{m-1}b$ and $y_a^+=ST^{m}b$. Lemma \ref{b-connected} holds with $y_b^-=T^{-1}Sa$ and $y_b^+=T^{-1}ST^{-1}b$ and all upper level will be connected by an argument similar to one described obove. To prove the bijectivity of $F$ on $A_{a,b}$ one proceeds the same way as above, the only modification being that level $L_2$ does not exist, and 
$L_3=\{(x,y)\in A_{a,b}^\ell, a\le y\le a+1\}$.
\end{proof}
The following corollary is evident from the proof of part (ii) of the above theorem.
%\foot{SK: changed statement slightly}
\begin{cor} \label{cor:bry} If both $a$ and $b$ have the strong cycle property, then for any boundary component
$h$ of $A_{a,b}$ (vertical or horizontal) there exists $N>0$ such that $F^N(h)$ is in the interior of $A_{a,b}$.%\foot{replaced ``inside" by ``in the interior of"}
\end{cor}

\section{Finite rectangular structure of the attracting set}\label{s:6}

Recall that the attracting set $D_{a,b}$ was defined by \eqref{def-atrac}: starting with the trapping region $\Theta_{a,b}$ described in Theorem \ref{Delta-trapping}, one has %the associated \emph{attractor set} 
\[ D_{a,b}=\bigcap_{n=0}^\infty D_n, \text{ with } D_n=\bigcap_{i=0}^n F^i(\Theta_{a,b})\,.
\]
%Let $\U_{a,b}=\U_a\cup\U_b$ and $ \LL_{a,b}=\LL_a\cup \LL_b$.
\begin{lem} \label{Dn} Suppose that the map $f$ satisfies the finiteness condition. Then,  
for each $n\geq 0$, $D_n$ is a region consisting of two connected components, the upper one, $D_n^u$, and the lower one, $D_n^\ell$, bounded by non-decreasing step-functions.
 %taking values in the sets $\U_{a,b}$ and $\LL_{a,b}$, respectively.
%\item All horizontal levels of the boundary of $D_n^u$ (resp., $D_n^\ell$) remain as horizontal levels of $D_{n+1}^u$ (resp., $D_{n+1}^\ell$); all levels of  $\U_{a,b}$ appear in the boundary of some $D_n^u$, and all levels of  $\LL_{a,b}$ appear in the boundary of some $D_n^\ell$.
%\end{enumerate}
\end{lem}
\begin{proof} 
The proof is by induction on $n$. The base of induction holds by the definition of the trapping region $\Theta_{a,b}$. For the induction step, let us assume that the region $D_n$ consists of two connected components, the upper one
$D_n^u$ and the lower one $D_n^\ell$, bounded by non-decreasing step-functions.
% and  $\LL_{a,b}$, respectively, and that all horizontal levels of $D_k$ for $k<n$ survive as the horizontal levels of $D_n$.
We will show that
the region $D_{n+1}$ consists of two connected components, %the upper one
$D_{n+1}^u$ and $D_{n+1}^\ell$, bounded by non-decreasing step-functions. 
%\item all horizontal levels of $D_n$ survive as the horizontal levels of $D_{n+1}$.
%\end{enumerate}

In what follows, we present the proof assuming that $0<b\le -a<1$. The situation $a\le -1$ is less complex due to the explicit cycle expressions described in Theorem \ref{ab-cycle} and can be treated similarly with some minor modifications.%\foot{IU: added}

%\foot{Labeling is compatible with the labeling on Fig. 2, probably need another picture here}
We decompose the regions $D_n^u$ and $D_n^\ell$ as follows
\[
\begin{split}
U_n^{11}&=\{(x,y)\in D_n^u\,:\, y\ge TSa\}\\
U_n^{12}&=\{(x,y)\in D_n^u\,:\, b\le y\le TSa\}\\
U_n^3&=\{(x,y)\in D_n^u\,:\, 0\le y\le b\}\\
U_n^{21}&=\{(x,y)\in D_n^u\,:\, \frac{a}{1-a}\le y\le 0\}\\
U_n^{22}&=\{(x,y)\in D_n^u\,:\, b-1\le y\le \frac{a}{1-a}\}\\
L_n^{11}&=\{(x,y)\in D_n^\ell\,:\, y\le T^{-1}Sb\}\\
L_n^{12}&=\{(x,y)\in D_n^\ell\,:\, T^{-1}Sb\le y\le a\}\\
L_n^3&=\{(x,y)\in D_n^\ell\,:\, a\le y\le 0\}\\
L_n^{21}&=\{(x,y)\in D_n^\ell\,:\, 0\le y\le \frac{b}{b+1}\}\\
L_n^{22}&=\{(x,y)\in D_n^\ell\,:\, \frac{b}{b+1}\le y\le a+1\}.\\
\end{split}
\]

By induction hypothesis, the regions $U_{12},\,U_n^3,\,U_n^{21}$ and $U_n^{22}$ are bounded below and above, and $U_n^{11}$ only below, by a ray and on the right by a non-decreasing step-function. Similarly, the regions $L_n^{12},\,L_n^3,\,L_n^{21}$ and $L_n^{22}$ bounded  above and below, and $L_n^{11}$ only above, by a ray and on the left by a non-decreasing step-function. 

If $B\subset D_n^u$
is one of the upper subregions, let $\partial B$ be the union of the boundary components of $B$ that belong to the boundary of $D_n^u$, and, similarly, if 
$B\subset D_n^\ell$
is one of the lower subregions, let $\partial B$ be the union of the boundary components of $B$ that belong to the boundary of $D_n^\ell$.

Since $\Theta_{a,b}$ is a trapping region, $F(\Theta_{a,b})\subset \Theta_{a,b}$, $D_{n+1}=F(D_n)\subset D_n$, and hence $D_{n+1}^u\subset D_n^u$ and $D_{n+1}^\ell\subset D_n^\ell$. 

The natural extension map $F$ is piecewise fractional-linear, hence it maps regions bounded by non-decreasing step-functions to regions bounded by non-decreasing step-functions. More precisely, we have
\[
\begin{split}
U_{n+1}^u&=S(U_n^{22}\cup U_n^{21})\cup T^{-1}(U_n^{11}\cup U_n^{12})\cup S(L_n^3)\\
U_{n+1}^\ell&=S(L_n^{22}\cup L_n^{21})\cup T(L_n^{11}\cup L_n^{12})\cup S(U_n^3)\,.
\end{split}
\]

In order to show that the region $D^u_{n+1}$, is connected, we notice
that the region $T^{-1}(U_n^{11}\cup U_n^{12})$ is inside the ``quadrant" $[-\infty,0]\times[b-1,\infty]$
 while $S(U_n^{22}\cup U_n^{21})$ is inside the strip $[0,1]\times[ST^{-1}b,\infty]$. Therefore, they either intersect by a ray of the $y$-axis, or are disjoint. In the first case, either $T^{-1}ST^{-1}b<Sa$, which implies that $S(L_n^3)$ is inside the connected region
 $S(U_n^{22}\cup U_n^{21})\cup T^{-1}(U_n^{11}\cup U_n^{12})$, or $Sa\leq T^{-1}ST^{-1}b$ which implies that the level $Sa$ belongs to the boundary of the trapping region, and again $S(L_n^3)$ is inside the connected region
 $S(U_n^{22}\cup U_n^{21})\cup T^{-1}(U_n^{11}\cup U_n^{12})$.
Now suppose that the regions  $T^{-1}(U_n^{11}\cup U_n^{12})$ and $S(U_n^{22}\cup U_n^{21})$ are disconnected. Notice that the right vertical boundary of the region $S(L_n^3)$ is a ray of the $y$-axis, thus $S(L_n^3)\cup S(U_n^{22}\cup U_n^{21})$ is a connected region bounded by a non-decreasing step-function. Since $T^{-1}(U_n^{12})\cap S(L_n^3)=\emptyset$, the non-connectedness situation may only appear from the intersection of $T^{-1}(U_n^{11})$ and $S(L_n^3)$,
i.e. inside the strip $[-1,0]\times[-1/a,\infty]$. 
Since $f$ satisfies the finiteness condition,
Theorem \ref{thm:recstructure} is applicable, and  the set $A_{a,b}$ constructed there belongs to each $D_n$. This is because $A_{a,b}\subset \Theta_{a,b}$, and if $A_{a,b}\subset D_n$, we have $A_{a,b}=F(A_{a,b})\subset F(D_n)=D_{n+1}$. The set
$A_{a,b}$ has finite rectangular structure
and contains the strip $[-1,0]\times[-1/a,\infty]$. Thus the connectedness
of the region $D^u_{n+1}$ is proved. Moreover, this argument shows that $\partial T^{-1}(U_n^{11})$ is inside $D^u_{n+1}$ and therefore does not contribute to its boundary, and
\[
\partial U_{n+1}^u=\partial (T^{-1}( U_n^{12}))\cup\partial (S(U_n^{22}\cup U_n^{21})\cup S(L_n^3)).
\]
Since $\partial (T^{-1}( U_n^{12})$  and $\partial (S(U_n^{22}\cup U_n^{21})\cup S(L_n^3))$ are given by non-decreasing step-functions, one $<Sa$, and the other $\ge Sa$, it follows that $\partial U_{n+1}^u$ is also given by a non-decreasing step-function. A similar argument proves that $D_{n+1}^\ell$ is connected and bounded by a non-decreasing step-function.
\end{proof}
\begin{lem}\label{l:all-levels} Suppose that, for each $n$, $D_n$  consists of two connected components as in Lemma \ref{Dn}.
Then 
\begin{enumerate}
\item all horizontal levels of the boundary of $D_n^u$ belong to $\U_{a,b}$ (resp., $D_n^\ell$ belong to $\LL_{a,b}$) and remain as horizontal levels of $D_{n+1}^u$ (resp., $D_{n+1}^\ell$); 
\item all levels of  $\U_{a,b}$ appear in the boundary of some $D_n^u$, and all levels of  $\LL_{a,b}$ appear in the boundary of some $D_n^\ell$;
\item
the attractor $D_{a,b}$ consists of two connected components bounded by non-decreasing step-functions; the upper boundary function takes all values from the set $\U_{a,b}$, and the lower boundary function takes all values from the set $\LL_{a,b}$. 
\item The map $F: D_{a,b}\to D_{a,b}$ is surjective.
\end{enumerate}
\end{lem}
\begin{proof} (1) We prove this by induction. For the base case, $D_0^u$ contains the horizontal levels $T^{-1}b$, $ST^{-1}b$ and $\min(T^{-1}ST^{-1}b,Sa)$. The levels $T^{-1}b$, $ST^{-1}b$ belong to the boundary of $D_1^u$. If $Sa<T^{-1}ST^{-1}b$, then $ST^{-1}b>TSa$ and therefore is the end of the cycle and does not belong to $\U_{a,b}$. If 
$Sa>T^{-1}ST^{-1}b$, then $T^{-1}ST^{-1}b$ appears as a boundary segment of $D_1^u$.
A similar argument  applies to $D_0^\ell$ that contains the horizontal levels $Ta$, $STa$, and either $TSTa$ or $Sb$.

For the induction step we assume that  (1) holds for $k=n-1$, and prove that it holds for $k=n$.
Let $y\in\partial D_n$ be a horizontal segment of the boundary, $y\geq ST^{-1}b$, and $y\in \U_{a,b}$. Then $y=Sy'$, where $y'\in\partial D_{n-1},\,b-1\leq y'<0$. %%%%START HERE!
By inductive hypothesis, $y'\in\partial D_n$, hence $y=Sy'\in\partial D_{n+1}$. Now let $y\in\partial D_n$ be a horizontal segment of the boundary, $b-1\leq y<Sa$. Then $y=T^{-1}y'$, where $y'\in\partial D_{n-1},\,0< y'<TSa$. By inductive hypothesis, $y'\in\partial D_n$, hence $y=Sy'\in\partial D_{n+1}$.

The level $y=Sa$ appears as a boundary segment of $D_n^u$ since $T^{-1}(\partial(U_{n-1}^{11})\cup \partial(U_{n-1}^{12}))$ and  $S(\partial(L_{n-1}^3))$ do not overlap. Then $y=Sy'$, where $y'=a$ is the $y$-coordinate of the horizontal lower boundary of $L_{n-1}^3$. Since $L_{n}^3\subset L_{n-1}^3$ and $U_{n}^{11}\cup U_{n}^{12}\subset U_{n-1}^{11}\cup U_{n-1}^{12}$, we get that $T^{-1}(\partial(U_{n}^{11})\cup \partial(U_{n}^{12}))$ and  $S(\partial(L_{n}^3))$  do not overlap, and $y=Sa$ will appear as a boundary segment of $D_{n+1}^u$.

On the other hand, assume $y\in\partial D_{n+1}$ was not a horizontal level of $\partial D_n$. Then $y=Sy'$ for some $y'\in\partial (U_n^{22}\cup U_n^{21})$, $y=T^{-1}y'$ for some $y'\in\partial (U_n^{12})$,  or $y=Sa$. In all cases $y\in \U_{a,b}$ 
by the structure of the sets $\U_{a}$ and $\U_b$ established in Theorems \ref{a-cycle} and \ref{b-cycle}.

(2)
We start with level $-\frac1{b-1}$ which belongs to the boundary of the trapping region $\Theta_{a,b}$ by definition. We have seen that if $T^{-1}ST^{-1}b\in \U_b$, %and $T^{-2}ST^{-1}b$ 
%\foot{SK: some small changes}
then the level appears in the boundary of $D_1^u$. Now, if $b-1<T^{-k}ST^{-1}b<\frac{a}{1-a}$ (for the smallest $k=2$ or $3$), then the expansion continues, each $T^{-i}ST^{-1}b,\,i\leq k$ appears for the first time in the boundary of $D_i^u$ for $i\leq k$, and
the next element in the cycle, $ST^{-k}ST^{-1}b$, appears in the boundary of $D_{k+1}^u$. Using the structure of the set $\U_b$
established in Theorem \ref{b-cycle} we see that all levels of the set $\U_b$ appear as boundary levels of some $D_n^u$. We use the same argument for level $-\frac1{a}$ which appears for the first time in the boundary of some $D_{n_0}^u$, to see that all elements of the set $\U_a$ appear as boundary levels of all successive sets $D_n^u$. The same argument works for the lower boundary.

(3) Thus starting with some $n$, all sets $D_n$ have two connected components bounded by non-decreasing step-functions whose $y$ levels coincide with the sets $\U_{a,b}$ and $\LL_{a,b}$.
Therefore, the attractor $D_{a,b}=\cap_{n=0}^\infty D_n$ has the same property. 

(4) The surjectivity of the map $F$ on $D_{a,b}$ follows from the nesting property of the sets $D_n$.
%Since the map $F$ is monotone increasing in each variable, mapping horizontal lines to horizontal and vertical lines to vertical, for each $n\geq 0$, $D_n$ is a region consisting of ...
\end{proof}
A priori the map $F$ on $D_{a,b}$ does not have to be injective, but in our case it will be since we will identify $D_{a,b}$ with an earlier constructed set $A_{a,b}$.
\begin{cor} \label{FRS}If the map $f$  satisfies the finiteness condition, then the attractor $D_{a,b}$ has finite rectangular structure, i.e. bounded by non-decreasing step-functions with a finite number of steps.
\end{cor}

%Now we can prove the main result of this section.
\begin{thm} \label{attractor} If the map $f$ satisfies the finiteness condition, then the set $A_{a,b}$ constructed in Theorem \ref{thm:recstructure} is the attractor for the map $F$.
\end{thm}
\begin{proof} We proved in Theorem \ref{thm:recstructure}  that the set $A_{a,b}$ constructed there is uniquely determined by the prescribed set of $y$-levels $\U_{a,b}\cup\LL_{a,b}$. By Corollary \ref{FRS}, the set $D_{a,b}$ has finite rectangular structure with the same set of $y$-levels. Now we look at the $x$-levels of the jumps of its boundary step-functions.
Take the vertex $(x, b-1)$ of $D_{a,b}$. From the surjectivity of $F$ on $D_{a,b}$, there is a point $z\in D_{a,b}$ s.t. $F(z)=(x, b-1)$.
Then $z$ must be the intersection of the ray at the level $b$ with the boundary of $D_{a,b}$, i.e. $z=(\tilde x_b, b)$, hence $x=\tilde x_b-1$. Continue the same argument: look at the vertex at the level $-1/(b-1)$. It must be $F(\tilde x_b-1, b-1)$, etc. Since each $y$-level of the boundary has a unique ``predecessor" in its orbit, all $x$-levels of the jumps obtained by ``transporting" the rays $[-\infty,\tilde x_b]$ and $[\tilde x_a,\infty]$ over the corresponding cycles, satisfy the same
 equations that defined the boundary of the set $A_{a,b}$ of Theorem \ref{thm:recstructure}. Therefore $\tilde x_a=x_a,\,\,\tilde x_b=x_b$,
the step-functions that define the boundaries are the same, and $D_{a,b}=A_{a,b}$.
\end{proof}

\section{Reduction theory conjecture}\label{s:naturalextension}\label{s:7}
\begin{figure}[thb]
%\vspace*{-.5cm}
  \begin{minipage}[b]{.3\textwidth}
 % \vspace*{-.5cm}
    \begin{center}
\includegraphics[height=1.5in]{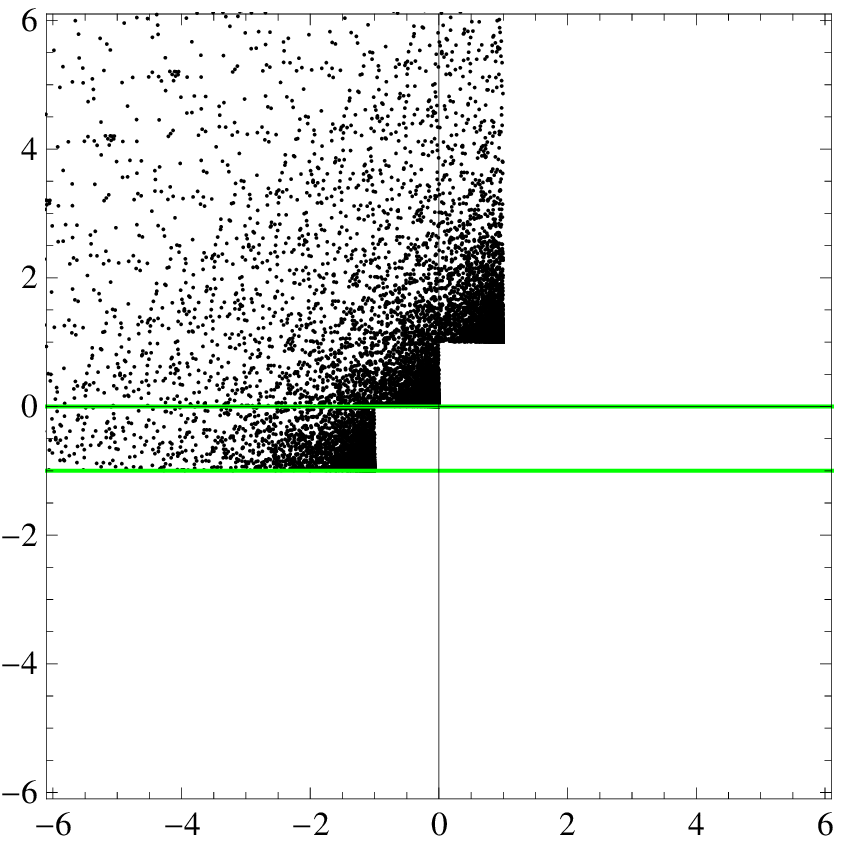}
\center{$A_{-1,0}$}
 \end{center}
  \end{minipage}
  \hfill
  \begin{minipage}[b]{.3\textwidth}\mbox{ }
    \begin{center}
\includegraphics[height=1.5in]{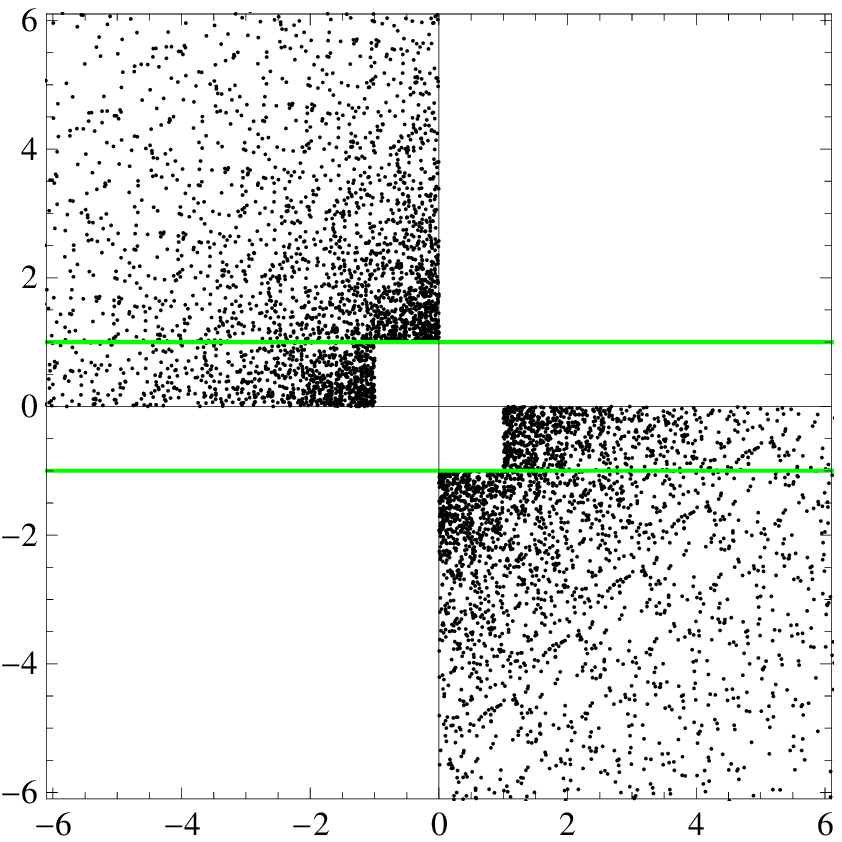}
\center{$A_{-1,1}$}
 \end{center}
  \end{minipage}
   \hfill
  %\vspace*{0.2in}
   \begin{minipage}[b]{.3\textwidth}\mbox{ }
   %\vspace*{0.2in}
    \begin{center}
\includegraphics[height=1.5in]{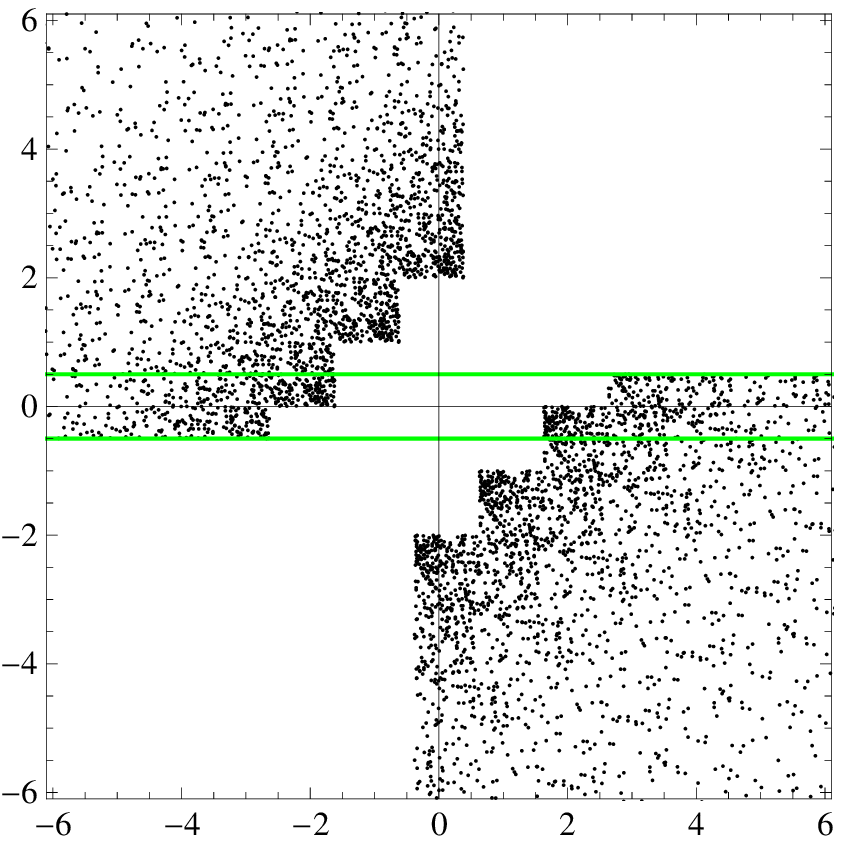}
\center{$A_{\tiny -1/2,1/2}$}
%Hurwitz- or 
%$H$-expansion}
 \end{center}
  \end{minipage}
  \caption{Attractors for the classical cases}\label{fig:A}
  \end{figure}

Don Zagier conjectured that the Reduction Theory properties, stated in the Introduction, hold for every $(a,b)\in\mathcal P$. 
He was motivated by the classical cases and computer experimentations with random parameter values $(a,b)\in\mathcal P$ (see Figures \ref{don-a} and \ref{fig:A} for  attractors obtained by iterating 
random points using Mathematica program). 

The following theorem gives a sufficient condition for the Reduction Theory conjecture to hold:
\begin{thm} \label{RTC}If both $a$ and $b$ have the strong cycle property, then for every point $(x,y)\in \bar\R^2\setminus \Delta$ there exists $N>0$ such that $F^N(x,y)\in D_{a,b}$.
\end{thm}
\begin{proof} Every point  $(x,y)\in\bar\R^2\setminus \Delta$ is mapped to the trapping region by some iterate $F^{N_1}$. Since the sets $D_n$ are nested and contain $D_{a,b}$, for large $N$, $F^N(x,y)$ will be close to the boundary of $D_{a,b}$. By Corollary \ref{cor:bry}, for any boundary component $h$ of $D_{a,b}$ there exists $N_2>0$ such that $F^{N_2}(h)$ is inside $D_{a,b}$. Therefore, there exists a large enough $N>0$ such that $F^N(x,y)$ will be in the interior of $D_{a,b}$.
%\foot{replaced ``inside" by "in the interior of"}
\end{proof}

The strong cycle property is not necessary for the Reduction theory conjecture to hold. For example, it holds for the two classical expansions $(-1,0)$ and $(-1,1)$ that satisfy only a weak cycle property. In the third classical expansion $(-1/2,1/2)$ that also satisfies a weak cycle property, property (3) does not hold for some points $(x,y)$ with $y$ equivalent to $r=(3-\sqrt{5})/2$.
\noindent

%We prove if $(a,b)\in \mathcal P$ satisfies the finiteness condition but not the strong cycle property, then the above result remains true for almost every point $(x,y)$ of the plane.

The next result shows that, under the finiteness condition, almost every point $(x,y)\in \bar\R^2\setminus \Delta$ lands in the attractor $D_{a,b}$ after finitely many iterations. %\foot{SK: will need to change the argument here}
\begin{prop}
If the map $f_{a,b}$ satisfies the finiteness condition, then for almost every point $(x,y)\in \bar\R^2\setminus \Delta$, there exists $N>0$ such that
$F_{a,b}^N(x,y)\in D_{a,b}$.
\end{prop}
\begin{proof}
Let $(x,y)\in \R^2$ with $y$ irrational and $y=\lfloor n_0,n_1,n_2,\dots\rceil_{a,b}$. In the proof of Theorem \ref{Delta-trapping}, we showed that there exists $k>0$ such that 
$$(x_{j+1},y_{j+1})=ST^{-n_j}\dots ST^{-n_1}ST^{-n_0}(x,y)\in [-1,1]\times ([-1/a,\infty]\cup [-\infty,-1/b])$$
for all $j\ge k$. The point $F^{N}_{a,b}(x,y)=(x_{k+1}, y_{k+1})$  is in $A_{a,b}$, if $(x_{k+1},y_{k+1})\in [-1,0]\times [-1/a,\infty]$ or
$(x_{k+1},y_{k+1})\in [0,1]\times [-\infty,-1/b]$. Also, $F^{N+1}(x,y)=F(x_{k+1},y_{k+1})$ is in $A_{a,b}$ if $(x_{k+1},y_{k+1})\in [0,1]\times [-1/a+1,\infty]$
or $(x_{k+1},y_{k+1})\in [-1,0]\times [-\infty,-1/b-1]$. Thus we are left with analyzing the situation when the sequence of iterates 
$$(x_{j+1},y_{j+1})=ST^{-n_j}\dots ST^{-n_1}ST^{-n_0}(x,y)$$ belongs to $[0,1]\times [-1/a,-1/a+1]$ for all $j\ge k$ (or $[-1,0]\times [-1/b,-1/b-1]$ for all $j\ge k$). Assume that we are in the first situation: $y_{j+1}\in [-1/a,-1/a+1]$ for all $j\ge k$. This implies that all digits $n_{j+1}$, $j\ge k$ are either $\lfloor -1/a\rceil$ or $\lfloor -1/a\rceil+1$.  In the second situation, the digits $n_{j+1}$, $j\ge k$ are either $\lfloor -1/b\rceil$ or $\lfloor -1/b\rceil-1$. Therefore the continued fraction expansion of $y$ is written with only two consecutive digits (starting from a certain position). By using Proposition \ref{bdigits1} and Remark \ref{bdigits2} we obtain that the set of all such points has zero Lebesgue measure. This proves our result.
\end{proof}
\begin{rem}
In the next section we show that there is a non-empty Cantor-like set $\E\subset \Delta$ belonging to the boundary segment $b=a+1$ of $\mathcal P$ such that for $(a,b)\in\E$ the set $\U_{a,b}\cup\LL_{a,b}$ is infinite. Therefore, for $(a,b)\in\E$ either the set $D_n^u$ or $D_n^\ell$ is disconnected for some $n>0$, or,
by Lemma \ref{l:all-levels}(3),  the attractor $D_{a,b}$ consists of two connected components whose boundary functions are not step-functions with finitely many steps.
\end{rem}

\section{Set of exceptions to the finiteness condition}\label{s:8}
In this section we study the structure of the  set $\E\subset\mathcal P$ of exceptions to the finiteness condition. 
%\foot{SK: change; in Theorems \ref{b-cycle}, \ref{a-cycle} only no cycle property is discussed}
%In this section we study the structure of the exceptional set $\E$ of $(a,b)$ that 
%appears in part (II) of Theorems \ref{b-cycle}, \ref{a-cycle}. 
We write $\E=\E_\bb\cup\E_\af$
%, where $\E_\bb=\cup_{m=3}^\infty \E_b^m$
% is in the union of blue ``tongues" 
% and $\E_a=\cup_{m=3}^\infty \E_a^m$
% is in the union of red ``tongues". 
where $\E_\bb$ (resp., $\E_\af$) consists of all points $(a,b)\in \P$ for which  $b$ (resp., $a$) does not satisfy the finiteness condition, i.e. either the truncated orbit $\U_b$ or $\LL_b$ is infinite (resp., $\U_a$ or $\LL_a$). 
%\foot{SK: added ``truncated"}

We analyze the set $\E_\bb$.  Recall that, by Proposition \ref{shift}(2), the set $\U_b$ is infinite if and only if $\LL_b$ is infinite, therefore it is sufficient to analyze the condition that the orbit $\U_b$  is not eventually periodic and its values belong to the interval $(\frac{b}{b+1},a+1)$. As before, we restrict our analysis (due to the %\foot{SK: replaced ``some" by ``the"}
 symmetry considerations) to the parameter subset of $\P$ given by $b\le -a$ and write 
$\E_\bb=\cup_{m=3}^\infty\E^m_\bb$ where $b\in\E^m_\bb$ if $b\in \E_\bb$ and $T^mSb\in (\frac{b}{b+1},a+1)$.
% and  $\hat f^{(k)}=T^h S,$ for $k\geq 3$, where $\h=m$ or $m+1$. 
By Theorem \ref{b-cycle} and its proof, it follows that if $b\in \E^m_\bb$, then the first digit of the $(a,b)$-continued fraction expansion of $Sb$ is $-m$ and all the other digits are either $-m$ or $-(m+1)$.
%\newpage

We describe a recursive construction of the exceptional set $\E^m_\bb$. One starts with the `triangular' set 
%(the blue ``tongue")
%\foot{SK:change!}
$$\T^m_\bb=\{(a,b)\in \P : \frac{b}{b+1}\le T^mSb\le a+1\}.$$
The range of possible values of $b$ in $\T^m_\bb$ is given by the interval $[\underline b, \bar b]$ where $T^mS\bar b=\bar b$ and $T^mS\underline b=\underline b/(\underline b+1)$. Since
\[
\frac{b}{b+1}\leq b \text{ for all } b\ge 0,
\]
and the function $T^mSb$ is monotone increasing, we obtain that $\underline b<\bar b$, and $\underline b$ is the horizontal boundary of $\T^m_\bb$, while $\bar b$ is the $b$-coordinate of its `vertex'.

%we look for all points whose future iterations under the map $\hat f$ belong to the interval $(\frac{b}{b+1},a+1)$. 
At the next stage we obtain the following  regions:
$$\T^{m,m}_\bb=\{(a,b)\in \T^m_\bb: \frac{b}{b+1}\le T^{m}ST^mSb\le a+1\}$$
$$\T^{m,m+1}_\bb=\{(a,b)\in \T^m_\bb: \frac{b}{b+1}\le T^{m+1}ST^mSb\le a+1\}\,.$$
%The defining inequality $T^mST^mSb\le a+1$ of $\T^{m,m}$ is satisfied by all points of $\T^m$ because the orientation preserving property of $T,S$ and the fact that $T^mSb\le a+1$ implies that
%$$T^mST^mSb\le T^mS(a+1)\le T^mSb\le a+1\,.$$
%Thus the left side of the region  $\T^{m,m}$ is part of the left side of $\T^m_b$. 
By the same argument as above each region is `triangular', i.e. the $b$-coordinate of its lower (horizontal) boundary is less than the $b$-coordinate of its vertex. We show that its intersection with the triangular region obtained on the previous step is either empty or has `triangular' shape. The horizontal boundary of $\T^{m,m}_\bb$ has the $b$-coordinate given by the relation $T^mST^mSb=b/(b+1)$ (call it $\tilde b$). 
%We will show that $\underline b<\tilde b<\bar b$. 
We have
\[
T^mST^mS\underline b=T^mS\left(\frac{\underline b}{\underline b+1}\right)=T^mS\underline b-1=-\frac1{\underline b+1}<\frac{\underline b}{\underline b+1},
\]
so $\underline b<\tilde b$. On the other hand, 
\[
T^mST^mS\bar b=T^mS\bar b=\bar b,
\]
which shows that the hyperbola $T^mST^mSb=b$ intersects the diagonal side $b=a+1$ at the point with $b$-coordinate $\bar b$. It follows that the region $\T^{m,m}_\bb$ is triangular and non-empty with $\underline b<\tilde b<\bar b$.

%Such a value of $b$  exists because the range of possible values of $b$ in $\T^m_\bb$ is given by the interval $(\underline b, \bar b)$ where $T^mS\bar b=\bar b$ and $T^mS\underline b=\underline b/(\underline b+1)$
%so
%$$ -1/(\underline b+1)= T^mS\underline b-1=T^mS({\underline b}/{\underline b+1})< T^mST^mSb< T^mS\bar b= \bar b\,.$$
%The defining inequality $T^{m+1}ST^b\ge \frac{b}{b+1}$ of $\T^{m,m+1}$ is satisfied by all points of $\T^m$ because $T^mSb\ge \frac{b}{b+1}$ implies that
%$$T^{m+1}ST^mSb\ge T^{m+1}S\frac{b}{b+1}=m-\frac{1}{b}=T^mSb\ge \frac{b}{b+1}\,.$$
%Thus the horizontal side of $T^{m,m+1}$ is part of the horizontal side of $\T^m_b$. 
The upper boundary of $\T^{m,m+1}_\bb$ is given by the hyperbola $T^{m+1}ST^mSb=a+1$. Notice that, if $\underline a+1=T^mS\underline b$, then the point $(\underline a,\underline b)$ lies on the curves $T^mSb=a+1$ (obviously) and $T^{m+1}ST^mSb=a+1$ because 
$$T^{m+1}ST^mS\underline b=T^{m+1}S(\underline b/(\underline b+1))=T^mS\underline b=\underline a+1\,.$$ 
This shows that the entire horizontal boundary of $\T^m_\bb$ belongs to that of $\T^{m,m+1}_\bb$. 
%\foot{SK: misprint: replaced $\T^{m,m}_bb$ by $\T^{m,m+1}_bb$}
Moreover, the hyperbola $T^{m+1}ST^mSb=a+1$ intersects the diagonal side $b-a=1$ at the point $\hat b$ satisfying $T^{m+1}ST^mS\hat b=\hat b$. Therefore,  $T^{m}ST^mS\hat b=\hat b-1<\frac{\hat b}{\hat b+1}$, i.e. $\hat b<\tilde b$.
In this case we have
$\underline b<\hat b<\tilde b<\bar b$, and the two triangular regions $\T^{m,m}_\bb$ and $\T^{m,m+1}_\bb$ are disjoint and non-empty. 

\begin{figure}[ht]
\psfrag{B}[r]{$\T^{n_1,n_2,\dots,n_k,m}_\bb$}
\psfrag{A}[l]{$\T^{n_1,n_2,\dots,n_k,m+1}_\bb$}
\psfrag{c}[c]{$\underline b$}
\psfrag{d}[c]{$\hat b$}
\psfrag{e}[c]{$\tilde b$}
\psfrag{f}[c]{$\bar b$}
\psfrag{u}[c]{$\underline a$}
\psfrag{v}[c]{$\underline a'$}
\centerline{\includegraphics[scale=0.8]{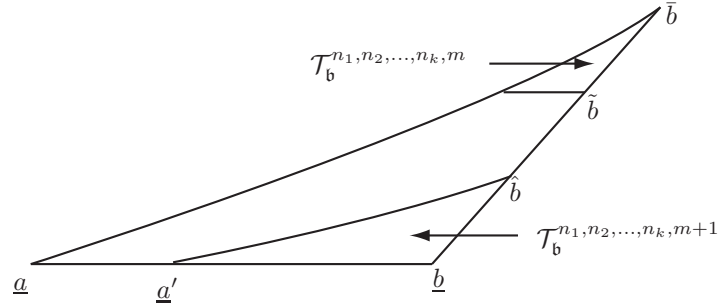}}
\caption{Set $\T^{n_1,n_2,\dots,n_k}_\bb$ and its two subregions}
\label{fig-exc}
%\end{center}
\end{figure}

The situation becomes more complicated as we proceed recursively. Let $\T^{n_1,n_2,\dots,n_k}_\bb$ be one of the regions obtained after $k$ steps of this construction, with $n_1=m$ and $n_i\in\{m,m+1\}$ for $2\le i\le k$. At the next step we get two new sets (possible empty) (see Figure \ref{fig-exc}):
$$\T^{n_1,n_2,\dots,n_k,m}_\bb=\{(a,b)\in \T^{n_1,n_2,\dots,n_k}_\bb: \frac{b}{b+1}\le T^{m}ST^{n_k}S\dots T^{n_1}Sb\le a+1\}$$
$$\T^{n_1,n_2,\dots,n_k,m+1}_\bb=\{(a,b)\in \T^{n_1,n_2,\dots,n_k}_\bb: \frac{b}{b+1}\le T^{m+1}ST^{n_k}S\dots T^{n_1}Sb\le a+1\}\,.$$

%Let us prove now the above statements. \foot{Important!!!}
%First, notice that 
%\foot{SK:change}
As in the base case, the inequality $T^{m}ST^{n_k}S\dots T^{n_1}Sb\le a+1$ of $\T^{n_1,n_2,\dots,n_k,m}_\bb$ is satisfied by all points of $\T^{n_1,n_2,\dots,n_k}_\bb$ because of the monotone increasing property of $T,S$ and the fact that $T^{n_k}S\dots T^{n_1}Sb\le a+1$ implies
$$T^{m}ST^{n_k}S\dots T^{n_1}Sb\le T^mS(a+1)\le T^mS(b)\le a+1\,.$$
Thus the upper boundary of the region  $\T^{n_1,n_2,\dots,n_k,m}_\bb$ (if nonempty) is part of the upper boundary of $\T^{n_1,n_2,\dots,n_k}_\bb$; it is the lower (horizontal) boundary that changes. Similarly, the defining inequality $\frac{b}{b+1}\le T^{m+1}ST^{n_k}S\dots T^{n_1}Sb$ of $\T^{n_1,n_2,\dots,n_k,m+1}_\bb$  is satisfied by al points of  $\T^{n_1,n_2,\dots,n_k}_\bb$ because
$$T^{m+1}ST^{n_k}S\dots T^{n_1}Sb\ge T^{m+1}S\frac{b}{b+1}=m-\frac{1}{b}=T^mSb\ge\frac{b}{b+1}\,.$$
Thus the lower boundary of $\T^{n_1,n_2,\dots,n_k,m+1}_\bb$ (if nonempty) is part of the lower boundary of $\T^{n_1,n_2,\dots,n_k}_\bb$. 
Therefore, we can describe the above sets as
\begin{eqnarray}
& &\T^{n_1,n_2,\dots,n_k,m}_\bb = \{(a,b)\in \T^{n_1,n_2,\dots,n_k}_\bb: \frac{b}{b+1}\le T^{m}ST^{n_k}S\dots T^{n_1}Sb\}\\
& &\T^{n_1,n_2,\dots,n_k,m+1}_\bb = \{(a,b)\in \T^{n_1,n_2,\dots,n_k}_\bb: T^{m+1}ST^{n_k}S\dots T^{n_1}Sb\le a+1\}\,.
\end{eqnarray}
%\foot{SK: added sentence}
By the same reason as in the base case, the two regions $\T^{n_1,\dots,n_k,m}_\bb$ and $\T^{n_1,\dots,n_k,m+1}_\bb$ do not overlap.

The set $\E^m_\bb$ is now obtained as the union of all sets of type
\begin{equation}\label{eni}
\E^{(n_i)}_\bb=\bigcap_{k=1}^\infty\T^{n_1,n_2,\dots,n_k}_\bb
\end{equation}
where
$n_1=m$, $n_i\in\{m,m+1\}$ if $i\ge 2$, and the sequence $(n_i)$ is not eventually periodic. If such a set $\E^{(n_i)}_\bb$ is non-empty and $(a,b)$ belongs to it, then $b$ is uniquely determined from the $(a,b)$-expansion of $Sb=\lfloor -n_1,-n_2,\dots\rceil$. 
%We shall prove that $\E^{(n_i}_b$ is either a line segment or a point on $b-a=1$.

First we need some additional lemmas:
%\foot{SK:fixed the notation $\lfloor Sb\rceil_{a,b}$}
\begin{lem}\label{lem-cf} $ $
\begin{itemize}
\item[(i)] A point $b\in [0,1]$ satisfying  $T^{n_k}S\dots T^{n_1}Sb=b$ with $|n_i|\ge 2$ can be written formally using a periodic ``$-$" continued fraction expansion 
\begin{equation}\label{bbar}
b=-1/(\overline{-n_1,-n_2,\dots,-n_k})=(0,\overline{-n_1,-n_2,\dots,-n_k})\,.
\end{equation}
If $b$ is in $\T^{n_1,n_2,\dots,n_k}_b$, then $Sb$ has the $(a,b)$-continued fraction expansion 
$$\lfloor Sb\rceil_{a,b}=\lfloor\overline{-n_1,-n_2,\dots,-n_k}\,\rceil\,.$$
\item[(ii)] A point $b$ in $[0,1]$ satisfying  $T^{n_k}S\dots T^{n_1}Sb=b/(b+1)$ can be written formally using the  periodic ``$-$'' continued fraction expansion  
\begin{equation}\label{underb}
b=(0,-n_1,\overline{-n_2,\dots,-n_k,-(m+1)})\,.
\end{equation}
If the point $b\in T^{n_1,n_2,\dots,n_k}_b$, then $\lfloor Sb\rceil_{a,b}=\lfloor -n_1,\overline{-n_2,\dots,-n_k,-(m+1)}\rceil$.
\end{itemize} 
\end{lem}
\begin{proof} 
One can verify directly that the point $b$ given by \eqref{bbar} is the fixed point of the hyperbolic transformation $T^{n_k}S\dots T^{n_1}S$ and $b\in [0,1]$ (see also \cite[Proposition 1.3]{KU1}).
%Since $b$ is a point in $\T^{n_1,n_2,\dots,n_k}_b$, its continued fraction expansion is given by $b=[0,-n_1,-n_2,\dots,-n_k,-n_{k+1},\dots]$. Notice that 
%$$T^{n_k}S\dots T^{n_1}Sb=[0,-n_{k+1},-n_{k+2},\dots]$$
%so the part (i) is true.

The equation in part (ii) can be written as  $STST^{n_k}S\dots T^{n_1}Sb=b$ and one verifies directly that the value $b$ given by \eqref{underb} is the fixed point of that hyperbolic transformation and $b\in [0,1]$.
\end{proof}

Notice that the relation $(0,-n_1,-n_2,\dots)=-(0,n_1,n_2,\dots)$ is satisfied, assuming that the formal ``$-$" continued fraction expansions are convergent (from the proof of Theorem \ref{convergence}, the convergence property holds if $|n_i|\ge 2$ for all $i\ge 1$).
%\foot{IU: added reference}

\begin{defn}
We say that two sequences (finite or infinite) $\sigma_1=(n_i)$ and $\sigma_2=(p_j)$ of positive integers are in lexicographic order, $\sigma_1\prec \sigma_2$, if on the first position $k$ where the two sequences differ one has $n_k<p_k$ ,or if the finite sequence $(n_i)$ is a starting subsequence of $(p_j)$. 
\end{defn}
The following property follows from the monotonicity of $T,S$.
\begin{lem}\label{order}
Given two infinite sequences $\sigma_1=(n_i)$ and $\sigma_2=(p_j)$ of integers $n_i\ge 2$ and $p_j\ge 2$ such that $\sigma_1\prec \sigma_2$ then
$$(0,n_1,n_2,\dots) < (0,p_1,p_2,\dots)\,.$$ 
\end{lem}

The next lemma provides necessary conditions for a set $\E^{(n_i)}_\bb$ to be non-empty. Denote by $\mathfrak{l}_m$ the length of the initial block of $m$'s and by $\mathfrak{l}_{m+1}$ the length of the first block of $(m+1)$'s in $(n_i)$.

\begin{lem}\label{lemE} $ $
\begin{itemize}
\item[(i)] If a set $\E^{(n_i)}_\bb$ in the upper region $\T^{m,m}_\bb$ is non-empty then the sequence $(n_i)$ contains no consecutive $(m+1)$'s and the length of any block of $m$'s is equal to $\frak{l}_m$ or $\frak{l}_m-1$.
\item[(ii)] If a set $\E^{(n_i)}_\bb$ in the lower region $\T^{m,m+1}_\bb$ is non-empty then the sequence $(n_i)$ contains no consecutive $m$'s and the length of any block of $(m +1)$'s is %\foot{SK: change according to the wording in (i)}
%at most the length of the first block of 
equal to $\frak{l}_{m+1}$ or $\frak{l}_{m+1}+1$.
\end{itemize}
\end{lem}

\begin{proof}
(i)  Assume that the sequence $(n_i)$ contains two consecutive $(m+1)$'s. 
%\foot{SK: edited}
Then some $\T^{n_1,n_2,\dots,n_k,m+1,m+1}_\bb$ (with $n_1=n_2=n_k=m$) is non-empty. The upper vertex of such a triangular set satisfies the inequality
%\foot{SK: replaced  $n_k$ by $m$}
%\foot{SK: we have to be consistent: $m+1$ or $(m+1)$ thereafter?}
\[
\begin{split}
\bar b&\le -(0,\overline{n_1,n_2,\dots,n_k,m+1,m+1})\\
&=-(0,m,m,\dots,m,m+1,\boxed{m+1},\dots) 
\end{split}
\]
while the lower (horizontal) boundary satisfies
%\foot{SK: replaced  $n_k$ by $m$}
\[
\begin{split}
\underline{b}&\ge -(0,n_1,\overline{n_2,\dots,n_k,m+1})\\
&=-(0,m,m,\dots,m,m+1,\boxed{m},\dots)\qquad
\end{split}
\]
This implies that $\underline{b}>\bar b$ because the entries of the  corresponding continued fractions 
with positive entries are in lexicographic order (they coincide on the first $k+1$ places, and on the $(k+2)^{th}$
%\foot{SK: change} 
position the first continued fraction has digit $m+1$ while the second one has digit $m$), i.e. the set $\T^{n_1,n_2,\dots,n_k,m+1,m+1}_\bb$ is empty.

Now assume that there exists a non-empty set $\T^{n_1,n_2\dots,n_k, m, m, \dots, m}$ ($n_k=m+1$) with the final block of $m$'s of length greater than $\frak{l}_m$. The upper vertex of this set is given by 
%\[
%\begin{split}
%\bar b&=[0,\overline{\underbrace{-m,\dots,-m}_p,-(m+1)}]\\
%&=[0,\underbrace{-m,-m,\dots,-m}_p,-(m+1),\underbrace{-m,-m,\dots,-m}_p,-(m+1),\dots]
%\end{split}
%\]
%while the lower horizontal segment is given by 
%\[
%\begin{split}
%\underline{b}&=[0,-m,\overline{-m,\dots,-m,-(m+1),\underbrace{-m,-m\dots, -m}_q,-(m+1)}]\\
%&=[0,\underbrace{-m,-m,\dots,-m}_{p},-(m+1),\underbrace{-m,-m\dots, -m}_q,-(m+1),\underbrace{-m,\dots -m}_{p-1},-(m+1),\dots]\,.
%\end{split}
%\]
%If $p<q$ then the $\bar b,\underline{b}$ coincide on the first $k+p$ entries. Looking at the $k+p+1$ entry, one has $\bar b_{k+p+1}=-(m+1)$ and $\underline{b}_{k+p+1}=-m$, so $\bar b<\underline{b}$. If $p=q$, then the first $k+We assumed that $p\le q$ therefore $\bar b<\underline b$, hence the contradiction. The very same argument shows that the conditions listed are also sufficient, since one has that $\bar b>\underline b$.
\[
\begin{split}
\bar b&\le -(0,\overline{n_1,n_2,\dots,n_k})=-(0,\overline{\underbrace{m,m,\dots,m}_{\frak{l}_m},m+1,\dots,n_k})\\
&=-(0,\underbrace{m,m,\dots,m}_{\frak{l}_m},m+1,\dots,n_k,\underbrace{m,m,\dots,m}_{\frak{l}_m},m+1,\dots)
\end{split}
\]
while the lower horizontal segment is given by 
\[
\underline{b}\ge -(0,n_1,\overline{n_2,\dots,n_k,\underbrace{m,m\dots, m}_q,m+1})\,.
\]
If $\frak{l}_m<q$ then the two continued fractions coincide on the first $k+p$ entries. Looking at the $k+p+1$ entry, we get that $\bar b<\underline{b}$, hence the set $\T^{n_1,n_2\dots,n_k, m, m, \dots, m}_\bb$ would be empty.

Assume now that there exists a non-empty set of type $\T^{n_1,n_2\dots,n_k, m, m, \dots, m,m+1}_\bb$ ($n_k=m+1$) with the last block of $m$'s of length $q$ strictly less than $\frak{l}_m-1$. Because $n_k=m+1$, $n_{k-1}=m$, and $\T^{n_1,n_2\dots,n_k, m, m, \dots, m,m+1}_\bb\subset \T^{n_1,n_2\dots,n_k}_\bb$ we have that the lower limit of the set $\T^{n_1,n_2\dots,n_k, m, m, \dots, m,m+1}_\bb$ satisfies the relation
\[
\begin{split}
\underline b&\ge -(0,n_1\overline{n_2,\dots,n_{k-1},m+1})=-(0,n_1\overline{n_2,\dots,n_{k-1},n_k})\\
&=-(0,\underbrace{m,m,\dots,m}_{\frak{l}_m},m+1,\dots,n_k,\underbrace{m,\dots,m}_{\frak{l}_m-1},m+1,\dots)
%&=[0,\underbrace{-m,-m,\dots,-m}_p,-(m+1),\dots,-n_k,\underbrace{-m,-m,\dots,-m}_p,-(m+1),\dots]
\end{split}
\]
while the upper limit of the same set satisfies the relation
\[
\bar b\le -(0,\overline{n_1, n_2,\dots,n_k,\underbrace{m,m\dots, m}_q,m+1}]\,.
\]
This implies that $\bar b<\underline b$ because the two continued fractions coincide on their first $k+q$ entries, and the $k+q+1$ entries are $m$, and $m+1$ respectively. Therefore the set $\T^{n_1,n_2\dots,n_k, m, m, \dots, m,m+1}_\bb$ is empty.

\medskip

(ii)  Assume that a set $\T^{n_1,n_2,\dots,n_k,m,m}_\bb$ (with $n_1=m$, $n_2=m+1$ and $n_k=m+1$) is non-empty. The upper vertex of such a set satisfies the inequality
\[
%\begin{split}
\bar b\le -(0,\overline{n_1,n_2,\dots,n_k})
=-(0,m,m+1,\dots,n_k,m,\boxed{m+1},\dots)
%\end{split}
\]
while the lower horizontal segment satisfies the relation 
\[
%\begin{split}
\underline{b}\ge -(0,n_1,\overline{n_2,\dots,n_k,m,m,m+1})
=-(0,m,m+1,\dots,n_k,m,\boxed{m},m+1,\dots).
%\end{split}
\]
%\foot{SK: edited}
Then $\underline{b}>\bar b$ because the sequences of the corresponding continued fractions with positive entries are in lexicographic order, i.e. the set  $\T^{n_1,n_2,\dots,n_k,m,m}_\bb$ is empty.

Now assume that there exists a non-empty set $\T^{n_1,n_2\dots,n_k, m+1, m+1, \dots, m+1}_\bb$ ($n_k=m$) with the final block of $(m+1)$'s of length $q$ greater than $\frak{l}_{m+1}+1$. The upper vertex of this set satisfies
%\[
%\begin{split}
%\bar b&=[0,\overline{\underbrace{-m,\dots,-m}_p,-(m+1)}]\\
%&=[0,\underbrace{-m,-m,\dots,-m}_p,-(m+1),\underbrace{-m,-m,\dots,-m}_p,-(m+1),\dots]
%\end{split}
%\]
%while the lower horizontal segment is given by 
%\[
%\begin{split}
%\underline{b}&=[0,-m,\overline{-m,\dots,-m,-(m+1),\underbrace{-m,-m\dots, -m}_q,-(m+1)}]\\
%&=[0,\underbrace{-m,-m,\dots,-m}_{p},-(m+1),\underbrace{-m,-m\dots, -m}_q,-(m+1),\underbrace{-m,\dots -m}_{p-1},-(m+1),\dots]\,.
%\end{split}
%\]
%If $p<q$ then the $\bar b,\underline{b}$ coincide on the first $k+p$ entries. Looking at the $k+p+1$ entry, one has $\bar b_{k+p+1}=-(m+1)$ and $\underline{b}_{k+p+1}=-m$, so $\bar b<\underline{b}$. If $p=q$, then the first $k+We assumed that $p\le q$ therefore $\bar b<\underline b$, hence the contradiction. The very same argument shows that the conditions listed are also sufficient, since one has that $\bar b>\underline b$.
\[
\bar b\le -(0,\overline{m,\underbrace{m+1,\dots,m+1}_{\frak{l}_{m+1}},m,\dots,n_k,\underbrace{m+1,\dots,m+1}_q})
%&=[0,-m,\underbrace{-(m+1),\dots,-(m+1)}_p,-m,\dots,m,\undebrace{m+1,dots,m+1}_q,\dots]
\]
while the lower horizontal segment satisfies the relation 
\[
\begin{split}
\underline{b}&\ge -(0,n_1,\overline{n_2,\dots,n_k,m+1})\\
&=-(0,m,{\underbrace{m+1,\dots,m+1}_{\frak{l}_{m+1}},m,\dots, n_k,\underbrace{m+1,\dots,m+1}_{\frak{l}_{m+1}+1},m,\dots}).
%&=-[0,m,\dots,n_k,\underbrace{(m+1),\dots, (m+1)}_{\frak{l}_{m+1}+1},m,\dots,n_k,\dots]
\end{split}
\]
Since the two continued fraction expansions with positive entries coincide on the first $k+\frak{l}_{m+1}+1$ entries and their $k+\frak{l}_{m+1}+2$ entries are $m+1$ and $m$, respectively, we obtain $\bar b<\underline{b}$, i.e. the set $\T^{n_1,n_2\dots,n_k, m+1, m+1, \dots, m+1}_\bb$.

%\foot{SK: added the last case; please, check!}
Finally, suppose that there exists a non-empty set  $\T^{n_1,n_2\dots,n_k, m+1, m+1, \dots, m+1,m}_\bb$ ($n_k=m$) with the final block of $(m+1)$'s of length $q$ less than $\frak{l}_{m+1}$. The upper vertex of this set satisfies
\[
\bar b\le -(0,\overline{m,\underbrace{m+1,\dots,m+1}_{\frak{l}_{m+1}},m,\dots,n_k,\underbrace{m+1,\dots,m+1}_{\frak{l}_{m+1}}})
%&=[0,-m,\underbrace{-(m+1),\dots,-(m+1)}_p,-m,\dots,m,\undebrace{m+1,dots,m+1}_q,\dots]
\]
while the lower horizontal segment satisfies the relation 
\[
\begin{split}
\underline{b}&\ge -(0,n_1,\overline{n_2,\dots,n_k,\underbrace{m+1,\dots,m+1}_{q},m,m+1})\\
&=-(0,m,{\underbrace{m+1,\dots,m+1}_{\frak{l}_{m+1}},m,\dots, n_k,\underbrace{m+1,\dots,m+1}_{q},m,\dots}).
%&=-[0,m,\dots,n_k,\underbrace{(m+1),\dots, (m+1)}_{\frak{l}_{m+1}+1},m,\dots,n_k,\dots]
\end{split}
\]
Since the two continued fraction expansions with positive entries coincide on the first $k+\frak{l}_{m+1}$ entries and their $(k+\frak{l}_{m+1}+1)^{th}$ entries are $(m+1)$ and $m$, respectively, we obtain $\bar b<\underline{b}$, i.e. the set $\T^{n_1,n_2\dots,n_k, m+1, m+1, \dots, m+1,m}_\bb$ is empty.
%A similar argument shows that the length of any block of $(m+1)$'s cannot be smaller than $\frak{l}_{m+1}$.
\end{proof}

%If a set $\E^{(n_i)}_\bb$ is nonempty and $(a,b)\in \E^{(n_i)}_\bb$, then $b$ is uniquely determined by the expansion $b=-\lfloor 0,n_1,n_2,\dots\rceil$. 
In what follows, we describe in an explicit manner the symbolic properties of a sequence $(n_i)$ for which $\E^{(n_i)}_\bb\ne \emptyset$. Notice that in both cases of Lemma \ref{lemE} there are two admissible blocks that can be used to express the admissible sequence $(n_i)$:

\medskip

case (i): $A^{(1)}=(\underbrace{m,\dots,m}_{\frak{l}_{m}},m+1)$ and $B^{(1)}=(\underbrace{m,\dots,m}_{\frak{l}_{m}-1},m+1)$;

\medskip

case (ii): $A^{(1)}=(m,\underbrace{m+1,\dots,m+1}_{\frak{l}_{m+1}})$ and $B^{(1)}=(m,\underbrace{m+1,\dots,m+1}_{\frak{l}_{m+1}+1})$. 

\medskip
\noindent with $\frak{l}_{m}\ge 2$, $\frak{l}_{m+1}\ge 1$. In both situations $A^{(1)}\prec B^{(1)}$. 
One could think of $A^{(1)}$ as being the new `$m$' and $B^{(1)}$ the new `$m+1$', and treat the original sequence of $m$'s and $m+1$'s as a sequence of $A^{(1)}$'s and $B^{(1)}$'s. Furthermore, the next lemma shows that such a substitution process can be continued recursively to construct blocks $A^{(n)}$ and $B^{(n)}$ (for any $n\ge 1$), so that the original sequence $(n_i)$ may be considered to be a sequence of $A^{(n)}$'s and $B^{(n)}$'s. Moreover, only particular blocks of $A^{(n)}$'s and $B^{(n)}$'s warrant non-empty triangular regions of the next generation.%\foot{SK29: formula corrected}

Let us also introduce the notations $A^{(0)}=m$ and $B^{(0)}=m+1$.  Assume that $\E^{(n_i)}_\bb$ is a nonempty set. We have:
\begin{lem}\label{lem-rec}
For every $n\ge 0$, there exist integers  $\frak{l}_{A^{(n)}}\ge 2$, $\frak{l}_{B^{(n)}}\ge 1$ such that the sequence $(n_i)$ can be written as a concatenation of blocks
\begin{equation}\label{case1}
A^{(n+1)}=(\underbrace{A^{(n)},\dots,A^{(n)}}_{\frak{l}_{A^{(n)}}},B^{(n)})  \,,\quad  B^{(n+1)}=(\underbrace{A^{(n)},\dots,A^{(n)}}_{\frak{l}_{A^{(n)}}-1},B^{(n)})
\end{equation}
or 
\begin{equation}\label{case2}
A^{(n+1)}=(A^{(n)},\underbrace{B^{(n)},\dots,B^{(n)}}_{\frak{l}_{B^{(n)}}}) \,,\quad B^{(n+1)}=(A^{(n)},\underbrace{B^{(n)},\dots,B^{(n)}}_{\frak{l}_{B^{(n)}}+1})\,.
\end{equation}
\end{lem}
%\foot{SK8: It would be nice to spell out exactly the conditions that we prove by induction and to add the additional one for $\sigma=(A^{(n)},\tau,A^{(n)})$, but I am afraid it will be too heavy - what do you think?}
\begin{proof}
Notice that Lemma \ref{lemE} proves the above result for $n=0$ with $\frak{l}_{A^{(0)}}=\frak{l}_{m}$, $\frak{l}_{B^{(0)}}=\frak{l}_{m+1}$. 
%We prove the inductive step. 
We show inductively that
\begin{equation}\label{eqorder}
A^{(n)}\prec B^{(n)}
\end{equation}
and if a finite sequence $\sigma$ starts with an $A^{(n)}$ block and ends with a $B^{(n)}$ block, $\sigma=(A^{(n)},\tau,B^{(n)})$, then the lower boundary $\underline b(\sigma)$ of $\T^\sigma_\bb$ (if nonempty) satisfies
\begin{equation}\label{bng}
\underline b(\sigma)\ge-(0,A^{(n)},\overline{\tau,B^{(n)}})\,.
\end{equation}
Relation \eqref{eqorder} is obviously true for $n=0$; \eqref{bng} is also satisfied if $n=0$, since one applies Lemma \ref{lem-cf} part (ii) to the sequence $\tilde\sigma=(A^{(0)},\tau)$ where $\T^{\tilde\sigma}_\bb\supset \T^\sigma_\bb$.

We point out that by applying Lemma \ref{lem-cf} part (i) to the region $\T^{\sigma}$ we have
\begin{equation}\label{ubng}
\bar b(\sigma)\le-(0,\overline{\sigma})=-(0,\overline{A^{(n)},\tau,B^{(n)}})\,.
\end{equation}

To prove the inductive step, suppose that for some $n\ge 1$, we can rewrite the sequence $(n_i)$ using blocks $A^{(n+1)}$ and $B^{(n+1)}$ as in case \eqref{case1} or \eqref{case2}. 

\medskip

\noindent\textbf{Case 1.} Assume $A^{(n+1)}$ and $B^{(n+1)}$ are given by \eqref{case1}.
It follows immediately that $A^{(n+1)}\prec B^{(n+1)}$ since $A^{(n)}\prec B^{(n)}$. Also, if a sequence $\sigma$ starts with an $A^{(n+1)}$ block and ends with a $B^{(n+1)}$ block (thus, implicitly, $\sigma$ starts with an $A^{(n)}$ block and ends with a $B^{(n)}$ block), 
$$\sigma=(A^{(n+1)},\tau,B^{(n+1)})=(\underbrace{A^{(n)},\dots,A^{(n)}}_{\frak{l}_{A^{(n)}}},B^{(n)},\tau,\underbrace{A^{(n)},\dots,A^{(n)}}_{\frak{l}_{A^{(n)}}-1},B^{(n)})\,$$ then, by applying \eqref{bng} to $\tilde \sigma=(\underbrace{A^{(n)},\dots,A^{(n)}}_{\frak{l}_{A^{(n)}}},B^{(n)},\tau)=(A^{(n)},B^{(n+1)},\tau)$ (which starts with $A^{(n)}$ and ends with $B^{(n)}$) we get 
$$\underline b(\sigma)\ge b(\tilde\sigma)\ge -(0,A^{(n)},\overline{B^{(n+1)},\tau})=-(0,A^{(n)},B^{(n+1)},\overline{\tau, B^{(n+1)}})\,.$$%=-(0,A^{(n+1)},\overline{\tau,B^{(n+1)}})$$
Therefore, \eqref{bng} holds for $n+1$, since $(A^{(n)},B^{(n+1)})=A^{(n+1)}$.

Now assume that $(n_i)$ starts with a block of $A^{(n+1)}$'s of length $\frak{l}_{A^{(n+1)}}>1$. We prove that the sequence $(n_i)$ cannot have two consecutive $B^{(n+1)}$'s and any sequence of consecutive blocks $A^{(n+1)}$ has length $\frak{l}_{A^{(n+1)}}$ or $\frak{l}_{A^{(n+1)}-1}$. Suppose the sequence $(n_i)$ contains two consecutive blocks of type  $B^{(n+1)}$: 
\[
(n_i)=(A^{(n+1)},A^{(n+1)},\dots,A^{(n+1)},B^{(n+1)},B^{(n+1)},\dots).
\]
 We look at the set
$$\T^{A^{(n+1)}A^{(n+1)}\dots A^{(n+1)}B^{(n+1)}B^{(n+1)}}$$
 and remark that the upper boundary satisfies (from \eqref{ubng})
\begin{equation}\label{tb}
\bar b\le -(0,\overline{A^{(n+1)},A^{(n+1)},\dots,A^{(n+1)},B^{(n+1)},B^{(n+1)}})
\end{equation}
and the lower boundary satisfies (from \eqref{bng})
\begin{equation}\label{lb}
\underline b\ge -(0,A^{(n+1)},\overline{A^{(n+1)},\dots,A^{(n+1)},B^{(n+1)}})\,.
\end{equation}
%The relation for the upper boundary follows from \eqref{bbar}. To obtain \eqref{lb} one looks at the set
%$$\T^{n_1,n_2,\dots,n_p}\supset \T^{A^{(1)}A^{(1)}\dots A^{(1)}}\supset\T^{A^{(1)}A^{(1)}\dots A^{(1)}B^{(1)}B^{(1)}}$$
%where  $n_p$ is the last $m$ in the last block $A^{(1)}$. Using relation \eqref{underb} one has that
%$$\underline b\ge -(0,n_1,\overline{n_2,\dots,n_p,m+1})= -(0,A^{(1)},\overline{A^{(1)},\dots,A^{(1)},B^{(1)}})$$
%because $A^{(1)}=(m,B^{(1)})$. 
But \eqref{tb} and \eqref{lb} imply that $\underline b>\bar b$, because the two corresponding continued fractions 
with positive entries are %\foot{SK:edited}%have sequences 
in lexicographic order. Thus, there cannot be two consecutive $B^{(n+1)}$ blocks in the sequence $(n_i)$. 

Now, let us check that the sequence $(n_i)$ cannot have a block of $A^{(n+1)}$'s of length $q>\frak{l}_{A^{(n+1)}}$. Assume the contrary,%\foot{SK: inserted $B^{(1)}$ before $\dots$}
$$(n_i)=(\underbrace{A^{(n+1)},\dots,A^{(n+1)}}_{\frak{l}_{A^{(n+1)}}},B^{(n+1)},\tau,B^{(n+1)}, \underbrace{A^{(n+1)},\dots,A^{(n+1)}}_{q},B^{(n+1)},\dots)\,.$$
Then the set $\T^{(n_i)}_\bb$ has the upper bound $\bar b$ satisfying
$$\bar b\le -(0,\overline{\underbrace{A^{(n+1)},\dots,A^{(n+1)}}_{\frak{l}_{A^{(n+1)}}},B^{(n+1)},\tau,B^{(n+1)}})$$
while the lower bound $\underline b$ satisfies by \eqref{bng}
$$
\underline b\ge -(0,A^{(n+1)},\overline{\underbrace{A^{(n+1)},\dots,A^{(n+1)}}_{\frak{l}_{A^{(n+1)}}-1},B^{(n+1)},\tau,B^{(n+1)},\underbrace{A^{(n+1)},\dots,A^{(n+1)}}_{q},B^{(n+1)}}).
$$
%\foot{SK: in the above formula replaced $\frak{l}^{(1)}$ by $\frak{l}^{(1)}-1$}
Comparing the two continued fractions, we get that $\bar b<\underline b$ (since $A^{(n+1)}\prec B^{(n+1)}$ and $q>\frak{l}_{A^{(n+1)}}$).

Now assume that $(n_i)$ starts with $A^{(n+1)}$ and then continues with a block of $B^{(n+1)}$'s of length $\frak{l}_{B^{(n+1)}}\ge 1$. We prove that the sequence $(n_i)$ cannot have two consecutive $A^{(n+1)}$'s and any sequence of consecutive blocks $B^{(n+1)}$ has length $\frak{l}_{B^{(n+1)}}$ or $\frak{l}_{B^{(n+1)}}+1$. Suppose the sequence $(n_i)$ contains two (or more) consecutive blocks of type  $A^{(n+1)}$: 
$$(n_i)=(A^{(n+1)},B^{(n+1)},\tau,B^{(n+1)},\underbrace{A^{(n+1)},\dots,A^{(n+1)}}_{q\ge 2},B^{(n+1)},\dots)\,.$$ 
We study the  region $\T^{A^{(n+1)},B^{(n+1)},\tau,B^{(n+1)},A^{(n+1)},\dots,A^{(n+1)},B^{(n+1)}}$
and remark that its upper boundary satisfies (from \eqref{ubng})
\begin{equation}\label{ttb}
\bar b\le -(0,\overline{A^{(n+1)},B^{(n+1)},\tau,B^{(n+1)}})
\end{equation}
and the lower boundary satisfies (from \eqref{bng})
\begin{equation}\label{llb}
\underline b\ge -(0,A^{(n+1)},\overline{B^{(n+1)},\tau,B^{(n+1)},\underbrace{A^{(n+1)},\dots,A^{(n+1)}}_{q\ge 2},B^{(n+1)}}).
\end{equation}
%The relation for the upper boundary follows from \eqref{bbar}. To obtain \eqref{lb} one looks at the set
%$$\T^{n_1,n_2,\dots,n_p}\supset \T^{A^{(1)}A^{(1)}\dots A^{(1)}}\supset\T^{A^{(1)}A^{(1)}\dots A^{(1)}B^{(1)}B^{(1)}}$$
%where  $n_p$ is the last $m$ in the last block $A^{(1)}$. Using relation \eqref{underb} one has that
%$$\underline b\ge -(0,n_1,\overline{n_2,\dots,n_p,m+1})= -(0,A^{(1)},\overline{A^{(1)},\dots,A^{(1)},B^{(1)}})$$
%because $A^{(1)}=(m,B^{(1)})$. 
But \eqref{ttb} and \eqref{llb} implie that $\underline b>\bar b$ because the two corresponding continued fractions 
with positive entries are %\foot{SK:edited}%have sequences 
in lexicographic order. Thus, there cannot be two consecutive $A^{(n+1)}$ blocks in the sequence $(n_i)$. 

Now, let us check that the sequence $(n_i)$ cannot have a block of $B^{(n+1)}$'s of length $q>\frak{l}_{B^{(n+1)}}+1$. Assume the contrary,%\foot{SK: inserted $B^{(1)}$ before $\dots$}
$$(n_i)=(A^{(n+1)},\underbrace{B^{(n+1)}\dots,B^{(n+1)}}_{\frak{l}_{B^{(n+1)}}},A^{(n+1)},\tau,A^{(n+1)}, \underbrace{B^{(n+1)},\dots,B^{(n+1)}}_{q},A^{(n+1)},\dots)\,.$$
Then the set $\T^{(n_i)}$ has the upper bound $\bar b$ satisfying
$$\bar b\le -(0,\overline{A^{(n+1)},\underbrace{B^{(n+1)}\dots,B^{(n+1)}}_{\frak{l}_{B^{(n+1)}}},A^{(n+1)},\tau,A^{(n+1)}, \underbrace{B^{(n+1)},\dots,B^{(n+1)}}_{q},A^{(n+1)}})
$$
while the lower bound $\underline b$ satisfies by \eqref{bng}
$$
\underline b\ge -(0,A^{(n+1)},\overline{\underbrace{B^{(n+1)}\dots,B^{(n+1)}}_{\frak{l}_{B^{(n+1)}}},A^{(n+1)},\tau,A^{(n+1)}, B^{(n+1)}}).
$$
%\foot{SK: in the above formula replaced $\frak{l}^{(1)}$ by $\frak{l}^{(1)}-1$}
Comparing the two continued fractions, we get that $\bar b<\underline b$.

\medskip

\noindent\textbf{Case 2.} Assume $A^{(n+1)}$ and $B^{(n+1)}$ are given by \eqref{case2}. It follows that $A^{(n+1)}\prec B^{(n+1)}$ since $A^{(n+1)}$ is the beginning block of $B^{(n+1)}$.
%\foot{SK30: isn't this automatic since $A^{(n+1)}$ is the beginning block of $B^{(n+1)}$?}
Also, if a sequence $\sigma$ starts with an $A^{(n+1)}$ block and ends with a $B^{(n+1)}$ block (thus, implicitly, $\sigma$ starts with an $A^{(n)}$ block and ends with a $B^{(n)}$ block), 
$$\sigma=(A^{(n+1)},\tau,B^{(n+1)})=(A^{(n)},\underbrace{B^{(n)},\dots,B^{(n)}}_{\frak{l}_{B^{(n)}}},\tau,A^{(n)},\underbrace{B^{(n)},\dots,B^{(n)}}_{\frak{l}_{B^{(n)}}+1})\,$$ 
then by applying \eqref{bng} to $\tilde \sigma=(A^{(n)},\underbrace{B^{(n)},\dots,B^{(n)}}_{\frak{l}_{B^{(n)}}},\tau,A^{(n)},B^{(n)})$, which starts with $A^{(n)}$ and ends with $B^{(n)}$, we get 
\begin{equation*}
\begin{split}
\underline b(\sigma)\ge \underline b(\tilde\sigma)\ge -(0,A^{(n)},\overline{\underbrace{B^{(n)},\dots,B^{(n)}}_{\frak{l}_{B^{(n)}}},\tau,A^{(n)},B^{(n)}})\\=-(0,A^{(n+1)},\overline{\tau, A^{(n)},\underbrace{B^{(n)},\dots,B^{(n)}}_{\frak{l}_{B^{(n)}+1}}})
\end{split}
\end{equation*}
%=-(0,A^{(n+1)},\overline{\tau,B^{(n+1)}})$$
so \eqref{bng} holds for $n+1$.

Assume that $(n_i)$ starts with a sequence of $A^{(n+1)}$'s of length $\frak{l}_{A^{(n+1)}}>1$. Similar to the analysis of the first case, one proves that the sequence $(n_i)$ cannot have two consecutive $B^{(n+1)}$'s and any sequence of consecutive blocks $A^{(n+1)}$ has length $\frak{l}_{A^{(n+1)}}$ or $\frak{l}_{A^{(n+1)}}-1$. 

If the sequence $(n_i)$ starts with $A^{(n+1)}$ and then continues with a sequence of $B^{(n+1)}$'s of length $\frak{l}_{B^{(n+1)}}\ge 1$, one can prove that the sequence $(n_i)$ cannot have two consecutive $A^{(n+1)}$'s and any sequence of consecutive blocks $B^{(n+1)}$ has length $\frak{l}_{B^{(n+1)}}$ or $\frak{l}_{B^{(n+1)}}+1$. 
\end{proof}
%\foot{SK30: should I trust you on the end of the proof?}
Additionally, we prove
\begin{lem}\label{lem-order}
%We have the relation $A^{(n)}\prec B^{(n)}$ for any $n\ge 1$. 
If the block $\tau_1=(n_i,\dots, n_l)$ is a tail of $A^{(n)}$ and $\tau_2=(p_j,\dots,p_h)$ is a tail of $B^{(n)}$, then
$A^{(n)}\prec \tau_1$ and $B^{(n)}\prec \tau_2$.
\end{lem}
\begin{proof}
%The first relation can be easily proved by an inductive argument. We provide an inductive proof for the second property. 
The statement is obviously true if $n=1$. Assume it is true for some $n$ both for $A^{(n)}$ and $B^{(n)}$. We analyze the case of $A^{(n+1)}$ being given by \eqref{case1}, $A^{(n+1)}=(\underbrace{A^{(n)},\dots, A^{(n)}}_{\frak{l}_{A^{(n)}}},B^{(n)})$. Consider an arbitrary tail $\tau$ of $A^{(n+1)}$; $\tau$ could start with a block $A^{(n)}$ or a tail of $A^{(n)}$ or $\tau$ coincides with $B^{(n)}$ or a tail of $B^{(n)}$. In all situations, the inductive hypothesis and the fact that $A^{(n)}\prec B^{(n)}$ prove that $A^{(n+1)}\prec \tau$.  The case of $A^{(n+1)}$ given by \eqref{case2} is treated similarly.
\end{proof}

\begin{rem}\label{rem-ineq}
Using the relations \eqref{bng} and \eqref{ubng}, notice that a set $\T^{A^{(n+1)}}_\bb$ (if nonempty) has the upper vertex satisfying
\begin{equation}\label{bn}
\bar b_{n+1}\le -(0,\overline{A^{(n+1)}})
\end{equation}
and a lower horizontal boundary that satisfies
\begin{equation}\label{ubn1}
\underline b_{n+1}\ge -(0,A^{(n+1)},\overline{B^{(n+1)}})
\end{equation}
if $A^{(n+1)}$ is given by the substitution rule \eqref{case1}, and 
\begin{equation}\label{ubn2}
\underline b_{n+1}\ge -(0,A^{(n)},\overline{B^{(n)}})
\end{equation}
if $A^{(n+1)}$ is given by \eqref{case2}.
\end{rem}
%Indeed, the inequality for $\bar b(n)$ is obtained by applying Lemma \ref{lem-cf} part (i) to the sequence $A^{(n)}$. The relation for $\underline b(n)$ is derived from the inductive procedure (see \eqref{bng}) where $\tau$ and the second $A^{(n)}$ are absent). \foot{SK: I still do not quite understand  the argument
%``we apply Lemma \ref{lem-cf} part (ii) to the starting block $\sigma(n)$ of $A^{(n)}$ obtained by eliminating all the $B$-type blocks from the right end of $A^{(n)}$.", and am quite satisfied by a reference to (\eqref{bng}); there is something here I do not understand: if you eliminate all $B$-type blocks, then your block will not end by $A^{(n)}$}
%Also, if we let $\tau(n)$ be the eliminated end block, then it follows from Lemma \ref{lem-order} that $\sigma\prec\tau(n)$ for any starting block $\sigma$ of $A^{(n)}$.
%\foot{SK6: I think $\sigma\prec\tau(n)$ follows simply because $\sigma$ starts with $m$ while $\tau(n)$ starts with $m+1$ - do you agree? But in Thm 8.8 you use that $\tau\prec\tau(n)$ for any tail $\tau$ of $A^{(n)}$, and this is true still, but the proof is less trivial - do you want to include it? }

We will prove that the above inequalities are actually equality relations. For that we construct a starting subsequence of $A^{(n+1)}$ defined inductively as:

$$
\sigma^{(1)}=
\begin{cases}
(\underbrace{m,\dots,m}_{{\frak l}_m})  & \text{if } A^{(1)}=(\underbrace{m,\dots,m}_{{\frak l}_m},m+1)\\
({m})  & \text{if } A^{(1)}=(m,\underbrace{m+1,\dots,m+1}_{{\frak l}_{m+1}})\\
\end{cases}
$$
\noindent{Case 1.} If $A^{(n)}$ is given by a relation of type \eqref{case1}, i.e. $A^{(n)}=(A^{(n-1)}, \dots,A^{(n-1)}, B^{(n-1)})$, then
\begin{equation}\label{snc1}
\sigma^{(n+1)}=
\begin{cases}
(\underbrace{A^{(n)},\dots,A^{(n)}}_{{\frak l}_A^{(n)}-1},\sigma^{(n)})  & \text{if } A^{(n+1)}=(\underbrace{A^{(n)},\dots,A^{(n)}}_{{\frak l}_{A^{(n)}}},B^{(n)})\\
\sigma^{(n)}  & \text{if } A^{(n+1)}=(A^{(n)},\underbrace{B^{(n)},\dots,B^{(n)}}_{{\frak l}_{B^{(n)}}})
\end{cases}
\end{equation}

\noindent{Case 2.} If $A^{(n)}$ is given by a relation of type \eqref{case2}, i.e. $A^{(n)}=(A^{(n-1)}, B^{(n-1)}, \dots, B^{(n-1)})$, then
\begin{equation}\label{snc2}
\sigma^{(n+1)}=
\begin{cases}
(\underbrace{A^{(n)},\dots,A^{(n)}}_{{\frak l}_A^{(n)}},\sigma^{(n)})  & \text{if } A^{(n+1)}=(\underbrace{A^{(n)},\dots,A^{(n)}}_{{\frak l}_{A^{(n)}}},B^{(n)})\\
(A^{(n)},\sigma^{(n)})  & \text{if } A^{(n+1)}=(A^{(n)},\underbrace{B^{(n)},\dots,B^{(n)}}_{{\frak l}_{B^{(n)}}})
\end{cases}
\end{equation}
%Let $\sigma(n)$ be the starting block of $A^{(n)}$ obtained by eliminating all the $B$-type blocks from the right end of $A^{(n)}$. For the case (i) substitution rule we have
%\[
%\begin{split}
%A^{(n)}=&(\sigma(n),B^{(0)},B^{(1)}, B^{(2)},\dots ,B^{(n-2)},B^{(n-1)})\\
%=&(\underbrace{A^{(n-1)},\dots,A^{(n-1)}}_{\frak{I}^{(n-1)}-1},\sigma(n-1),B^{(0)},B^{(1)}, B^{(2)},\dots ,B^{(n-2)},B^{(n-1)}).
%\end{split}
%\]
We introduce the notation $f^{\sigma}$ to denote the transformation $T^{n_k}S\dots T^{n_1}S$ if $\sigma=(n_1,\dots,n_k)$.
\begin{lem}\label{relbn}
Let $\sigma^{(n+1)}$ be the starting block of $A^{(n+1)}$ defined as above. Then the equation 
$$f^{\sigma^{(n+1)}}b=\frac{b}{b+1}$$ has a unique solution $b\in [0,1]$ given by
\begin{equation}\label{eqbn}
b_{n+1}=\begin{cases}
 -(0,A^{(n+1)},\overline{B^{(n+1)}}) & \text{if } A^{(n+1)} \text{ given by } \eqref{case1}\\
  -(0,A^{(n)},\overline{B^{(n)}}) & \text{if } A^{(n+1)} \text{ given by } \eqref{case2}
\end{cases}
\end{equation}
%Then if a point $b\in[0,1]$ satisfies $f^{\sigma(n)}b=\frac{b}{b+1}$ and
%$b\in\T_{\bb}^{\sigma(n)}$, then $b=-(0,A^{(n)},\overline{B^{(n)}})$.
\end{lem}
\begin{proof} 
We proceed with an inductive proof, and as part of it we also show that 
\begin{equation}\label{sigma}
(\sigma^{(n+1)},m+1,\tilde A^{(n)})=
\begin{cases}
A^{(n+1)} & \text{if } A^{(n+1)}=(\underbrace{A^{(n)},\dots,A^{(n)}}_{{\frak l}_{A^{(n)}}},B^{(n)})\\
(A^{(n)},B^{(n)}) & \text{if } A^{(n+1)}=(A^{(n)},\underbrace{B^{(n)},\dots,B^{(n)}}_{{\frak l}_{B^{(n)}}})
\end{cases}
\end{equation}
where $A^{(n)}=(m,\tilde A^{(n)})$.

The relation \eqref{eqbn} is true for $n=0$ due to Lemma \ref{lem-cf}(ii). Also, \eqref{sigma} follows immediately.
%\foot{SK3: references fixed}
Suppose now that the inductive relations hold for some $n$. We analyze the solution of $f^{\sigma^{(n+2)}}b=\frac{b}{b+1}$. 

Assume that $A^{(n+1)}=(\underbrace{A^{(n)},\dots,A^{(n)}}_{{\frak l}_{A^{(n)}}}, B^{(n)})$. We look at the two possible cases:

\noindent (i)  If $A^{(n+2)}=(\underbrace{A^{(n+1)},\dots,A^{(n+1)}}_{{\frak l}_{A^{(n+1)}}}, B^{(n+1)})$, 
$\sigma^{(n+2)}=(\underbrace{A^{(n+1)},\dots,A^{(n+1)}}_{{\frak l}_{A^{(n+1)}}-1},\sigma^{(n+1)})$.
Using Lemma \ref{lem-cf}(ii), we have that the solution to $f^{\sigma^{(n+2)}}b=\frac{b}{b+1}$ is given by
%\foot{SK3: brace fixed}
\begin{equation*}
\begin{split}
b_{n+2}&=-(0,m,\overline{\tilde A^{(n+1)},\underbrace{A^{(n+1)},\dots,A^{(n+1)}}_{{\frak l}_{A^{(n+1)}}-2},\sigma^{(n+1)},m+1})\\
&= -(0,m,\tilde A^{(n+1)},\overline{\underbrace{A^{(n+1)},\dots,A^{(n+1)}}_{{\frak l}_{A^{(n+1)}}-2},\sigma^{(n+1)},m+1,\tilde A^{(n+1)}})\\
&=-(0,m,\tilde A^{(n+1)},\overline{\underbrace{A^{(n+1)},\dots,A^{(n+1)}}_{{\frak l}_{A^{(n+1)}}-2},\sigma^{(n+1)},m+1,\tilde A^{(n)},B^{(n+1)}})\\
&=-(0,A^{(n+1)},\overline{\underbrace{A^{(n+1)},\dots,A^{(n+1)}}_{{\frak l}_{A^{(n+1)}}-2},A^{(n+1)},B^{(n+1)}}\\
&=-(0,A^{(n+1)},\overline{B^{(n+2)}}=-(0,A^{(n+2)},\overline{B^{(n+2)}}).
\end{split}
\end{equation*}
Also,
\begin{equation*}
\begin{split}
(\sigma^{(n+2)},m+1,\tilde A^{(n+1)})&=(\underbrace{A^{(n+1)},\dots,A^{(n+1)}}_{{\frak l}_{A^{(n+1)}}-1},\sigma^{(n+1)},m+1,\tilde A^{(n)}, B^{(n+1)})\\
&=(\underbrace{A^{(n+1)},\dots,A^{(n+1)}}_{{\frak l}_{A^{(n+1)}}-1},A^{(n+1)},B^{(n+1)})=A^{(n+2)}.
\end{split}
\end{equation*}

\noindent (ii)  If $A^{(n+2)}=(A^{(n+1)},\underbrace{B^{(n+1)},\dots,B^{(n+1)}}_{{\frak l}_{B^{(n+1)}}})$, then
$
\sigma^{(n+2)}=\sigma^{(n+1)}\,,
$
and the induction step gives us the solution of $f^{\sigma(n+2)}b=\frac{b}{b+1}$ as
$
b_{n+2}=-(0,A^{(n+1)},\overline{B^{(n+1)}})
$. Also,
\begin{equation*}
(\sigma^{(n+2)},m+1,\tilde A^{(n+1)})=(\sigma^{(n+1)},m+1,\tilde A^{(n)},B^{(n+1)})=(A^{(n+1)},B^{(n+1)}).
\end{equation*}

\noindent Now assume that $A^{(n+1)}=(A^{n},\underbrace{B^{(n)},\dots,B^{(n)}}_{{\frak l}_{B^{(n)}}})$. 
%\foot{SK3: upper index changed}
We look again at the two possible cases:

\noindent (i)  If $A^{(n+2)}=(\underbrace{A^{(n+1)},\dots,A^{(n+1)}}_{{\frak l}_{A^{(n+1)}}}, B^{(n+1)})$,
$
\sigma^{(n+2)}=(\underbrace{A^{(n+1)},\dots,A^{(n+1)}}_{{\frak l}_{A^{(n+1)}}},\sigma^{(n+1)})
$.
Using Lemma \ref{lem-cf}(ii), we have that the solution to $f^{\sigma(n+2)}b=\frac{b}{b+1}$ is given by
\begin{equation*}
\begin{split}
b_{n+2}&=-(0,m,\overline{\tilde A^{(n+1)},\underbrace{A^{(n+1)},\dots,A^{(n+1)}}_{{\frak l}_{A^{(n+1)}}-1},\sigma^{(n+1)},m+1})\\
&= -(0,m,\tilde A^{(n+1)},\overline{\underbrace{A^{(n+1)},\dots,A^{(n+1)}}_{{\frak l}_{A^{(n+1)}}-1},\sigma^{(n+1)},m+1,\tilde A^{(n+1)}})\\
&=-(0,m,\tilde A^{(n+1)},\overline{\underbrace{A^{(n+1)},\dots,A^{(n+1)}}_{{\frak l}_{A^{(n+1)}}-1},\sigma^{(n+1)},m+1,\tilde A^{(n)},\underbrace{B^{(n)},\dots,B^{(n)}}_{{\frak l}_{B^{(n)}}})}\\
&=-(0,m,\tilde A^{(n+1)},\overline{\underbrace{A^{(n+1)},\dots,A^{(n+1)}}_{{\frak l}_{A^{(n+1)}}-1},A^{(n)},B^{(n)},\underbrace{B^{(n)},\dots,B^{(n)}}_{{\frak l}_{B^{(n)}}})}\\
&=-(0,A^{(n+1)},\overline{\underbrace{A^{(n+1)},\dots,A^{(n+1)}}_{{\frak l}_{A^{(n+1)}}-1},B^{(n+1)}})=-(0,A^{(n+2)},\overline{B^{(n+2)}}).
\end{split}
\end{equation*}
A similar approach gives us that $(\sigma^{(n+2)},m+1,\tilde A^{(n+1)})=A^{(n+2)}$.

\medskip

\noindent (ii)  If $A^{(n+2)}=(A^{(n+1)},\underbrace{B^{(n+1)},\dots,B^{(n+1)}}_{{\frak l}_{B^{(n+1)}}})$, then
$
\sigma^{(n+2)}=(A^{(n+1)},\sigma^{(n+1)}).
$
Using Lemma \ref{lem-cf}(ii),
%\foot{SK3: ref changed to (i)}, 
we have that the solution to $f^{\sigma(n+2)}b=\frac{b}{b+1}$ is given by

\begin{equation*}
\begin{split}
b_{n+2}&=-(0,m,\overline{\tilde A^{(n+1)},\sigma^{(n+1)},m+1})\\
&= -(0,m,\tilde A^{(n+1)},\overline{\sigma^{(n+1)},m+1,\tilde A^{(n+1)}})\\
&=-(0,m,\tilde A^{(n+1)},\overline{\sigma^{(n+1)},m+1,\tilde A^{(n)},\underbrace{B^{(n)},\dots,B^{(n)}}_{{\frak l}_{B^{(n)}}}})\\
&=-(0,m,\tilde A^{(n+1)},\overline{A^{(n)},B^{(n)},\underbrace{B^{(n)},\dots,B^{(n)}}_{{\frak l}_{B^{(n)}}}})\\
&=-(0,A^{(n+1)},\overline{B^{(n+1)}}).
\end{split}
\end{equation*}
Also, 
\begin{equation*}
\begin{split}
(\sigma^{(n+2)},m+1,\tilde A^{(n+1)})&=(A^{(n+1)},\sigma^{(n+1)},m+1,\tilde A^{(n)},\underbrace{B^{(n)},\dots,B^{(n)}}_{{\frak l}_{B^{(n)}}})\\
&=(A^{(n+1)},A^{n},B^{(n)},\underbrace{B^{(n)},\dots,B^{(n)}}_{{\frak l}_{B^{(n)}}})=(A^{(n+1)},B^{(n+1)}).
\end{split}
\end{equation*}

%Recall that $A^{(0)}=m,\,A^{(0)}=m+1$. By the substitution rule we have 
%\[
%A^{(n-1)}=(A^{(n-2)},B^{(n-1)})=\cdots =(m,B^{(1)},\dots ,B^{(n-1)})=(A^{(0)},B^{(1)}, \dots ,B^{(n-1)}).
%\]
%An application of Lemma \ref{lem-cf} part (ii) to the block $\sigma(n)$ gives
%\[
%\begin{split}
%b=&-(0,m,\overline{\sigma(n),(m+1)}\\
%=&-(0, A^{(0)}, \overline{B^{(1)}, B^{(2)},\dots ,B^{(n-1)},\underbrace{A^{(n-1)},\dots,A^{(n-1)}}_{\frak{I}^{(n-1)}-2},\dots ,B^{(0)}})\\
%=&-(0,A^{(n-1)},\overline{\underbrace{A^{(n-1)},\dots,A^{(n-1)}}_{\frak{I}^{(n-1)}-2},\dots,B^{(0)},B^{(1)}, B^{(2)},\dots ,B^{(n-2)},B^{(n-1)}})\\
%=&-(0,A^{(n-1)},\overline{\underbrace{A^{(n-1)},\dots,A^{(n-1)}}_{\frak{I}^{(n-1)}-1}, B^{(n-1)}},\\
%=&-(0,A^{(n-1)},\overline{B^{(n)}})=-(0,A^{(n)},\overline{B^{(n)}}).
%\end{split}
%\]
\end{proof}

\begin{thm}\label{thm-exact}
Any sequence $(n_i)$ constructed recursively using relations \eqref{case1} and \eqref{case2} provides a non-empty set $\mathcal E^{(n_i)}_b$.
\end{thm}

\begin{proof}
We prove inductively that any set $\T^{A^{(n+1)}}_\bb$ is nonempty and the relations \eqref{bn} and \eqref{ubn1} or \eqref{ubn2} 
are actual equalities, i.e.
\begin{equation}\label{bne}
\bar b_{n+1}= -(0,\overline{A^{(n+1)}})
\end{equation}
and a lower horizontal boundary that satisfies
\begin{equation}\label{ube1}
\underline b_{n+1}= -(0,A^{(n+1)},\overline{B^{(n+1)}})
\end{equation}
if $A^{(n+1)}$ is given by the substitution rule \eqref{case1} or 
\begin{equation}\label{ube2}
\underline b_{n+1}= -(0,A^{(n)},\overline{B^{(n)}})
\end{equation}
if $A^{(n+1)}$ is given by \eqref{case2}. As part of the inductive proof, we also show that any tail block $\tau$ of $A^{(n+1)}$, $\tau\ne \tau^{(n+1)}$  satisfies $\tau\prec \tau^{(n+1)}$, where $\tau^{(n+1)}$ denotes the tail block of $A^{(n+1)}$ obtained by eliminating the starting block $\sigma^{(n+1)}$ defined by \eqref{snc1} or \eqref{snc2}.

%if $A^{(n+1)}=(\sigma ^{(n+1)},\tau^{(n+1)})=(\sigma, \tau)$, then $\tau\prec \tau^{(n+1)}$.

%\begin{equation}\label{eqbnbar}
%\bar b(n)=-(0,\overline{A^{(n)}})\,,\quad
%\underline b(n)=-(0,A^{(n)},\overline{B^{(n)}})\,.
%\end{equation}
%Moreover, the value $\underline b(n)$ is obtained by writing relation \eqref{lb} to the starting block \foot{SK: replaced ``subsequence" by ``block". May be change ``starting block" to "head" as opposed to tail? Then need to give both definitions.} $\sigma(n)$ of $A^{(n)}$ obtained by eliminating all the $B$-type blocks from the right end of $A^{(n)}$. Let $\tau(n)$ be the eliminated tail. It is easy to see that $\sigma\prec\tau(n)$ for any starting block $\sigma$ of $A^{(n)}$.

Indeed for $n=0$, one can check directly that the sets $\T^{m,m,\dots,m,m+1}_\bb$ and $\T^{m,m+1,\dots,m+1}_\bb$ satisfy the above equalities using the fact that an ``$m$" digit does not change the position of the upper vertex, while an  ``$m+1$'' digit does not change the position of the horizontal segment of such a triangular set. Also, for any tail $\tau\ne \tau^{(1)}$ of $A^{(1)}$,  $\tau\prec \tau^{(1)}$. 

Now, let us assume that $\T_\bb^{A^{(n+1)}}$ obtained from $A^{(n+1)}=(\underbrace{A^{(n)},\dots,A^{(n)}}_{\frak{l}^{(n)}},B^{(n)})$ is nonempty and satisfies \eqref{bne} and \eqref{ube1}.
%Suppose that $A^{(n+1)}$ is given by \eqref{case1}, $A^{(n+1)}=(\underbrace{A^{(n)},\dots,A^{(n)}}_{\frak{l}^{(n)}},B^{(n)})$. %The other situation described by \eqref{case2} can be treated similarly. 
For $\T^{A^{(n+2)}}_\bb$ we look at the two possible cases:

\noindent (i)  $A^{(n+2)}=(\underbrace{A^{(n+1)},\dots,A^{(n+1)}}_{{\frak l}_{A^{(n+1)}}}, B^{(n+1)})$. By Remark \ref{rem-ineq}, 
\[
\bar b_{n+2}\le -(0,\overline{A^{(n+2)}})=-(0,\overline{\underbrace{A^{(n+1)},\dots,A^{(n+1)}}_{\frak{l}_{A^{(n+1)}}},B^{(n+1)}})=:\hat b
\]
and 
\[
\underline b_{n+2}\ge -(0,A^{(n+2)},\overline{B^{(n+2)}})=-(0,A^{(n+1)},\overline{\underbrace{A^{(n+1)},\dots,A^{(n+1)}}_{\frak{l}_{A^{(n+1)}-1}},B^{(n+1)}})=:\tilde b
\]
where $\tilde b$ was obtained by applying Lemma \ref{lem-cf} part (ii) to the starting block 
$$
\sigma^{(n+2)}=(\underbrace{A^{(n+1)},\dots,A^{(n+1)}}_{{\frak l}_{A^{(n+1)}}-1},\sigma^{(n+1)})
$$
of $A^{(n+2)}$.
%We have $\tilde b\le \underline b_{n+2}$, $\bar b_{n+2}\le \hat b$. 

%We prove that any point $(b-1,b)$ with $\tilde b<b<\hat b$ is part of the set $\T^{A^{(n+2)}}$.

%We introduce the notation $f^{\sigma}$ to denote the transformation $T^{n_k}S\dots T^{n_1}S$ if $\sigma=(n_1,\dots,n_k)$.

We prove first the other inductive step: any tail block $\tau$ of $A^{(n+2)}$, $\tau\ne\tau^{(n+2)}$ satisfies $\tau\prec \tau^{(n+2)}$. Notice that $\tau^{(n+2)}=(\tau^{(n+1)}, B^{(n+1)})$.  There exists $\tau'$ a tail block of $A^{(n+1)}$ with the property that
$$\tau=(\tau',\underbrace{A^{(n+1)},\dots,A^{(n+1)}}_l,B^{n+1}), \quad 0\le l\le \frak{l}_{A^{(n+1)}-1}$$ %\foot{SK6: indices fixed}
or $\tau=\tau'$.
%Indeed, the former holds when the tail block starts within the $(A^{(n+1)},\dots, A^{(n+1)})$ block of $A^{(n+2)}$ , 
The latter case holds when $\tau$ is just a tail of $B^{(n+1)}$ (which itself is a tail of $A^{(n+1)}$). It is possible that $\tau'=\emptyset$, but in this case  $\tau\prec\tau^{(n+2)}$ %\foot{SK6: fixed $\prec$}
because $A^{(n+1)}\prec \tau^{(n+1)}$ by Lemma \ref{lem-order}. If $\tau'\ne\emptyset$, we also get %\foot{SK6: extra words deleted}
 that $\tau\prec \tau^{(n+2)}$ by using the inductive hypothesis relation $\tau'\prec \tau^{(n+1)}$.

Now we show that the points $(\tilde b-1,\tilde b)$ and $(\hat b-1,\hat b)$ belong to the set $\T^{A^{(n+2)}}_\bb$. The point $(\hat b-1,\hat b)$  belongs to $\T^{A^{(n+1)}}_\bb$ so $f^{A^{(n+1)}}\hat b\le \hat b$. If $\sigma$ is an intermediate block between $A^{(n+1)}$ and $A^{(n+2)}$, $A^{(n+1)}\subset\sigma\subset A^{(n+2)}$, then %\foot{IU: added relation}
$$f^\sigma(\hat b)=-(0,\tau,\overline{A^{(n+2)}})\le -(0,\overline{A^{(n+2)}})=\hat b$$
The inequality is due to the fact that $\tau$ is a tail block of $A^{(n+2)}$ obtained by eliminating $\sigma$, so
$A^{(n+2)}\prec \tau$.

Now we show that $f^\sigma(\tilde b)\ge \tilde b/(\tilde b+1)$  for any intermediate block $\sigma$ between $A^{(n+1)}$ and $A^{(n+2)}$.  We have that $f^{\sigma^{(n+2)}}(\tilde b)=\tilde b/(\tilde b+1)$ by Lemma \ref{relbn}, and 
%\foot{SK6: space inserted}
\[
f^{\sigma^{(n+2)}}(\tilde b)=-(0,\tau^{(n+2)},\overline{B^{(n+2)}})\, ,
\]
where $\tau^{(n+2)}=(\tau^{(n+1)},B^{(n+1)})$. Also $f^\sigma(\tilde b)=-(0,\tau,\overline{B^{(n+2)}})$ with $\tau$ being the tail block of $A^{(n+2)}$ obtained by eliminating $\sigma$. But $\tau\prec\tau^{(n+2)}$ as we have just proved, hence $f^\sigma(\tilde b)\ge f^{\sigma(n+2)}(\tilde b)$. 

In conclusion, any intermediate block 
%\foot{SK6: replaced ``sequence" by ``block"} 
$\sigma$ between $A^{(n+1)}$ and $A^{(n+2)}$ satisfies
\[
\tilde b/(\tilde b+1)\le f^\sigma(\tilde b)\le f^\sigma(\hat b)\le \hat b\,,
\]
therefore the points $(\tilde b-1,\tilde b)$ and $(\hat b-1,\hat b)$ belong to the intermediate set $\T^\sigma_\bb$. This proves the induction step for $\T^{A^{(n+2)}}_\bb$.

\medskip

\noindent (ii)   $A^{(n+2)}=(A^{(n+1)},\underbrace{B^{(n+1)},\dots,B^{(n+1)}}_{{\frak l}_{B^{(n+1)}}})$. By Remark \ref{rem-ineq}, we have that 
\[
\bar b_{n+2}\le  -(0,\overline{A^{(n+2)}})=-(0,\overline{A^{(n+1)},\underbrace{B^{(n+1)},\dots,B^{(n+1)}}_{\frak{l}_{B^{(n+1)}}}})=:\hat b
\]
and 
\[
\underline b_{n+2}\ge -(0,A^{(n+1)},\overline{B^{(n+1)}})=:\tilde b
\]
where $\tilde b$ was obtained by applying Lemma \ref{lem-cf} part (ii) to the starting block 
$
\sigma^{(n+2)}=\sigma^{(n+1)}
$
of $A^{(n+2)}$.

We prove first the other inductive step: any tail block $\tau$ of $A^{(n+2)}$, $\tau\ne\tau^{(n+2)}$, satisfies $\tau\prec \tau^{(n+2)}$. There exists $\tau'$ a tail block of $A^{(n+1)}$ with the property that 
$$\tau=(\tau',\underbrace{B^{(n+1)},\dots,B^{(n+1)}}_l), \quad 0\le l\le \frak{l}_{B^{(n+1)}}$$ %\foot{SK6: indices fixed}
(again, using the fact that $B^{(n+1)}$ is a tail block of $A^{(n+1)}$). 
Since 
$$\tau^{(n+2)}=(\tau^{(n+1)}, \underbrace{B^{(n+1)},\dots,B^{(n+1)}}_{{\frak l}_{B^{(n+1)}}}))\,,$$ 
we get that $\tau\prec \tau^{(n+2)}$ by using the inductive hypothesis $\tau'\prec \tau^{(n+1)}$.

Now we show that the points $(\tilde b-1,\tilde b)$ and $(\hat b-1,\hat b)$ belong to the set $\T^{A^{(n+2)}}_\bb$. The point $(\hat b-1,\hat b)$  belongs to $\T^{A^{(n+1)}}_\bb$ so $f^{A^{(n+1)}}\hat b\le \hat b$. If $\sigma$ is an intermediate block between between $A^{(n+1)}$ and $A^{(n+2)}$ then
$$f^\sigma(\hat b)=-(0,\tau,\overline{A^{(n+2)}})\le -(0,\overline{A^{(n+2)}})=\hat b$$
because $\tau$ is a tail block of $A^{(n+2)}$ obtained by eliminating $\sigma$, so
$A^{(n+2)}\prec \tau$.

Now we show that $f^\sigma(\tilde b)\ge \tilde b/(\tilde b+1)$. We have that $f^{\sigma^{(n+2)}}(\tilde b)=\tilde b/(\tilde b+1)$ by Lemma \ref{relbn}, and 
\[
f^{\sigma^{(n+2)}}(\tilde b)=-(0,\tau^{(n+1)},\overline{B^{(n+1)}})\, \quad f^\sigma(\tilde b)=-(0,\tau,\overline{B^{(n+1)}})
\] 
with $\tau$ being the end block of $A^{(n+2)}$ obtained by eliminating $\sigma$. But $\tau\prec\tau^{(n+2)}$ as we have just proved, hence $f^\sigma(\tilde b)\ge f^{\sigma(n+2)}(\tilde b)$. In conclusion, any intermediate sequence $\sigma$ between $A^{(n+1)}$ and $A^{(n+2)}$ satisfies
\[
\tilde b/(\tilde b+1)\le f^\sigma(\tilde b)\le f^\sigma(\hat b)\le \hat b\,,
\]
therefore the points $(\tilde b-1,\tilde b)$ and $(\hat b-1,\hat b)$ belong to the intermediate set $\T^\sigma_\bb$. 

We proved the induction step for $\T^{A^{(n+2)}}_\bb$, when $A^{(n+1)}$ is given by \eqref{case1}. A similar argument can be provided for the case when $A^{(n+1)}$ is given by \eqref{case2}, so the conclusion of the theorem is true.
%\foot{SK6: I am happy to trust you here!}
\end{proof}

We prove now that each set nonempty set $\E^{(n_i)}$ with $(n_i)$ not eventually aperiodic sequence is actually a singleton. 
\begin{thm}
Assume that $(n_i)$ is a not eventually periodic sequence such that the set $\E^{(n_i)}_b$ is nonempty. Then the set $\E^{(n_i)}_\bb$ is a point on the line segment $b-a=1$.
\end{thm}

\begin{proof}
The sequence $(n_i)$ satisfies the recursive relations \eqref{case1} or \eqref{case2}. We look at the set $\T^{A^{(n+1)}}_\bb$ and estimate the length of its lower base. In case \eqref{case1} %\foot{SK6: fixed ref} 
its upper vertex is given by \eqref{bne} and its lower base satisfies \eqref{ube1}. The lower base is a segment whose right end coordinate is
$$\underline a^r_{n+1}= -(0,A^{(n+1)},\overline{B^{(n+1)}})-1$$ and left end coordinate is
\[
\underline a^l_{n+1}=f^{A^{(n +1)}}(- (0,A^{(n+1)},\overline{B^{(n+1)}})) -1=-(0,\overline{B^{(n+1)}}) -1\,.
\]
Hence the length of the lower base is given by 
$$
L_{n+1}=\underline a^r_{n}-\underline a^l_{n+1}=(0,\overline{B^{(n+1)}})-(0,A^{(n+1)},\overline{B^{(n+1)}})\,.
$$
In case \eqref{case2}, the lower base is a segment whose right end coordinate is
$$
\underline a^r_{n+1}=-(0,A^{(n)},\overline{B^{(n)}})-1
$$  and the left end coordinate is given by 
\[
\underline a^l_{n+1}=f^{A^{(n +1)}}(- (0,A^{(n)},\overline{B^{(n)}})) -1=-(0,\overline{B^{(n)}}) -1\,.
\]
Hence the length of the lower base is given by
$$
L_{n+1}=\underline a^r_{n+1}-\underline a^l_{n+1}=(0,\overline{B^{(n)}})-(0,A^{(n)},\overline{B^{(n)}})\,.
$$
Notice that in the first case the two continued fraction expansions have in common at least the block $A^{(n)}$, while in the second case they have in common at least the block $A^{(n-1)}$. This implies that in both cases $L_{n+1}\rightarrow 0$ as $n\rightarrow \infty$. 
Moreover, the bases of the sets $\T^{n_1,\dots n_k}_\bb$ have non-increasing length and we have found a subsequence of these bases whose lengths converge to zero. Therefore the set $\E^{(n_i)}_\bb$ consists of only one point $(b-1,b)$, where
$b=-(0,n_1,n_2,\dots)$.
\end{proof}
The above result gives us a complete description of the set of exceptions $\E_\bb$ to the finiteness condition. It is a subset of the boundary segment $b=a+1$ of $\P$. Moreover, each set $\E^m_\bb$ is uncountable because the recursive construction of a nonempty set $\E^{(n_i)}_\bb$ allows for an arbitrary number of successive blocks $A^{(k)}$ at step $(k+1)$. Formally, one constructs a surjective map $\mathfrak j: \E^m_\bb\ra \mathbb N^{\mathbb N}$ by associating to a singleton set $\E^{(n_i)}_b$ a sequence of positive integers defined as%\foot{IU:added}

\medskip

\centerline{$\mathfrak j(\E^{(n_i)}_b)(k)=\#$ of consecutive $A^{(k)}$-blocks at the beginning of $(n_i)$.}

\medskip

The set $\E_\bb$ has one-dimensional Lebesgue measure $0$. The reason is that all associated formal continued fractions expansions of $b=-(0,n_1,n_2,\dots)$  have only two consecutive digits; such formal expansions $(0,n_1,n_2,\dots)$ are valid (-1,0)-continued fractions. Hence the set of such $b$'s has measure zero by Proposition \ref{bdigits1}. %\foot{IU:changed}
Analogous conclusions hold for $\E_\af$. Thus we have
%\foot{SK: We need to justify all the claims of the above theorem. The set has cardinality continuum since the set of sequences obtained by our substitution process is continuum -ref? How do we know it is of measure zero? Can we estimate the proportion that is taken out on each step?}
\begin{thm} \label{ex}
For any $(a,b)\in\P$, $b\neq a+1$, the finiteness condition holds. The set of exceptions $\E$ to the finiteness condition
is an uncountable set of one-dimensional Lebesgue measure $0$ that lies on the boundary $b=a+1$ of $\P$. 
\end{thm}
%\foot{SK16: ADDED}
Now we are able to provide the last ingredient in the proof of part (b) of the Main Result:
\begin{prop}\label{open-dense}
The strong cycle property is an open and dense condition. 
\end{prop}\begin{proof}
It follows from Theorems \ref{b-cycle} and \ref{a-cycle} that the condition is open. Theorem \ref{ex} asserts that for all $(a,b)\in\P$, $b\neq a+1$ the finiteness condition holds, i.e. all we need to show is that if $b$ has the week cycle property  or the $(a,b)$-expansions of $Sb$ and $T^{-1}b$  are eventually periodic, then in any neighborhood of it there is a $b$ with the strong cycle property. For,  if $b$ has the weak cycle property, it is a rational number obtained from the equation $\hat{f}^nT^mSb=0$, and any small perturbation of it will have the strong cycle property. Similarly, if the $(a,b)$-expansions of $Sb$ and $T^{-1}b$  are eventually periodic, then $b$ is a quadratic irrationality (see Remark \ref{quadratic}), and for any neighborhood of $b$ will contain values
%$\varepsilon>0$ there exists $b'$ such that $|b-b'|<\varepsilon$ 
satisfying the strong cycle property. A similar argument holds for $Sa$ and $Ta$.
\end{proof}

\section{Invariant measures and ergodic properties}\label{s:9}

Based on the finite rectangular geometric structure of the domain $D_{a,b}$ one can study the measure-theoretic properties of the Gauss-type map $\hat f_{a,b}:[a,b)\rightarrow [a,b)$,
\begin{equation}\label{1dGauss}
\hat f_{a,b}(x)=-\frac{1}{x}-\left\lfloor -\frac{1}{x}\right\rceil_{a,b}\,,\quad  \hat f_{a,b}(0)=0
\end{equation}
%where
%$$\lfloor x\rceil_{a,b}=\begin{cases}
%\lfloor x-a\rfloor &\text{ if }x<a\\
%0 & \text{ if } a\le x<b\\
%\lceil x-b\rceil & \text{ if } x\ge b
%\end{cases}
%$$
and its associated natural extension map $\hat F_{a,b}:\hat D_{a,b}\ra \hat D_{a,b}$
\begin{equation}\label{2dGauss}
\hat F _{a,b}=\left(\hat f_{a,b}(x),-\frac{1}{y-\lfloor -1/x\rceil_{a,b}}\right).
\end{equation}
%\foot{SK6: added period at the end of the formula (9.2)}
We remark that $\hat F_{a,b}$ is obtained from the map $F_{a,b}$ induced on the set $D_{a,b}\cap \{(x,y)| a\le y <b\}$ by a change of coordinates $x'=y$, $y'=-1/x$. Therefore the domain
$\hat D_{a,b}$ is easily identified knowing $D_{a,b}$ and may be considered its ``compactification".
%\foot{SK6: added to the end of the sentence}

We present the simple case when $\displaystyle 1\le -\frac{1}{a}\le b+1$ and $a-1\le -\frac{1}{b}\le -1$. The general 
theory  is the subject  our paper in preparation \cite{KU6}.

The truncated orbits of $a$ and $b$ are
\begin{eqnarray*}
& &\mathcal{L}_a =\left\{a+1,-\frac{1}{a+1}\right\}, \quad \mathcal{U}_a=\left\{-\frac{1}{a},-\frac{a+1}{a}\right\} \\
& &\mathcal{L}_b =\left\{-\frac{1}{b},\frac{b-1}{b}\right\}, \quad \mathcal{U}_b=\left\{b-1,-\frac{1}{b-1}\right\}
\end{eqnarray*}
and the end points of the cycles are $c_a=\frac{a}{a+1}$, $c_b=\frac{b}{1-b}$.
\begin{thm} If $1\le -\frac{1}{a}\le b+1$ and $a-1\le -\frac{1}{b}\le -1$, then the domain $\hat D_{a,b}$ of $\hat F_{a,b}$ is given by
\[
\begin{split}
\hat D_{a,b}&=[a,-\frac{1}{b}+1]\times[-1,0]\cup [-\frac{1}{b}+1,a+1]\times[-1/2,0]\\
&\quad \cup [b-1,-\frac{1}{a}-1]\times [0,1/2]\cup [-\frac{1}{a}-1,b]\times [0,1]
\end{split}
\]
and $\hat F_{a,b}$ preserves the Lebesgue equivalent probability measure
\begin{equation}\label{dnu}
d\nu_{a,b}=\frac{1}{\log[(1+b)(1-a)]}\frac{dxdy}{(1+xy)^2}\,.
\end{equation}
%\item[(c)] $G_{a,b}$ is ergodic with respect to $\nu$.
\end{thm}

\begin{proof}
The description of $\hat D_{a,b}$ follows directly from the cycle relations and the finite rectangular structure. It is a standard computation that the measure $\frac{dxdy}{(1+xy)^2}$ is preserved by $\hat F_{a,b}$, by using the fact any M\"obius transformation, hence $F_{a,b}$, preserves the measure $\frac{du\,dw}{(w-u)^2}$, and $\hat F_{a,b}$ is obtained from $F_{a,b}$ by coordinate changes $x=w, y=-1/u$.  
%\foot{IU: added explanation}

%{SK6: I suggest to explain this a little better as in my Pisa lectures pp.83-84 and give a reference to Adler-Flatto for use of these coordinates. Please, make changes.}
Moreover, the density $\frac{1}{(1+xy)^2}$ is bounded away from zero on $\hat D_{a,b}$ and
$$
\int_{\hat D_{a,b}}\frac{dxdy}{(1+xy)^2}=\log[(b+1)(1-a)]<\infty
$$
hence the last part of the theorem is true.
\end{proof}
\begin{figure}[htb]
\includegraphics[scale=0.8]{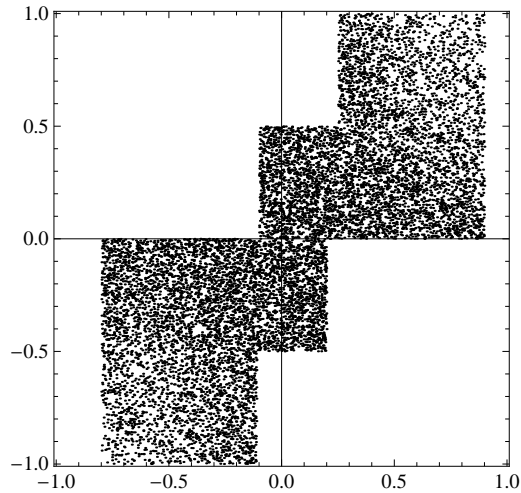}
\caption{Typical domain $\hat D_{a,b}$ for the case studied}
%\center{The map $F$ is injective}
%\label{g8}
\end{figure}

The Gauss-type map $\hat f_{a,b}$ is a factor of $\hat F_{a,b}$ 
(projecting on the $x$-coordinate), so one can obtain its smooth 
invariant measure $d\mu_{a,b}$ by integrating $d\nu_{a,b}$ over $\hat 
D_{a,b}$ with respect to the $y$-coordinate as explained in 
\cite{AF3}. Thus, if we know the exact shape of the set $D_{a,b}$, we  
can calculate the invariant measure precisely.

The measure $d\mu_{a,b}$ is ergodic and the measure-theoretic entropy 
of $\hat f_{a,b}$ can be computed explicitly using Rokhlin's formula.

%The Gauss-type map $\hat f_{a,b}$ is a factor of $\hat F_{a,b}$ (projecting on the $x$-coordinate) so one can obtain its smooth invariant measure $d\mu_{a,b}$ by integrating $d\nu_{a,b}$ over $\hat D_{a,b}$ with respect to the $y$-coordinate. The measure $d\mu_{a,b}$ is ergodic and the measure-theoretic entropy of $\hat f_{a,b}$ can be computed explicitly using Rokhlin's formula.
\begin{thm} $ $
The map $\hat f_{a,b}:[a,b)\ra [a,b)$ is ergodic with respect to Lebesgue equivalent invariant probability measure
\begin{equation}
d\mu_{a,b}=\frac{1}{C_{a,b}}\left(\frac{\chi_{(a,-\frac{1}{b}+1)}}{1-x}+\frac{\chi_{(-\frac{1}{b}+1,a+1)}}{2-x}+\frac{\chi_{(b-1,-\frac{1}{a}-1)}}{x+2}+
\frac{\chi_{(-\frac{1}{a}-1,b)}}{x+1}\right)dx
\end{equation}
where $C_{a,b}=\log[(1+b)(1-a)]$. The measure-theoretic entropy of $\hat f_{a,b}$ is given by
\begin{equation}\label{entropy}
h_{\mu_{a,b}}(\hat f_{a,b})=\frac{\pi^2}{3\log[(1-a)(1+b)]}\,.
\end{equation}
\end{thm}
\begin{proof}
The measure $d\mu_{a,b}$ is obtained by integrating $d\nu_{a,b}$ over $\hat D_{a,b}$. Ergodicity follows from a more general result concerning one-dimensional expanding maps (see \cite{AF3,Zw}). To compute the entropy, we use Rokhlin's formula 
%$h_{\mu_{a,b}}(\hat f_{a,b})=\int_a^b\log|\hat f'_{a,b}|d\mu_{a,b}=-2\int_a^b\log|x|d\mu_{a,b}$, so
\[
\begin{split}
h_{\mu_{a,b}}(\hat f_{a,b})&=\int_a^b\log|\hat f'_{a,b}|d\mu_{a,b}=-2\int_a^b\log|x|d\mu_{a,b}\\
&=\frac{-2}{C_{a,b}}\left(\int_a^{-\frac{1}{b}+1}\frac{\log|x|}{1-x}dx+\int_{-\frac{1}{b}+1}^{a+1}\frac{\log|x|}{2-x}dx\right.\\
&\qquad+\left.\int_{b-1}^{-\frac{1}{a}-1}\frac{\log|x|}{x+2}dx+\int_{-\frac{1}{a}-1}^b\frac{\log|x|}{x+1}dx \right)
\end{split}
\]

Let $I(a,b)$ denote the sum of the four integrals. The function depends smoothly on $a,b$, hence we can compute the partial derivatives $\partial I/\partial a$ and $\partial I/\partial b$. We get that both partial derivatives are zero, hence $I(a,b)$ is constant. Using $a=-1, b=1$, we get
$$I(a,b)=I(-1,1)=2\int_0^1\frac{\log|x|}{1+x}dx=-\pi^2/6\,,$$
and the entropy formula \eqref{entropy}.
\end{proof}

%%%%%%

\end{document}